\documentclass[11pt]{article}
\usepackage{amsmath, amssymb, amscd, epsfig}

\newtheorem{thm}{Theorem}
\newtheorem{prop}[thm]{Proposition}
\newtheorem{lem}[thm]{Lemma}

\newtheorem{cor}[thm]{Corollary}
\newtheorem{rem}[thm]{Remark}
\newtheorem{df}[thm]{Definition}
\newtheorem{ex}[thm]{Example}

\renewcommand{\epsilon}{\varepsilon}

\newcommand{\BB}{\mathbb}
\newcommand{\g}{\mathfrak}
\newcommand{\pf}{\noindent {\it Proof. }}
\newcommand{\qed}{\nopagebreak $\qquad$ $\square$ \vskip5pt}
\newcommand{\separate}{\vskip5pt}

\newcommand{\supp}{\operatorname{supp}}

\newcommand{\tr}{\operatorname{Tr}}
\newcommand{\End}{\operatorname{End}}
\newcommand{\lefttorightarrow}{\circlearrowright}

\begin{document}

\title{\bf Lecture Notes on Equivariant Cohomology}
\author{Matvei Libine}
\maketitle

\begin{section}
{Introduction}
\end{section}

These are the lecture notes for the introductory graduate course I taught
at Yale during Spring 2007. I mostly followed \cite{GS, BGV, AB, Par2},
and there are no original results in these notes.

Let $G$ be a compact Lie group acting on a topological space $M$.
For the topologists, the equivariant cohomology of $M$ is defined to be
the ordinary cohomology of the space $(M \times E) /G$,
where $E$ is any contractible topological space on which $G$ acts freely.
(This definition does not depend on the choice of $E$.)
If $M$ is a finite-dimensional manifold there is an alternative way of
defining the equivariant cohomology groups of $M$ involving de Rham theory.

The notion of equivariant cohomology plays an important role in
symplectic geometry, algebraic geometry, representation theory and
other areas of mathematics.

The central topics in this course are the localization theorems
due to Atiyah-Bott and Berline-Vergne and their consequences such as the
Kirillov's character formula for the characters of compact groups.
In particular, we will prove the Berline-Vergne localization formula
(which can be regarded as a generalization of the Duistermaat-Heckman Theorem)
expressing integrals of equivariant forms as sums over fixed points.

Other topics covered in these notes are the twisted de Rham complex,
the equivariant vector bundles and the equivariant characteristic classes,
the equivariant Thom class, equivariant formality,
Goresky-Kottwitz-MacPherson Theorem (GKM Theory), Paradan-Witten localization.


\separate

\begin{section}
{Equivariant Cohomology via Algebraic Topology}
\end{section}

In this section we give a topological definition of equivariant cohomology;
we essentially reproduce Chapter 1 of \cite{GS}.
Let $G$ be a compact Lie group acting continuously on a topological space $M$.
We write $G \lefttorightarrow M$ to denote this action.

\separate

\subsection{Basic Definitions}

\begin{df}
The action $G \lefttorightarrow M$ is {\em free} if for all $x \in M$
and $g \in G$ $g \cdot x = x$ implies $g=e$ (the identity element in $G$).
\end{df}

If $G$ acts on $M$ freely, then the quotient space $M/G$ is usually as nice
a topological space as $M$ is. Say, if $M$ is a manifold, then so is $M/G$.
heuristically, if $G \lefttorightarrow M$ freely, then the equivariant
cohomology $H^*_G(M)$ should be $H^*(M/G)$.
For example, if $G$ acts on itself by left translations, then
$H^*_G(G) = H^*(pt)$.
On the other hand, if the action $G \lefttorightarrow M$ is {\em not} free,
then the space $M/G$ may be pathological and $H^*_G(M)$ provides the
``right'' substitute for $H^*(M/G)$.

\begin{ex}  \label{Bott}
Let $G$ be the circle group $S^1$ act on the sphere $M=S^2$ by rotation about
the vertical axis. Since the poles of $S^2$ remain fixed by all elements of
$S^1$, this action is not free. The quotient $S^2/S^1$ can be identified with
a closed segment, which is contractible.
We will see later that $H^*_{S^1}(S^2) \ne H^*(pt)$.
\end{ex}

Recall that if $M$ is homotopy equivalent to another topological space
$\tilde M$, then $H^*(M) = H^*(\tilde M)$.
Starting with an action $G \lefttorightarrow M$ which may not be free
we expect to find a homotopy equivalent space $\tilde M$
on which $G$ acts freely and set $H^*_G(M) = H^*(\tilde M/G)$.
One way to find such a space $\tilde M$ is to take $\tilde M = M \times E$,
where $E$ is a contractible space on which $G$ acts freely.

\begin{df}
The equivariant cohomology of $M$ is defined to be
$H^*_G(M) = H^*((M \times E)/G)$.
\end{df}

\begin{rem}
Some people advocate the notation $H^*(G \lefttorightarrow M)$ instead of
$H^*_G(M)$ to emphasize the dependence of equivariant cohomology on the way
$G$ acts on $M$.
\end{rem}

For example, $H^*_G(pt) = H^*(E/G)$.

We still need to verify that such space $E$ exists and that
$H^*_G(M)$ is independent of the choice of the space $E$.

Note that if $G$ acts on $M$ freely, then the projection
$M \times E \twoheadrightarrow M$ induces a map
$(M \times E)/G \twoheadrightarrow  M/G$.
If $M$ is a manifold, then this is a fibration with contractible fiber $E$.
So, in this case, $H^*_G(M) = H^*((M \times E)/G) = H^*(M/G)$.

\begin{ex}
Take $G = S^1$ (the circle group) realized as
$S^1 = \{ z \in \BB C ;\: |z|=1 \}$.
And let $S^1$ act on the unit sphere $S^{2n+1}$ sitting inside $\BB C^{n+1}$
by $z \cdot (w_0,\dots,w_n) = (zw_0,\dots,zw_n)$.
The action is free and $S^{2n+1}/S^1 = \BB CP^n$,
but $S^{2n+1}$ is not contractible.
However, $S^{\infty}$ -- the unit sphere in $\BB C^{\infty}$ is contractible,
so for $G=S^1$ we may take $E= S^{\infty}$, and
$S^{\infty} \twoheadrightarrow S^{\infty}/S^1 = \BB CP^{\infty}$
is the classifying (or universal) bundle for $S^1$.
\end{ex}

\separate

\subsection{Classifying Spaces}

Fix a Lie group $G$ and assume that it acts on a contractible space $E$ freely.
We can form a principal $G$-bundle $E \twoheadrightarrow B=E/G$.

\begin{df}
If $Y$ is a topological space with a map $f: Y \to B$, then we can form the
{\em pull-back bundle}
$$
\begin{CD}
f^*E   @>>>   E  \\
@VVV      @VV{p}V  \\
Y   @>{f}>>   B
\end{CD}
$$
by setting $f^*E = \{(y,e) \in Y \times E ;\: f(y)=p(e) \}$,
then the map $(y,e) \mapsto y$ makes $f^*Y \to Y$ into a principal
$G$-bundle over $Y$.
\end{df}

The following result is proved in most introductory courses on algebraic
topology.

\begin{thm} [Classification Theorem]
Let $\pi: X \to Y$ be a principal $G$-bundle. Then there are a map
$f: Y \to B$ and an isomorphism of principal $G$-bundles $\Phi: X \to f^*E$.
Moreover, $f$ and $\Phi$ are unique up to homotopy.
\end{thm}

\begin{rem}
In other words we have a natural bijection
$$
\Bigl\{ \begin{matrix} \text{equivalence classes of}  \\
\text{principal $G$-bundles}\end{matrix}\Bigr\}
\longleftrightarrow
\{f:Y \to B\} / \text{homotopy}.
$$
For this reason the bundle $E \twoheadrightarrow B$ is called the
{\em classifying bundle}.
\end{rem}

The first consequence of the classification theorem is that the classifying
bundle is unique:

\begin{thm}
If $E_1$ and $E_2$ are contractible spaces on which $G$ acts freely,
then they are equivalent as $G$-spaces. This means that there exist
$G$-equivariant maps $\phi: E_1 \to E_2$, $\psi: E_2 \to E_1$
and $G$-equivariant homotopies $\psi \circ \phi \simeq Id_{E_1}$,
$\phi \circ \psi \simeq Id_{E_2}$.
\end{thm}

\begin{cor}
The definition of equivariant cohomology is independent of the choice of
the space $E$, i.e.
$$
H^*((X \times E_1)/G) \simeq H^*((X \times E_2)/G).
$$
\end{cor}

We make the following observation. If $H \subset G$ is a subgroup,
then the space $E$ which ``works'' for $G$ also ``works'' for $H$.
We have a map $(X \times E)/H \twoheadrightarrow (X \times E)/G$
with fiber $G/H$. Hence we obtain a canonical map
$H^*_G(X) \to H^*_H(X)$:
\begin{equation}  \label{restriction_map}
H^*_G(X) = H^*((X \times E)/G) \to H^*((X \times E)/H) = H^*_H(X).
\end{equation}

Every compact group has a faithful finite-dimensional representation
$(\pi,V)$, i.e. the map $\pi: G \to Aut(V)$ is injective.
Hence $G$ can be realized as a Lie subgroup of $U(n)$,
for some $n \in \BB N$.

We still need to show that there exists a contractible space $E$ on which $G$.
For this purpose we assume that $G$ is a Lie subgroup of some $U(n)$.
Here we reproduce the construction of the classifying space $E$
given in Section 1.2 in \cite{GS}.

Let $L^2[0,\infty)$ denote the space of square integrable $\BB C$-valued
functions on the positive real numbers relatively to the standard Lebesgue
measure. This space comes equipped with an inner product which gives it
the norm topology and makes it a Hilbert space.
Consider the space of orthonormal $n$-frames
$$
E = \{ {\bf f} = (f_1, \dots,f_n) \in
L^2[0,\infty) \times \dots \times L^2[0,\infty)
;\: \langle f_i, f_j \rangle = \delta_{ij}\}
$$
which is an open subset of $L^2[0,\infty) \times \dots \times L^2[0,\infty)$
equipped with the product topology.
We define the action of $G$ on $E$ by
$$
g \cdot (f_1, \dots,f_n) = (\tilde f_1, \dots, \tilde f_n),
\qquad
\tilde f_i = \sum_{j=1}^n a_{ij} f_j,
$$
where $g \in G \subset U(n)$ is represented by an invertible
$n \times n$ matrix $(a_{ij})$.
Clearly, $G$ acts on $E$ continuously and freely.

Let $E' \subset E$ denote the subset of $n$-tuples of functions which all
vanish on the interval $[0,1]$.

\begin{lem}
The subset $E'$ is a deformation retract of $E$.
\end{lem}

\pf
For any $f \in L^2[0,\infty)$ define $T_tf$ by
$$
T_tf(x) =
\begin{cases}
0 & \text{for } 0 \le x < t;  \\
f(x-t) & \text{for } t \le x < \infty.
\end{cases}
$$
Define
$$
{\bf T}_t {\bf f} = (T_tf_1, \dots,T_tf_n) \qquad \text{for }
{\bf f} = (f_1, \dots,f_n) \in L^2[0,\infty) \times \dots \times L^2[0,\infty).
$$
Since $T_t$ preserves the inner product on $ L^2[0,\infty)$,
we see that ${\bf T}_t$ is a deformation retract of $E$ into $E'$.
\qed

\begin{prop}
The space $E$ is contractible and the projection
$E \twoheadrightarrow E/G$ is a (locally trivial) fiber bundle.
\end{prop}

\pf
To prove that $E$ is contractible it is sufficient to show that
$E'$ is contractible to a point within $E$.
Pick an $n$-frame ${\bf g} = (g_1, \dots, g_n) \in E$ such that all its
components are supported in $[0,1]$.
For each ${\bf f} = (f_1, \dots , f_n) \in E'$ we define
$$
{\bf R}_t {\bf f} = \bigl( tg_1 + (1-t)f_1, \dots tg_n + (1-t)f_n \bigr),
\qquad t \in [0,1].
$$
Clearly, ${\bf R}_t (E') \subset E$ for all $t$
and ${\bf R}_t$ is a continuous deformation of $E'$ to ${\bf g}$ within $E$.

It remains to show that $p: E \twoheadrightarrow E/G$ is a fiber bundle.
Suppose first that $G=U(n)$.
Pick a Hilbert space orthonormal basis $\{v_1,v_2,\dots\}$ of $L^2[0,\infty)$.
For each sequence of $n$ different integers
$1 \le i_1 < i_2 < \dots < i_n < \infty$ we set
$\BB C^n_{\{i_1,\dots,i_n\}} = \BB C$-span of $v_{i_1}, \dots, v_{i_n}$
in $L^2[0,\infty)$, and define
$$
U_{\{i_1,\dots,i_n\}} = \biggl\{ {\bf f} = (f_1, \dots,f_n) \in E ;\:
\begin{matrix}
\text{the orthogonal projections of $f_1,\dots,f_n$}  \\
\text{onto $\BB C^n_{\{i_1,\dots,i_n\}}$ form a basis in
$\BB C^n_{\{i_1,\dots,i_n\}}$}
\end{matrix} \biggr\}.
$$
The sets $U_{\{i_1,\dots,i_n\}}$ are $G$-invariant and form an open
covering of $E$. Hence their projections
$\{p(U_{\{i_1,\dots,i_n\}}) \}_{1 \le i_1 < i_2 < \dots < i_n < \infty}$
form an open covering of $E/G$.
Note that every element ${\bf f} \in U_{\{i_1,\dots,i_n\}}$ can be uniquely
written as ${\bf f} = g_{\bf f} \cdot {\bf f'}$, where $g_{\bf f} \in U(n)$
and ${\bf f'} = (f_1',\dots,f_n') \in E$ is such that the orthogonal
projections of $f_1',\dots,f_n'$ onto $\BB C^n_{\{i_1,\dots,i_n\}}$ are 
$v_{i_1}, \dots, v_{i_n}$ respectively.
The maps
$\varphi_{\{i_1,\dots,i_n\}}: U_{\{i_1,\dots,i_n\}} \to
G \times p(U_{\{i_1,\dots,i_n\}})$ defined by
$$
\varphi_{\{i_1,\dots,i_n\}} : {\bf f} \mapsto (g_{\bf f}, p({\bf f})),
$$
provide a trivialization of $E \twoheadrightarrow E/G$ over each
$p(U_{\{i_1,\dots,i_n\}})$.

Finally, if $G$ is a proper subgroup of $U(n)$, then 
$p: E \twoheadrightarrow E/G$ is a pull-back of $E \twoheadrightarrow E/U(n)$:
$$
\begin{matrix}
E & \longrightarrow & E  \\
\downarrow & \quad & \downarrow  \\
E/G & \longrightarrow & E/U(n)
\end{matrix}
$$
hence a locally trivial fiber bundle too.
\qed

\separate

\begin{section}
{Equivariant Cohomology via Differential Forms}
\end{section}

\subsection{Equivariant Differential Forms}

Let $M$ be a smooth manifold. Since we will integrate over $M$,
let us also suppose that $M$ is oriented and compact.
Let $G$ be a compact Lie group acting on $M$ smoothly.
We denote by $\Omega^*(M)$ the algebra of smooth
$\BB C$-valued differential forms on $M$.

If $N$ is another smooth manifold and $f:N \to M$ is a
${\cal C}^{\infty}$ map, we get a pull-back map
$f^*: \Omega^*(M) \to \Omega^*(N)$.
In particular, given an element $g \in G$, we get the pull-back map
induced by the map $M \to M$, $x \mapsto g^{-1} \cdot x$.
By abuse of notation, we write $\alpha \mapsto g \cdot \alpha$,
$\alpha \in \Omega^*(M)$.
This way we get an action of $G$ on $\Omega^*(M)$.
The smooth $\BB C$-valued functions on $M$ can be regarded as $0$-forms
and the group $G$ acts on those by
$$
(g \cdot f) (x) = f(g^{-1} \cdot x),
\qquad g \in G, \: f \in \Omega^0(M), \: x \in M.
$$
Note that if we omit the inverse in $g^{-1} \cdot x$
we do not get a left action.

Let $\mathfrak{g}$ denote the Lie algebra of $G$. We have an adjoint action
of $G$ on $\mathfrak{g}$. If $G$ is a Lie subgroup of $GL(n)$, this action is
$$
Ad(g): X \mapsto Ad(g)X = gXg^{-1},
\qquad X \in \mathfrak{g}, \: g \in G.
$$
In general, the Lie algebra $\mathfrak{g}$ is identified with the tangent
space at the identity element of $G$: $\mathfrak{g} \simeq T_eG$, and the
adjoint map $Ad(g)$ is the tangent map $T_eG \to T_eG$ induced by the
conjugation map $G \to G$, $h \mapsto g h g^{-1}$.
Dually, we get a coadjoint action of $G$ on $\mathfrak{g}^*$ --
the dual of $\mathfrak{g}$ -- denoted by $Ad^{\#}$:
$$
\langle Ad^{\#}(g)l, Ad(g)X \rangle = \langle l, X \rangle,
\qquad l \in \mathfrak{g}^*, \: X \in \mathfrak{g}.
$$

\begin{df}
A ${\cal C}^{\infty}$ $G$-equivariant differential form on $M$ is a
${\cal C}^{\infty}$ map $\alpha: \mathfrak{g} \to \Omega^*(M)$ which is
$G$-equivariant, i.e. the following diagram commutes:
$$
\begin{CD}
\mathfrak{g}   @>{\alpha}>>   \Omega^*(M)  \\
@V{Ad(g)}VV      @VV{g}V  \\
\mathfrak{g}   @>{\alpha}>>   \Omega^*(M)
\end{CD}
\qquad \forall g \in G.
$$
In other words,
\begin{equation}  \label{eq-def}
\alpha(X) = g^{-1} \cdot \alpha (Ad(g)X),
\qquad \forall X \in \mathfrak{g}, \: g \in G.
\end{equation}
\end{df}

Note that the wedge product of two equivariant forms also is an equivariant
form. Thus ${\cal C}^{\infty}$ equivariant forms form an algebra denoted by
$\Omega^{\infty, *}_G(M)$. The equivariant forms of our interest usually
are not the ``homogeneous'' maps $\mathfrak{g} \to \Omega^*(M)$.
We define the $\BB Z_2$-grading on $\Omega^{\infty, *}_G(M)$ by setting
$$
\Omega^{even}(M) = \bigoplus_{\text{$n$-even}} \Omega^n(M),
\qquad
\Omega^{odd}(M) = \bigoplus_{\text{$n$-odd}} \Omega^n(M)
$$
and
\begin{align*}
\Omega^{\infty,even}_G(M) &= \{ \alpha \in \Omega^{\infty,*}_G(M) ;\:
\alpha: \mathfrak{g} \to \Omega^{even}(M) \},  \\
\Omega^{\infty,odd}_G(M) &= \{ \alpha \in \Omega^{\infty,*}_G(M) ;\:
\alpha: \mathfrak{g} \to \Omega^{odd}(M) \}.
\end{align*}

\begin{rem}
If the group $G$ is abelian, i.e. it is the circle group $S^1$ or the
torus $T^n = S^1 \times \dots \times S^1$ ($n$ times), then the adjoint
action is trivial, so (\ref{eq-def}) becomes
$$
\alpha(X) = g \cdot \alpha(X)
$$
and the $G$-equivariant form really is a $G$-invariant form.
\end{rem}

\begin{rem}
If $M_1$ and $M_2$ are two spaces with $G$-action, then one can define
the action of $G$ on $\operatorname{Map}(M_1,M_2)$ by
$$
(g \cdot f)(x) = g \cdot (f (g^{-1} \cdot x)),
\qquad g \in G, \: f \in \operatorname{Map}(M_1,M_2), \: x \in M_1.
$$
Thus we obtain an action of $G$ on
$\operatorname{Map}(\mathfrak{g},\Omega^*(M))$:
$$
(g \cdot \alpha)(X) = g \cdot ( \alpha(Ad(g^{-1})X)).
$$
Switching $g$ and $g^{-1}$ we see that $\alpha \in \Omega^{\infty, *}_G(M)$
if and only if $\alpha \in \operatorname{Map}(\mathfrak{g},\Omega^*(M))$,
$\alpha$ is smooth and $G$-invariant.
\end{rem}

\separate

\subsection{The Twisted de Rham Differential}

Recall that the group $G$ acts on $M$ which in turn induces action on
the space of functions $\Omega^1(M)$
$$
(g \cdot f)(x) = f(g^{-1} \cdot x).
$$
Differentiating this action we obtain an action by the Lie algebra
$\mathfrak{g}$. Thus, for each element $X \in \mathfrak{g}$ we obtain
a vector field denoted by $L_X$:
$$
(L_X f)(x) =_{def} \frac d{dt} f(\exp(-tX) \cdot x) \Bigr|_{t=0}.
$$
By construction, these vector fields satisfy $[L_X,L_Y] = L_{[X,Y]}$
and we have a map of Lie algebras
$\mathfrak{g} \to \{ \text{smooth vector fields on $M$} \}$,
$X \mapsto L_X$.
Note that if we did not have the negative sign in the definition of $L_X$
then $[L_X,L_Y] = L_{[X,Y]}$ would not hold.

In the same fashion, the action $G \lefttorightarrow M$ induces (left) actions
$G \lefttorightarrow \Omega^*(M)$ and 
$G \lefttorightarrow \{ \text{smooth vector fields on $M$} \}$.
Hence we get infinitesimal actions of $\mathfrak{g}$ on differential forms
and vector fields which we still denote by $L_X$.
Note that if $F$ is a vector field on $M$, then $X \in \mathfrak{g}$ acts on
it by sending it to the Lie bracket of vector fields $[L_X, F]$.

Next we recall contractions by vector fields. If $F$ is a vector field on $M$
and $\omega \in \Omega^1(M)$, get a function $\iota(F)\omega$ on $M$
such that its value at each point $x \in M$ is given by evaluation of
$\omega_x \in T^*_xM$ at $F_x \in T_xM$.
The contraction map uniquely extends to a map
$\iota(F): \Omega^*(M) \to \Omega^*(M)$ lowering degrees by $1$
so that it sends $\Omega^0(M)$ into $0$ and
$$
\iota(F) (\omega_1 \wedge \omega_2) =
( \iota(F) \omega_1) \wedge \omega_2
+ (-1)^k \omega_1 \wedge (\iota(F) \omega_2),
\qquad \forall \omega_1 \in \Omega^k(M), \: \omega_2 \in \Omega^*(M).
$$
It follows that $(\iota_X)^2=0$.

For $X \in \mathfrak{g}$ we denote by $\iota_X$ the contraction
$\iota(L_X) : \Omega^*(M) \to \Omega^*(M)$.

\begin{df}
The twisted de Rham differential
$d_{eq} : \Omega^{\infty, *}_G(M) \to \Omega^{\infty, *}_G(M)$
is defined as follows. Given an $\alpha \in \Omega^{\infty, *}_G(M)$,
$d_{eq} \alpha$ is a map $\mathfrak{g} \to \Omega^*(M)$ given by
$$
(d_{eq}(\alpha)(X) = d(\alpha(X)) - \iota_X(\alpha(X)),
$$
where $d$ denotes the ordinary de Rham differential.
\end{df}

\begin{lem}
\begin{enumerate}
\item
$d_{eq} \alpha \in \Omega^{\infty, *}_G(M)$;

\item
$(d_{eq})^2=0$.
\end{enumerate}
\end{lem}

\pf
To prove the first part we need to show that
$d_{eq} \alpha : \mathfrak{g} \to \Omega^*(M)$ is $G$-equivariant.
Since the de Rham differential commutes with pull-backs,
the first component $d(\alpha(X))$ is $G$-equivariant.
To prove that $\iota_X(\alpha(X))$ is $G$-equivariant we use the identity
$$
\iota_{Ad(g)X} = g \cdot \iota_X \cdot g^{-1}.
$$
So,
$$
\iota_{Ad(g)X} ( \alpha (Ad(g)X)) =
(g \cdot \iota_X \cdot g^{-1}) \cdot (g \cdot \alpha(X))
= g \cdot (\iota_X \alpha(X)),
$$
proving $d_{eq} \alpha$ is $G$-equivariant.

Since $d^2=0$, $(\iota_X)^2=0$ and
$$(d \iota_X + \iota_X d) \omega = L_X \omega,
\qquad \forall \omega \in \Omega^*(M),
$$
we have:
\begin{multline*}
((d_{eq})^2 \alpha)(X)
= d^2(\alpha(X)) - ( d \iota_X + \iota_X d) \alpha(X)
+ (\iota_X)^2 \alpha(X)  \\
= -L_X \alpha(X)
= - \frac d{dt} ( \exp(tX) \cdot \alpha(X) ) \Bigr|_{t=0}
= - \frac d{dt} \alpha ( Ad(\exp(-tX)) X ) \Bigr|_{t=0}  \\
= - \frac d{dt} \alpha ( e^{-ad(tX)} X ) \Bigr|_{t=0}
= - \frac d{dt} \alpha (X) \Bigr|_{t=0} =0.
\end{multline*}
This proves $(d_{eq})^2=0$.
\qed

\separate

\subsection{Equivariant Cohomology}

Recall that $d: \Omega^*(M) \to \Omega^{*+1}(M)$ and
$\iota_X : \Omega^*(M) \to \Omega^{*-1}(M)$, hence
$d_{eq}: \Omega^{\infty,even}_G(M) \leftrightarrows \Omega^{\infty,odd}_G(M)$.

\begin{df}
The equivariant cohomology with ${\cal C}^{\infty}$ coefficients is the
cohomology of $(\Omega^{\infty,*}_G(M), d_{eq})$ with $\BB Z_2$-grading.
That is,
\begin{align*}
H^{\infty, even}_G(M) &= H^{even}(\Omega^{\infty,*}_G(M), d_{eq}) =
\{ \alpha \in \Omega^{\infty,even}_G(M) ;\: d_{eq}\alpha=0 \}/
\{ d_{eq}\beta ;\: \beta \in \Omega^{\infty,odd}_G(M) \},   \\
H^{\infty,odd}_G(M) &= H^{odd}(\Omega^{\infty,*}_G(M), d_{eq}) =
\{ \alpha \in \Omega^{\infty,odd}_G(M) ;\: d_{eq}\alpha=0 \}/
\{ d_{eq}\beta ;\: \beta \in \Omega^{\infty,even}_G(M) \}.
\end{align*}
\end{df}

We are really interested in a complex computing $H^*_G(M)$.
Consider a subalgebra of $\Omega^{\infty,*}_G(M)$ consisting of
{\em polynomial} maps $\alpha: \mathfrak{g} \to \Omega^*(M)$, i.e.
$(S(\mathfrak{g}^*) \otimes \Omega^*(M))^G$, where the superscript
$G$ means we are taking the $G$-invariant elements.
In some sense, $\Omega^{\infty,*}_G(M)$ is the closure of
$(S(\mathfrak{g}^*) \otimes \Omega^*(M))^G$.
We will use the notation $\BB C[\mathfrak{g}]$ for $S(\mathfrak{g}^*)$
-- the algebra of polynomials on $\mathfrak{g}$.

\begin{df}
The degree of a polynomial map $\alpha: \mathfrak{g} \to \Omega^*(M)$
or an element of $\BB C[\mathfrak{g}] \otimes \Omega^*(M)$ is
$$
(\text{differential form degree}) + 2(\text{polynomial degree}).
$$
\end{df}

\begin{lem}
The twisted de Rham differential increases this degree by $1$.
\end{lem}

\pf
Recall that $d_{eq} = d - \iota_X$. The first component $d$ increases the
differential form degree by one and leaves the polynomial degree unchanged.
And the second component $\iota_X$ lowers the differential form degree by $1$
and at the same time increases the polynomial degree by $1$,
so the total degree is increased by $1$.
\qed

From now on the $\BB Z$-graded complex
$((\BB C[\mathfrak{g}] \otimes \Omega^*(M))^G, d_{eq})$ will be denoted by
$\Omega^*_G(M)$ and called the {\em twisted de Rham complex}.

\begin{thm}[Cartan, 1950]
The twisted de Rham complex $\Omega^*_G(M)$ computes the equivariant cohomology
of $M$, i.e.
\begin{equation}  \label{Cartan_thm}
H^*_G(M, \BB C) \simeq H^*(\Omega^*_G(M), d_{eq}).
\end{equation}
The assumptions are: $G$ is a compact Lie group acting on
a smooth compact manifold $M$.
\end{thm}

Let $\alpha: \mathfrak{g} \to \Omega^*(M)$ be a ${\cal C}^{\infty}$
equivariant form. As usual, $X \in \mathfrak{g}$ and let $\alpha(X)_{[j]}$
denote the component of $\alpha(X)$ lying in $\Omega^j(M)$.
We have the following diagram:
$$
\begin{matrix}
\begin{matrix}
\alpha(X) &= & \alpha(X)_{[0]} & + & \alpha(X)_{[1]} & + \\
\: & & & \searrow \hspace{-.16in} \swarrow
& & \searrow \hspace{-.16in} \swarrow &  \\
(d_{eq} \alpha)(X) & = & \beta(X)_{[0]} & + & \beta(X)_{[1]} & +
\end{matrix}
&
\begin{matrix}
\dots & + & \alpha(X)_{[n-1]} & + & \alpha(X)_{[n]} & &  \\
\: & \searrow \hspace{-.16in} \swarrow
& & \searrow \hspace{-.16in} \swarrow & & \searrow  & \\
\dots & + & \beta(X)_{[n-1]} & + & \beta(X)_{[n]} & + & \beta(X)_{[n+1]},
\end{matrix}
\end{matrix}
$$
where
\begin{align*}
\beta(X)_{[0]} &= -\iota_X \alpha(X)_{[1]}  \\
\beta(X)_{[1]} &= d\alpha(X)_{[0]} -\iota_X \alpha(X)_{[2]}  \\
& \dots  \\
\beta(X)_{[j]} &= d\alpha(X)_{[j-1]} -\iota_X \alpha(X)_{[j+1]}  \\
& \dots  \\
\beta(X)_{[n-1]} &= d\alpha(X)_{[n-2]} -\iota_X \alpha(X)_{[n]}  \\
\beta(X)_{[n]} &= d\alpha(X)_{[n-1]}  \\
\beta(X)_{[n+1]} &= d\alpha(X)_{[n]}
\end{align*}
Note that in the equation $d_{eq}\alpha=0$ even components ``interact''
with even ones and odd components ``interact'' with odd ones.
So $d_{eq}\alpha=0$ implies $\alpha=\alpha^{even}+\alpha^{odd}$ with
$d_{eq} \alpha^{even}=0$, $d_{eq} \alpha^{odd}=0$. Moreover,
$$
d_{eq} (\alpha \wedge \beta) = (d_{eq} \alpha) \wedge \beta
+ (-1)^{\text{parity of $\alpha$}} \alpha \wedge (d_{eq} \beta)
$$
if $\alpha \in \Omega^{\infty, even}$ or $\Omega^{\infty, odd}$.
Hence the wedge product induces a product structure on
equivariant cohomology, and the theorem states that the isomorphism
(\ref{Cartan_thm}) is a ring isomorphism.

Observe that if $d_{eq}\alpha=0$, then $d\alpha(X)_{[n]}=0$,
$d\alpha(X)_{[n-1]}=0$ and $\iota_X\alpha(X)_{[1]}=0$
for all $X \in \mathfrak{g}$.
For a fixed $X \in \mathfrak{g}$ and $p \in M$ such that $L_X \bigr|_p =0$,
the evaluation map $\alpha \mapsto \alpha(X)_{[0]} \bigr|_p$
descends to cohomology.
Similarly, the integration map
$\int_M: \Omega^{\infty,*}_G(M) \to {\cal C}^{\infty} (\mathfrak{g})$
so that
$\alpha \mapsto \int_M \alpha =_{def} \int_M \alpha(X)_{[\dim M]}$
descends to cohomology as well.
(We assume that the manifold $M$ is oriented and compact.)

\begin{rem}
Another commonly considered variant of equivariant cohomology is
equivariant cohomology with distributional or ${\cal C}^{-\infty}$
coefficients. This is the cohomology of the $\BB Z_2$-graded complex
of $\Omega^*(M)$-valued distributions on $\mathfrak{g}$ denoted by
$\Omega^{-\infty,*}_G(M)$.
The definition of the twisted de Rham differential has to be modified
slightly in order to extend it to $\Omega^{-\infty,*}_G(M)$.
Fix a basis $\{X_1,\dots,X_l\}$ of $\mathfrak{g}$, and let
$\{X^1,\dots,X^l\}$ be the dual basis of $\mathfrak{g}^*$; we can think of
$X^1,\dots,X^l$ as linear functions on $\mathfrak{g}$.
Let $\alpha \in \Omega^{-\infty,*}_G(M)$ an equivariant form and let
$\phi \in {\cal C}^{\infty}_0(\mathfrak{g})$ be a test function, then
$$
\langle d_{eq} \alpha, \phi \rangle =_{def}
d \langle \alpha, \phi \rangle -
\sum_{j=1}^l \iota_{X_j} \langle \alpha, X^j \phi \rangle.
$$
This definition of $d_{eq}$ does not depend on the choice of basis
$\{X_1,\dots,X_l\}$ of $\mathfrak{g}$ and extends the old one.
For example, let $\beta \in \Omega^*(M)$ be a closed $G$-invariant form,
then the distribution $\beta: \phi \mapsto  \phi(0)\beta$
is an equivariantly closed element of $\Omega^{-\infty,*}_G(M)$.
Perhaps surprisingly, this distributional cohomology
$H^{-\infty,*}_G(M)$ may be quite different from $H^{\infty,*}_G(M)$.
\end{rem}

\separate

\subsection{Equivariant Cohomology with a Parameter}

Recall that $d_{eq} = d - \iota_X$. Sometimes people consider a slightly more
general twisted de Rham differential $d_c = d + c \iota_X$, where
$c \in \BB C$ is a fixed parameter.
Thus $d_{eq}=d_{-1}$; for $c=0$ we get the ordinary de Rham differential
(but the complex is different from the usual de Rham complex).

\begin{prop}
For $c \ne 0$, the cohomology rings of $(\Omega^{\infty,*}_G(M),d_c)$
(or $(\Omega^*_G(M),d_c)$) are isomorphic to each other.
\end{prop}

\pf
We will show that $H^*(\Omega^*_G(M),d_c) \simeq H^*(\Omega^*_G(M),d_{eq})$.
Pick an $s \in \BB C$ such that $s^2=-c$.
To an $\alpha \in \Omega^*_G(M)$,
$\alpha(X)=\alpha(X)_{[0]} + \alpha(X)_{[1]} + \dots + \alpha(X)_{[n]}$,
we associate
$$
\tilde \alpha(X) = \alpha(X)_{[0]} + s \alpha(X)_{[1]}
+ s^2 \alpha(X)_{[2]} + \dots + s^n \alpha(X)_{[n]}
\qquad \in \Omega^*_G(M).
$$
Then
\begin{multline*}
(d_{eq}\tilde\alpha(X))_{[j]}
= d(\tilde\alpha(X)_{[j-1]}) - \iota_X(\tilde\alpha(X)_{[j+1]})
= s^{j-1} d(\alpha(X)_{[j-1]}) - s^{j+1} \iota_X(\alpha(X)_{[j+1]})  \\
= s^{j-1} ( d(\alpha(X)_{[j-1]}) +c \iota_X(\alpha(X)_{[j+1]}))
= s^{j-1} (d_c\alpha(X))_{[j]} = s^{-1} (\widetilde{d_c \alpha}(X))_{[j]}
\end{multline*}
and $d_{eq}\tilde\alpha=s^{-1}\widetilde{d_c\alpha}$.
In particular, the scaling map
$\alpha \mapsto \tilde\alpha$ sends $d_c$-closed and $d_c$-exact forms
into respectively $d_{eq}$-closed and $d_{eq}$-exact forms.
Hence we get an isomorphism of cohomology
$H^*(\Omega^*_G(M),d_c) \simeq H^*(\Omega^*_G(M),d_{eq})$.
\qed

\separate

\begin{section}
{Symplectic Manifolds}
\end{section}

In this section and the next one we describe two of the most common sources
of equivariant forms. This section shows how a Hamiltonian system gives raise
to an equivariant symplectic form.

\separate

\subsection{The Equivariant Symplectic Form}

Let $(M,\omega)$ be a symplectic manifold, i.e. $M$ is equipped with a
$2$-form $\omega$ which is closed $(d\omega=0$) and satisfies a certain
non-degeneracy condition: for each $x \in M$, the bilinear form on the
tangent space $T_xM$ determined by $\omega$ is non-degenerate.
If $f \in {\cal C}^{\infty}(M)$ is a smooth function, the
{\em Hamiltonian vector field} generated by $f$ is the unique vector field
$H_f$ such that
$$
df = \iota(H_f) \omega.
$$

Now let $G$ act on $M$ and suppose that the action preserves $\omega$.

\begin{df}
$(G \lefttorightarrow M, \omega, \mu)$ is a \em{Hamiltonian system} if
$(M,\omega)$ is a symplectic manifold, $G$ preserves $\omega$ and the
map $\mu: M \to \mathfrak{g}^*$ called the {\em symplectic moment map}
satisfies:
\begin{enumerate}
\item
The map $\mu: M \to \mathfrak{g}^*$ is $G$-equivariant, i.e.
$$
\mu(g \cdot x) = Ad^{\#}(g) \mu(x),
\qquad \forall g \in G, \: x \in M;
$$

\item
For each $X \in \mathfrak{g}$, we obtain a function
$\langle \mu(x), X \rangle$ on $M$ which we denote by $\mu(X)$
and the Hamiltonian vector field by this function is $L_X$, i.e.
$$
d\mu(X)=\iota(L_X) \omega = \iota_X \omega.
$$
\end{enumerate}
\end{df}

\begin{ex}  \label{Bott2}
Let us return to the setting of Example \ref{Bott}.
Thus $G=S^1$ act on the sphere
$M=S^2 = \{(x,y,z) \in \BB R^3;\:x^2+y^2+z^2=1\}$
by rotation about the $z$-axis.
The sphere $S^2$ has a natural symplectic form
$\omega = \iota(N) dx \wedge dy \wedge dz$,
where $N$ is the outward-pointing unit normal vector field,
$N_{(x,y,z)}=(x,y,z)$. Thus
$$
\omega = x dy\wedge dz - y dx \wedge dz + z dx \wedge dy.
$$
It is clear from the construction that $\omega$ is rotational-invariant.

The Lie algebra of $S^1$ is one-dimensional and can be identified with $\BB R$.
The vector fields $L_X$, $X \in Lie(S^1)$, are proportional to
$$
-y \frac{\partial}{\partial x} + x \frac{\partial}{\partial y}.
$$
We choose the identification $Lie(S^1) \simeq \BB R$ so that $t \in \BB R$
generates the vector field
$L_t = ty \frac{\partial}{\partial x} - tx \frac{\partial}{\partial y}$.
Then
$$
\iota(L_t) \omega = t(x^2 + y^2) dz - tz(xdx+ydy) = t(x^2+y^2+z^2) dz = tdz
$$
since $x^2+y^2+z^2=1$ and $xdx+ydy+zdz=0$ on $S^1$.
If we set $\mu(t): S^2 \to \BB R$, $\mu(t)=tz$, we obtain
the symplectic moment map which makes $S^1 \lefttorightarrow S^2$
into a Hamiltonian system. Note that if we shift the moment map by
$c \in \BB R$ and set $\mu_c(t)=t(z-c)$ we get another moment map
such that $(S^1 \lefttorightarrow S^2, \omega,\mu_c)$ is a Hamiltonian system.
\end{ex}

Now we are ready to make an important observation.
The last equation in the definition is equivalent to saying that
$$
\tilde \omega =_{def} \mu + \omega : \quad
\mathfrak{g} \ni X \mapsto \mu(X) + \omega
\quad \in \Omega^0(M) \oplus \Omega^2(M)
$$
is an equivariantly closed form called the {\em equivariant symplectic form}.
In fact, $\tilde \omega = \mu + \omega \in \Omega^2_G(M)$.
Thus it makes sense to discuss the cohomology class
$[\tilde\omega] \in H^2_G(M)$ represented by $\tilde\omega$.

Equivariantly closed forms form a ring, so the forms
$$
\tilde \omega^n = (\mu + \omega)^n, \qquad n \in \BB N,
$$
are equivariantly closed.
Moreover, the form
\begin{align*}
e^{\tilde\omega} = e^{\mu+\omega} &=_{def}
1 + \frac{\tilde\omega}{1!} + \frac{\tilde\omega^2}{2!}
+ \frac{\tilde\omega^3}{3!} + \dots  \\
&= e^{\mu} \Bigl(1 + \frac{\omega}{1!} + \frac{\omega^2}{2!}
+ \frac{\omega^3}{3!} + \dots + \frac{\omega^k}{k!} \Bigr),
\qquad \text{where } k = \frac{\dim M}2,
\end{align*}
is also equivariantly closed.
This equivariant form plays a crucial role in the Duistermaat-Heckman Theorem
and demonstrates that working with just polynomial equivariant forms
$\Omega^*_G(M)$ is not sufficient and one has to consider equivariant
forms with ${\cal C}^{\infty}$ coefficients as well.

\separate

\subsection{Cotangent Spaces}  \label{sigma}

In this subsection we describe the most well-known Hamiltonian system.
Let $G$ act on a manifold $N$, this action induces a vector bundle action
on $T^*N$. We set $M=T^*N$ (of course such $M$ is not compact).
Then $M=T^*N$ has a canonical symplectic form $\sigma$.
Since we need to pin down all signs we spell out the definition
of $\sigma$ here.
Let $(x_1,\dots,x_n)$ be a system of local coordinates on $N$.
The $x_j$'s are smooth functions defined on some open subset $U \subset N$,
satisfying $dx_1 \wedge \dots \wedge dx_n \ne 0$ on $U$.
At each point $x \in U$, $(dx_1,\dots,dx_n)$ defines a basis of the vector
space $T^*_xN$, and a vector $\xi \in T^*_xN$ can be uniquely written as
$\xi = \sum_{j=1}^n \xi_jdx_j$. Then $(x_1,\dots,x_n;\xi_1,\dots,\xi_n)$
is a system of coordinates of $T^*N$ associated to the coordinate system
$(x_1,\dots,x_n)$. The form $\sigma$ is given by
$\sum_{j=1}^n d\xi_j \wedge dx_j$.
One can check that $\sigma$ is independent of the choice of coordinates
and thus gives a canonical symplectic structure on $M=T^*N$.
The moment map $\mu: M \to \mathfrak{g}^*$ is defined by
$$
\mu(X)(\xi) = - \langle \xi, L_X \rangle,
\qquad X \in \mathfrak{g},\: \xi \in M=T^*N.
$$
We leave the proof of the following lemma to the reader.

\begin{lem}
We have $d\mu(X)=\iota_X\sigma$, and
$(G \lefttorightarrow T^*N, \sigma, \mu)$ is a Hamiltonian system.
\end{lem}

\separate

\subsection{Coadjoint Orbits}

For details on this subsection see \cite{Ki} and \cite{BGV}`.
Recall that $G$ acts on its Lie algebra $\mathfrak{g}$ by the adjoint action
which in turn induces the coadjoint action on $\mathfrak{g}^*$.
Pick a point $\lambda \in \mathfrak{g}^*$, and let $M= G \cdot \lambda$
be the coadjoint orbit.
Like on any homogeneous space, the vector fields $L_X$, $X \in \mathfrak{g}$,
are sections of $TM$ and span the tangent space at each point.
We define the canonical $2$-form $\sigma_{\lambda}$ on $M$ by setting for
$l \in M$, $\tilde X, \tilde Y \in T_lM$,
$$
\sigma_{\lambda} \bigr|_l (\tilde X, \tilde Y) = - l([X,Y]),
$$
where $X,Y \in \mathfrak{g}$ are such that $L_X \bigr |_l = \tilde X$ and
$L_Y \bigr|_l = \tilde Y$.
Since
\begin{equation}  \label{z}
L_X \bigr|_l = 0 \quad \text{if and only if} \quad l([X,Z])=0 \quad
\forall Z \in \mathfrak{g},
\end{equation}
the $2$-form $\sigma_{\lambda}$ is well-defined.
It is called the Kirillov-Kostant-Souriau form.
The form $\sigma_{\lambda}$ is $G$-invariant. One can check
that $d\sigma_{\lambda}=0$ directly using the Jacobi identity,
but we will provide an alternative proof later.
Finally, $\sigma_{\lambda}$ is non-degenerate because (\ref{z})
implies $\sigma_{\lambda} \bigr|_l (L_X, \tilde Y) =0$ for all
$\tilde Y \in T_lM$ if and only if $L_X \bigr|_l = 0$.
In particular, all coadjoint orbits are even-dimensional.


\begin{prop}  \label{symplectic-form}
The Kirillov-Kostant-Souriau form $\sigma_{\lambda}$ defines a symplectic
$G$-invariant form on the coadjoint orbit $M= G \cdot \lambda$.
Moreover, $(G \lefttorightarrow M, \sigma_{\lambda},\mu)$, where
$\mu: M \hookrightarrow \mathfrak{g}^*$ is the inclusion map,
is a Hamiltonian system.
\end{prop}

\pf
First we check that $d\mu(X)= \iota_X \sigma_{\lambda}$.
Since $\mu(X)$ is the restriction of a linear function, we have:
$$
d\mu(X) (L_Y) \bigr|_l = (L_Y) \bigr|_l \mu(X) = - l([X,Y])
= \sigma(L_X,L_Y) \bigr|_l = (\iota_X \sigma_{\lambda}) (L_Y) \bigr|_l.
$$

Next we check $d\sigma_{\lambda}=0$:
$$
0 = d^2 \mu(X) = d \iota_X \sigma_{\lambda} = - \iota_X d\sigma_{\lambda},
$$
where we used
$(d \iota_X + \iota_X d) \sigma_{\lambda} = L_X \sigma_{\lambda} =0$
since $\sigma_{\lambda}$ is $G$-invariant.
Since the vector fields $L_X$ span the tangent spaces at each point,
this implies $d \sigma_{\lambda} =0$.
\qed

\separate

\subsection{Symplectic Reduction}

For details on this subsection see \cite{GS}.
Let $(G \lefttorightarrow M, \omega, \mu)$ be a Hamiltonian system.
Since the moment map $\mu : M \to \mathfrak{g}$ is $G$-equivariant,
the set $M_0 =_{def} \mu^{-1}(0)$ is $G$-invariant.
If the group $G$ acts on $M$ freely, then one can show that $\mu$ is
a submersion, i.e. the tangent map at each point $x \in M$
$$
d\mu_x: T_xM \to T_{\mu(x)}\mathfrak{g}^* \cong \mathfrak{g}^*
$$
is surjective.
Then $M_0$ is a smooth submanifold of $M$ on which $G$ acts freely.
Let $\tilde M_0 = M_0/G$, it is a smooth manifold with trivial $G$-action
of dimension $(\dim M -2\dim \mathfrak{g})$.
We have a natural projection $\pi_0$ and an inclusion $i_0$:
$$
\begin{CD}
M_0   @>{i_0}>>   M  \\
@V{\pi_0}VV    \\
\tilde M_0.
\end{CD}
$$
The next theorem claims that the manifold $\tilde M_0$ has a canonical
symplectic form.

\begin{thm} [Marsden-Weinstein]
Let $G$ act on $M$ freely. Then there exists a unique symplectic form
$\tilde \omega_0$ on $\tilde M_0$ with the property
$$
\pi_0^* \tilde \omega_0 = i_0^*\omega.
$$
\end{thm}

The operation of passing from $(M,\omega)$ to $(\tilde M_0, \tilde\omega_0)$
is known as {\em symplectic reduction} or {\em Marsden-Weinstein reduction}.

If the group $G$ is abelian, then instead of zero level of the moment map
we can consider any other level as well.
For every $a \in \mathfrak{g}^*$, the submanifold $M_a =_{def} \mu^{-1}(a)$
is $G$-invariant. Let $\tilde M_a = M_a/G$, it is also a smooth manifold
with trivial $G$-action of dimension $(\dim M -2\dim \mathfrak{g})$.
As before, we have a projection $\pi_a$ and an inclusion $i_a$:
$$
\begin{CD}
M_a   @>{i_a}>>   M  \\
@V{\pi_a}VV    \\
\tilde M_a.
\end{CD}
$$
The manifold $\tilde M_a$ also has a unique symplectic form defined by
the property $\pi_a^* \tilde \omega_a = i_a^*\omega$.

\separate

\begin{section}
{Vector Bundles}
\end{section}

In this section we describe another source of equivariant forms.
Vector bundles on manifolds have cohomological invariants called
characteristic classes which have an explicit description in terms
of differential forms. If a vector bundle happens to be $G$-equivariant,
then one can associate to it elements of equivariant cohomology
which can be described in terms of equivariant forms.
For simplicity we will describe the equivariant Chern-Weil and Euler
classes only, although other classes, notably the Thom class,
also have their equivariant analogues. See \cite{BGV} for details.
The equivariant Thom class construction is due to V.~Mathai and D.~Quillen
\cite{MQ}.

\separate

\subsection{Characteristic Classes via the Curvature Form}

In this subsection we recall how the Chern-Weil classes of a complex vector
bundle can be described in terms of the curvature form.
For the moment we ignore the $G$-action.

Let ${\cal E} \to M$ be a complex vector bundle.
We denote by $\Omega^*(M, {\cal E}) = \Omega^*(M) \otimes {\cal E}$
the bundle of ${\cal E}$-valued (or twisted) differential forms on $M$.
If ${\cal F}$ is a vector bundle, $\Gamma(M, {\cal F})$ denotes the vector
space of global sections of ${\cal F}$ on $M$.

\begin{df}
A {\em covariant derivative} on ${\cal E}$ is a differential operator
$$
\nabla: \Gamma(M, {\cal E}) \to \Gamma(M, T^*M \otimes {\cal E})
$$
which satisfies the Leibnitz's rule; that is, if $s \in \Gamma(M, {\cal E})$
and $f \in {\cal C}^{\infty}(M)$, then
$$
\nabla(fs) =  df \otimes s + f \nabla s.
$$
\end{df}

Note that a covariant derivative extends in a unique way to a map
$$
\nabla: \Omega^*(M, {\cal E}) \to \Omega ^{*+1}(M, {\cal E})
$$
that satisfies the Leibnitz's rule: if $\alpha \in \Omega^k(M)$
and $\beta \in \Omega^*(M, {\cal E})$, then
$$
\nabla(\alpha \wedge \beta) =
d\alpha \wedge \beta + (-1)^k \alpha \wedge \nabla \beta.
$$

Let $U \subset M$ be an open subset over which the bundle ${\cal E}$
trivializes: ${\cal E} \bigr|_U = U \times E$. Then, over $U$,
the covariant derivative may be written as
$$
\nabla s = ds + \omega_U \cdot s,
$$
for some
$\omega_U \in \Omega^1(U,\End({\cal E})) = \Omega^1(U) \otimes \End(E)$.

The {\em curvature} of a covariant derivative $\nabla$ is the
$\End({\cal E})$-valued $2$-form on $M$ defined by
$$
\Theta(F_1,F_2) = [\nabla_{F_1},\nabla_{F_2}]-\nabla_{[F_1,F_2]}
$$
for two vector fields $F_1$, $F_2$.
It is easy to check that in the above trivialization
${\cal E} \bigr|_U = U \times E$ the curvature will be given by the
$\End({\cal E})$-valued $2$-form on $U$
$$
\Theta = d\omega_U + \omega_U \wedge \omega_U.
$$

\begin{prop}
The operator $\nabla^2: \Omega^*(M, {\cal E}) \to \Omega ^{*+2}(M, {\cal E})$
is given by the action of $\Theta \in \Omega^2(M, \End({\cal E}))$ on
$\Omega^*(M,{\cal E})$.
\end{prop}

(This proposition is proved by verifying that over $U$ the operator
$\nabla^2$ is also given by $d\omega_U + \omega_U \wedge \omega_U$.)

The trace map induces a map $\tr: \Omega^*(M, \End({\cal E})) \to \Omega^*(M)$.
Let $P(z)$ be a polynomial in variable $z$. Then we can form elements
$P(\Theta) \in \Omega^*(M, \End({\cal E}))$ and
$\tr P(\Theta) \in \Omega^*(M)$.
The differential form $\tr P(\Theta)$ is called the {\em characteristic form}
(or Chern-Weil form) of $\nabla$ corresponding to the polynomial $P(z)$.

\begin{prop}
The characteristic form $\tr P(\Theta)$ is a closed differential form.
Moreover, if $\nabla_1$ and $\nabla_2$ are two connections on ${\cal E}$,
the differential forms $\tr P(\Theta_1)$ and $\tr P(\Theta_2)$ lie
in the same de Rham cohomology class.
\end{prop}

\begin{rem}
Topologists prefer to consider differential forms defined by the formula
$\tr P(-\Theta/2\pi i)$ because the Chern class, the cohomology class of
$c(\Theta) = \det(1-\Theta/2\pi i)$, defines an element of integral
cohomology class of $M$ only if these negative powers of $2\pi i$ are included.
\end{rem}

Instead of polynomials we can also allow power series $P(z)$ with infinite
radius of convergence.
If we take $P(z) = \exp(-z)$, then we get the characteristic form called the
{\em Chern character} form:
$$
ch(\nabla) = \tr(e^{-\Theta}).
$$
The cohomology class of $ch(\nabla)$ is independent of the choice of
connection $\nabla$ on ${\cal E}$ and will be denoted by $ch({\cal E})$.
We have:
\begin{align*}
ch({\cal E}_1 \oplus {\cal E}_2) &= ch({\cal E}_1) + ch({\cal E}_2),  \\
ch({\cal E}_1 \otimes {\cal E}_2) &= ch({\cal E}_1) \wedge ch({\cal E}_2).
\end{align*}
These two formulas explain why the cohomology class $ch({\cal E})$
is called a character.

\separate

\subsection{Equivariant Characteristic Classes}

Now let ${\cal E} \to M$ be a complex $G$-equivariant vector bundle.
This means that $G$ acts on ${\cal E}$ by vector bundle maps.
We define the space of $G$-equivariant differential forms with values in
${\cal E}$ by
$$
\Omega^*_G(M,{\cal E}) = ( \BB C[\mathfrak{g}] \otimes \Omega^*(M,{\cal E}) )^G
$$
with $\BB Z$-grading defined as on $\Omega^*_G(M)$.
Let $\nabla :  \Omega^*(M,{\cal E}) \to \Omega^{*+1}(M,{\cal E})$ be a
connection which commutes with $G$-action on $\Omega^*(M,{\cal E})$.
We say that $\nabla$ is a {\em $G$-invariant connection}.
Since the group $G$ is compact, starting with a possibly non-invariant
connection and averaging it over the Haar measure on $G$ one can obtain a
$G$-invariant connection. In particular, $G$-invariant connections exist.

\begin{df}
The {\em equivariant connection $\nabla_{eq}$} corresponding to a
$G$-invariant connection $\nabla$ is the operator on
$\BB C[\mathfrak{g}] \otimes \Omega^*(M,{\cal E})$ defined by the formula
$$
(\nabla_{eq} \alpha)(X) = (\nabla - \iota_X) \alpha(X),
\qquad X \in \mathfrak{g}.
$$
\end{df}

This definition is justified by the following version of the Leibnitz's rule:
$$
\nabla_{eq} (\alpha \wedge \beta) =
d_{eq} \alpha \wedge \beta + (-1)^{\deg \alpha} \alpha \wedge \nabla_{eq} \beta
$$
for all $\alpha \in \BB C[\mathfrak{g}] \otimes \Omega^*(M)$
and $\beta \in \BB C[\mathfrak{g}] \otimes \Omega^*(M,{\cal E})$.
The operator $\nabla_{eq}$ preserves the subspace
$\Omega^*_G(M,{\cal E}) \subset
\BB C[\mathfrak{g}] \otimes \Omega^*(M,{\cal E})$.

The $G$-action on ${\cal E}$ and $\Omega^*(M,{\cal E})$ induces
infinitesimal actions of $\mathfrak{g}$ which we denote by $L_X^{\cal E}$
for $X \in \mathfrak{g}$.

\begin{df}
The {\em equivariant curvature $\Theta_{eq}$} of an equivariant 
connection $\nabla_{eq}$ is defined to be
$$
\Theta_{eq}(X) = \nabla_{eq}(X)^2 + L_X^{\cal E},
\qquad X \in \mathfrak{g}.
$$
Thus $\Theta_{eq}$ is an element of $\Omega^*_G(M,\End({\cal E}))$.
\end{df}

The appearance of the term $L_X^{\cal E}$ can be justified as follows.
Suppose that ${\cal E} = M \times \BB C$ and $\nabla=d$ is the ordinary
de Rham operator. Then $\nabla_{eq} = d_{eq}$ and
$(d_{eq})^2 = -(d \iota_X + \iota_X d) = -L_X$, hence 
$\Theta_{eq}(X) = \nabla_{eq}(X)^2 + L_X =0$ in this case, as expected.

If we expand the definition of $\Theta_{eq}$ we get:
$$
\Theta_{eq}(X) = (\nabla - \iota_X)^2 + L_X^{\cal E}
= \Theta - (\iota_X \nabla + \nabla \iota_X) + L_X^{\cal E}
= \Theta - \nabla_X + L_X^{\cal E},
$$
where $\Theta = \nabla^2$ is the curvature of $\nabla$.
In particular, $\Theta_{eq}(0) = \Theta$.
It is convenient to define
$$
\mu(X) = \Theta_{eq}(X) - \Theta = L_X^{\cal E} - \nabla_X,
$$
$\mu(X) \in \Gamma(M,\End({\cal E}))$ and
$$
\mu \in ( \mathfrak{g}^* \otimes \Gamma(M, \End({\cal E})))^G \subset
\Omega^2_G(M,\End({\cal E})).
$$
Observe that if $X \in \mathfrak{g}$, $p \in M$ are such that
$L_X \bigr|_p =0$, then
\begin{equation}  \label{zero}
\mu(X) \bigr|_p = L_X^{\cal E} \bigr|_p.
\end{equation}

Now we can construct the equivariant analogues of characteristic forms.
The trace map induces a map
$\tr: \Omega^*_G(M, \End({\cal E})) \to \Omega^*_G(M)$.
Let $P(z)$ be a polynomial in variable $z$. Then we can form elements
$P(\Theta_{eq}) \in \Omega^*_G(M, \End({\cal E}))$ and
$\tr P(\Theta_{eq}) \in \Omega^*_G(M)$.
The differential form $\tr P(\Theta_{eq})$ is called the {\em equivariant
characteristic form} of $\nabla_{eq}$ corresponding to the polynomial $P(z)$.

\begin{prop}
The equivariant differential form $\tr P(\Theta_{eq})$ is equivariantly closed,
and its equivariant cohomology class is independent of the choice of
$G$-invariant connection $\nabla$ on ${\cal E}$.
\end{prop}

\pf
We will sketch a proof that the form $\tr P(\Theta_{eq})$
is equivariantly closed.
In the non-equivariant case we have
$$
d \tr \alpha = \tr ([\nabla, \alpha]),
\qquad \forall \alpha \in \Omega^*(M,\End({\cal E})),
$$
which we can verify locally.
So let $U \subset M$ be an open subset over which ${\cal E}$ trivializes,
and write the covariant derivative as $\nabla = d + \omega_U$,
for some $\omega_U \in \Omega^1(U,\End({\cal E}))$.
Then
$$
\tr ([\nabla, \alpha]) = \tr ([d, \alpha]) + \tr ([\omega_U, \alpha])
= \tr (d\alpha).
$$

Similarly, in the equivariant case we have
$$
d_{eq} \tr \alpha = \tr (\nabla_{eq}\alpha),
\qquad \forall \alpha \in \Omega^*_G(M,\End({\cal E})).
$$
The equation $d_{eq} (\tr P(\Theta_{eq}))=0$ now follows from
the equivariant Bianchi identity $\nabla_{eq} \Theta_{eq} =0$.
\qed

We can also allow $P(z)$ be a power series of infinite radius of convergence.
Taking $P(z) = \exp(-z)$, we get the {\em equivariant Chern character form
$ch_{eq}(\nabla)$} of an equivariant bundle ${\cal E}$:
$$
ch_{eq}(\nabla) = \tr(\exp(-\Theta_{eq})),
$$
in other words, $ch_{eq}(\nabla)(X) = \tr(\exp(-\Theta_{eq}(X)))$.

\separate

\subsection{Equivariant Euler Class}

In this subsection we define the equivariant Euler class which will
be needed in the statement of the Berline-Vergne integral localization formula.

Let $V$ be an oriented Euclidean space of dimension $n$.
We define the {\em Berezin integral} $\int: \Lambda^* V \to \BB R$ as follows:
first of all $\int \alpha =0$ if $\alpha \in \Lambda^k V$ for $k <n$;
secondly, $\int e_1 \wedge \dots \wedge e_n =1$ whenever $\{e_1,\dots,e_n\}$
is a positively oriented orthonormal basis for $V$; and thirdly,
if $\alpha \in \Lambda^n V$, then
$\alpha = a \cdot e_1 \wedge \dots \wedge e_n$ and the coefficient
$a \in \BB R$ does not depend on the choice of a positively oriented
orthonormal basis $\{e_1,\dots,e_n\}$ for $V$, we set $\int \alpha =a$.

\begin{df}
The {\em Pfaffian} of an element $\alpha \in \Lambda^2 V$ is the number
$$
\operatorname{Pf}_{\Lambda} (\alpha) = \int \exp \alpha =
\int (1 + \alpha/1! + \alpha^2/2! + \alpha^3/3! + \dots) =
\begin{cases}
\int \frac {\alpha^{n/2}}{(n/2)!}  & \text{if $n$ is even};\\
0 & \text{if $n$ is odd}.
\end{cases}
$$
\end{df}

If $A \in \mathfrak{so}(V)$, then it determines an anti-symmetric bilinear
form on $V$, namely $\langle Av_1, v_2 \rangle$, $v_1,v_2 \in V$, i.e. an
element $\alpha^* \in \Lambda^2 V^*$.
Identifying $V$ with $V^*$ we get an element $\alpha \in \Lambda^2V$.
Explicitly, choosing a positively oriented orthonormal basis
$\{e_1,\dots,e_n\}$ of $V$, we get
$$
\alpha = \sum_{i<j} \langle Ae_i, e_j \rangle e_i \wedge e_j
\quad \in \Lambda^2 V.
$$
We define the Pfaffian of an $A \in \mathfrak{so}(V)$ as
$$
\operatorname{Pf}(A) = \operatorname{Pf}(\alpha) =
\int ( \exp \sum_{i<j} \langle Ae_i, e_j \rangle e_i \wedge e_j).
$$
For example, if $V=\BB R^2$ with a positively oriented orthonormal basis
$\{e_1,e_2\}$ and if
$A = \begin{pmatrix} 0 & -\theta  \\ \theta & 0 \end{pmatrix}$,
then $\alpha = \theta \cdot e_1 \wedge e_2$ and
$\operatorname{Pf}(A) = \theta$.
Observe that reversing the orientation of $V$ changes the sign of 
$\operatorname{Pf}(A)$.

In general, the Pfaffian of an $A \in \mathfrak{so}(V)$ can be computed
by the following procedure. First of all assume that $n=\dim V$ is even
(otherwise $\operatorname{Pf}(A) = 0$). Then there exists a positively
oriented orthonormal basis $\{e_1,\dots,e_n\}$ of $V$ such that $A$
has a block-diagonal form, and each block is a $2 \times 2$ matrix of the type
$\begin{pmatrix} 0 & -c_j  \\ c_j & 0 \end{pmatrix}$ for some $c_j \in \BB R$.
Then, just as the $2$-dimensional example,
$\alpha = c_1 \cdot e_1 \wedge e_2 + \dots + c_{n/2} \cdot e_{n-1} \wedge e_n$
and $\operatorname{Pf}(A) = c_1 \dots c_{n/2}$.
In particular, it follows $\operatorname{Pf}(A)^2 = \det(A)$.
For this reason the notation
$$
{\det}^{1/2}(A) = \operatorname{Pf}(A)
$$
is often used.

Now let ${\cal E} \to M$ be a real $G$-equivariant oriented vector bundle
with $G$-invariant metric and $G$-invariant connection $\nabla$ compatible
with the metric. Recall that a connection is compatible with the metric if
it is torsion-free and, for any sections $s_1$, $s_2$ of ${\cal E}$, we have:
$$
d \langle s_1, s_2 \rangle = \langle \nabla s_1, s_2 \rangle +
\langle s_1, \nabla s_2 \rangle.
$$
The curvature $\Theta$ and the moment $\mu$ are both elements of
$\Omega^*_G(M, \mathfrak{so}({\cal E}))$. The {\em equivariant Euler form}
of ${\cal E}$ is defined by
$$
\chi_{eq}(\nabla)(X) = \operatorname{Pf}(-\Theta_{eq}(X))
= {\det}^{1/2}(-\Theta_{eq}(X)).
$$
Then $\chi_{eq}(\nabla)(X)$ is an equivariantly closed form, and its class
in equivariant cohomology depends neither on the connection nor on the
Euclidean structure of ${\cal E}$, but only on the orientation of ${\cal E}$.
The observation (\ref{zero}) in this context becomes:
if $X \in \mathfrak{g}$, $p \in M$ are such that $L_X \bigr|_p=0$, then
$$
\chi_{eq}(\nabla)(X)_{[0]} \bigr|_p =
\operatorname{Pf}(-L_X^{\cal E} \bigr|_p).
$$

\separate

\subsection{Topological Construction of Characteristic Classes}

In this subsection we give a topological construction of equivariant
characteristic classes following \cite{GS}.

Suppose first that $M$ and $N$ are two manifolds upon the group $G$ acts,
and let $f: M \to N$ be a smooth map intertwining the $G$-actions.
Then from the topological definition of equivariant cohomology we have
a natural pull-back map $H^*_G(N) \to H^*_G(M)$:
$$
H^*_G(N) = H^*((N \times E)/G) \to H^*((M \times E)/G) = H^*_G(M).
$$
On the level of equivariant differential forms we see that
$f^*: \Omega^*(N) \to \Omega^*(M)$ induces a pull-back map of complexes
$$
f^*: (\Omega^*_G(N), d_{eq}) \longrightarrow (\Omega^*_G(M), d_{eq})
$$
such that $(d_{eq} f^*(\alpha))(X) = (f^*(d_{eq}\alpha))(X)$,
$X \in \mathfrak{g}$.
Now taking $N=pt$ with trivial $G$-action we obtain a canonical map
$$
H^*_G(pt)= H^*(E/G) = \BB C[\mathfrak{g}]^G \longrightarrow H^*_G(M)
$$
which makes $H^*_G(M)$ into an algebra over
$H^*_G(pt) = \BB C[\mathfrak{g}]^G$.

Fix a compact group $K$ and consider a principal $K$-bundle
$P \twoheadrightarrow M = P/K$. From above considerations
we get a map
$$
\BB C[\mathfrak{k}]^K = H^*_K(pt) \longrightarrow
H^*_K(P) = H^*(M).
$$
The image of $\BB C[\mathfrak{k}]^K$ gives a subring of $H^*(M)$
called called the ring of characteristic classes of the principal
$K$-bundle $P \twoheadrightarrow M$.

Let ${\cal E} \twoheadrightarrow M$ be a complex vector bundle so that
each fiber is a complex space of dimension $n$.
Choose a Hermitian structure on ${\cal E}$, and let $P= {\cal F}({\cal E})$
denote the bundle of unitary frames. The group $K=U(n)$ acts on
${\cal F}({\cal E})$ by multiplications on the right freely making
${\cal F}({\cal E}) \twoheadrightarrow M$ into a principal
$U(n)$-bundle over $M$. (The reason for choosing multiplication on the right
rather than on the left will be apparent later.)
Hence we again get a map $\BB C[\mathfrak{u}(n)]^{U(n)} \to H^*(M)$
and its image consists of the Chern classes of the vector bundle ${\cal E}$.
This map is independent of the choice of Hermitian metric on ${\cal E}$.

We can do the same procedure with ${\cal E} \twoheadrightarrow M$
a real vector bundle with fibers isomorphic to $\BB R^n$ and $K=O(n)$
thus obtaining the Pontryagin classes.
Taking ${\cal E} \twoheadrightarrow M$ an oriented real vector bundle
with fibers isomorphic to $\BB R^{2n}$ and $K=SO(2n)$ we obtain the
Pontryagin classes as before plus another class called the Euler class
(or Pfaffian).

Now suppose that we have two compact groups $G$ and $K$ acting on a manifold
$P$ so that their actions commute and the $K$-action is free. Let $M=P/K$,
then $M$ is a manifold with $G$-action.
Pick a classifying space $E$ for the group $G$ and make $K$ act on it
trivially. Then the projection
$$
(P \times E)/G \twoheadrightarrow (P \times E)/(K \times G) = (M \times E)/G
$$
is a principal $K$-bundle, and we get a natural map
$\BB C[\mathfrak{k}]^K \to H^*_G(M)$ as a composition of
$$
\BB C[\mathfrak{k}]^K = H^*_K(pt) \longrightarrow
H^*_K((P \times E)/G) = H^*((M \times E)/G) = H^*_G(M).
$$

More concretely, let ${\cal E} \twoheadrightarrow M$ be a $G$-equivariant
complex vector bundle with fibers isomorphic to $\BB C^n$.
We choose a $G$-invariant Hermitian structure on ${\cal E}$ and repeat the
frame bundle construction observing that now everything becomes
$G$-equivariant.
Thus $P= {\cal F}({\cal E})$ is the bundle of unitary frames, and
the group $K=U(n)$ acts on ${\cal F}({\cal E})$ by multiplications on the
right freely making ${\cal F}({\cal E}) \twoheadrightarrow M$ into a principal
$U(n)$-bundle over $M$. Now this $K$-action commutes with $G$-action,
so by the above considerations we get a map
$\BB C[\mathfrak{u}(n)]^{U(n)} \to H^*_G(M)$
and its image consists of the equivariant Chern classes of the vector bundle
${\cal E}$. As before, this map is independent of the choice of a
$G$-invariant Hermitian metric on ${\cal E}$.
Similarly one can get $G$-equivariant Pontryagin and Euler classes.

If the vector bundle ${\cal E} \twoheadrightarrow M$ is topologically trivial,
its characteristic classes vanish. But this need not be true of its
equivariant characteristic classes. For example, consider a complex vector
bundle over a point $E \to pt$. This is just an ordinary $n$-dimensional
vector space on which $G$ acts by linear transformations.
Equipping $E$ with a $G$-invariant Hermitian metric, we can regard the
representation of $G$ on $E$ as a homomorphism $G \to K=U(n)$.
This homomorphism induces a map of Lie algebras
$\mathfrak{g} \to \mathfrak{k}$ which gives a homomorphism of the rings
of invariants:
$$
\BB C[\mathfrak{k}]^K \to \BB C[\mathfrak{g}]^G = H^*_G(pt),
$$
and the equivariant Chern classes are just the image of
$\BB C[\mathfrak{k}]^K$.

\separate

\begin{section}
{The Equivariant Thom Class}
\end{section}

\separate

\subsection{The Classical (Non-equivariant) Thom Class}

Our goal in this section is to describe the equivariant version of the
Thom form. As a motivation, we briefly recall the properties of the
classical Thom class.
Let $M$ be a $d$-dimensional oriented manifold (not necessarily compact),
and let $N$ be a compact oriented submanifold of codimension $k$.
The integration over $N$ defines a linear function on the de Rham
cohomology group $H^{d-k}(M)$. On the other hand, Poincar\'e duality
asserts that the pairing
$$
(\: , \: ): \: H^{d-k}(M) \times H^k_c(M) \to \BB C,
\qquad (a,b) = \int_M a \wedge b
$$
is non-degenerate, where $H^k_c(M)$ denotes the compactly supported cohomology
groups.
In particular, there exists a unique cohomology class $\tau(N) \in H^k_c(M)$
such that
$$
(2\pi)^{k/2} \int_N \alpha = (\alpha,\tau(N)) = \int_M \alpha \wedge \tau(N)
\qquad \forall \alpha \in H^{d-k}(M).
$$
This class is called the {\em Thom class} associated with $N$, and any
closed differential form $\tau_N$ representing this class is called a
{\em Thom form}. Thus a Thom form is a closed form with the
``reproducing property''
$$
\int_M \alpha \wedge \tau_N = (2\pi)^{k/2} \int_N \alpha
\qquad \text{for all $\alpha \in \Omega^{d-k}(M)$ with $d\alpha=0$.}
$$
Clearly any closed form of degree $k$ with this property is a Thom form.

Unfortunately, Poincar\'e duality is not available in equivariant cohomology.
For example, let $M$ be a compact oriented $G$-manifold, then we have a
natural bilinear map
$$
(\: , \: ): \: H^*_G(M) \times H^*_G(M) \to \BB C [\mathfrak{g}]^G,
\qquad (a,b)= \int_M a \wedge b.
$$
This pairing, however, can be highly singular.
For instance, if the action of $G$ on $M$ is free, according to the
Berline-Vergne integral localization formula (Theorem \ref{BV-loc})
the integral over $M$ of any equivariantly closed form is zero,
so the pairing is trivial. Therefore, one cannot define the equivariant
Thom class in equivariant de Rham theory simply by invoking Poincar\'e
duality as in the non-equivariant case.
For this reason we need to rely on an explicit construction of the Thom form.

\separate

\subsection{Fiber Integration of Equivariant Forms}

Let $\pi: V \to N$ be a $G$-equivariant fiber bundle; we assume that both
$V$ and $N$ are oriented. Let $m = \dim V$ and $n = \dim N$ so that
$$
k =_{def} m -n \ge 0.
$$
If $k=0$ we also assume that $\pi$ is orientation-preserving.
Let $\Omega^l_c(V)$ denote the space of compactly supported $l$-forms
on $V$ with a similar notation for $N$. If $l \ge k$, there is map
$$
\pi_*: \Omega^l_c(V) \to \Omega^{l-k}_c(N)
$$
called {\em fiber integration} where, for $\beta \in \Omega^l_c(V)$,
$\pi_* \beta$ is uniquely characterized by
\begin{equation}  \label{thom-eqn}
\int_V \pi^* \alpha \wedge \beta = \int_N \alpha \wedge \pi_* \beta
\qquad \forall \alpha \in \Omega^{m-l}_c(N).
\end{equation}
It is clear that $\pi_*\beta$ is uniquely determined by this condition,
since if there are two such forms their difference $\nu$ is
an $(l-k)$-form on $N$ with the property that
$$
\int_N \alpha \wedge \nu = 0
\qquad \alpha \in \Omega^{m-l}_c(N)
$$
hence it must vanish.

Once we know that our condition determines $\pi_*$ uniquely, it is sufficient
(by partitions of unity) to prove the existence in a coordinate patch where
$(x_1,\dots,x_n,t_1,\dots,t_k)$ are coordinates on $V$ extending coordinates
$(x_1,\dots,x_n)$ on $N$ so that $\pi$ is given by $\pi(x,t)=x$.
Then if
$$
\beta = \sum_I b_I(x,t) dx_I \wedge dt_1 \wedge \dots \wedge dt_k + \dots
$$
where the remaining terms involve fewer than $k$ wedge products of the $dt_i$,
we can easily check that
$$
\pi_* \beta = \sum_I
\Bigl( \int b_I(x,t) dt_1 \wedge \dots \wedge dt_k \Bigr) dx_I
$$
satisfies (\ref{thom-eqn}).

We have elementary but important properties of fiber integration:

\begin{prop}
\begin{enumerate}
\item
$d\pi_* = \pi_* d$.

\item
If $\tilde f : V \to V$ and $f: N \to N$ are proper orientation-preserving
maps making the diagram
$$
\begin{CD}
V @>{\tilde f}>> V \\
@V{\pi}VV      @VV{\pi}V \\
N @>f>> N
\end{CD}
$$
commute ($\pi \circ \tilde f = f \circ \pi$), then
$$
\pi_* \circ \tilde f^* = f^* \circ \pi_*.
$$
In particular, $\pi_*$ commutes with the $G$-action.

\item
The infinitesimal version of the last property is
$$
L_X \circ \pi_* = \pi_* \circ L_X
\qquad \forall X \in \mathfrak{g}.
$$

\item
Let $F_V$ and $F_N$ be vector fields on $V$ and $N$ respectively
such that $\pi_* F_V =F_N$, then
$$
\pi_* \iota(F_V) = \iota(F_N) \pi_*
$$

\item
$d_{eq} \pi_* = \pi_* d_{eq}$.

\item
The fiber integration map induces a map of complexes of
equivariant differential forms
$$
\pi_*: (\Omega^*_G(V), d_{eq}) \to (\Omega^{*-k}_G(N), d_{eq}).
$$

\item
If we have a pull back square of fiber bundles
$$
\begin{CD}
f^*V @>{\tilde f}>> V \\
@V{\pi'}VV      @VV{\pi}V \\
N' @>f>> N
\end{CD}
$$
then $\pi'_* \tilde f^* = f^* \pi_*$.

\item
If $\pi: V \to N$ is a vector bundle, then the induced map on de Rham
cohomology
$$
\pi_*: H^l_c(V) \to H^{l-k}_c(N)
$$
is an isomorphism for all $l$. In particular, $ H^l_c(V) =0$ for $l<k$.
\end{enumerate}
\end{prop}

The proof of this proposition relies on (\ref{thom-eqn}).
See \cite{GS} for details.

Next we observe that left multiplication of forms makes
$\Omega^*_c(N)$ into a module over $\Omega^*(N)$.
Since $\pi^*: \Omega^*(N) \to \Omega^*(V)$ is a homomorphism
(in fact an injection), this makes $\Omega^*_c(V)$ into an
$\Omega^*(N)$-module as well.

\begin{prop}
$\pi_*$ is a homomorphism of $\Omega^*(N)$-modules:
$$
\pi_*(\pi^* \alpha \wedge \beta) = \alpha \wedge \pi_* \beta
\qquad \forall \alpha \in \Omega^*(N).
$$
\end{prop}

We conclude this subsection by pointing out the relation between fiber
integration and Thom forms. Suppose that $N$ is a compact submanifold
of a manifold $M$, $i: N \hookrightarrow M$ the inclusion map,
and $U$ a tubular neighborhood of $N$ which we identify
with an open neighborhood of $N$ in the normal bundle ${\cal N}$ of $N$.
Let $\tau \in \Omega^k_c(U)$ be a closed form, where $k$ is the codimension
of $N$. Then $\pi_* \tau$ is a zero-form, i.e. a function on $N$.
Suppose that $\pi_* \tau = (2\pi)^{k/2}$ (that is integrating over every
fiber results in the value $(2\pi)^{k/2}$).
(It is known that differential forms with these properties exist.)
Extending by zero, we can regard $\tau$ as a compactly-supported $k$-form
on $M$. Then $\tau$ is a Thom form.
Indeed, since $N$ is a deformation retract of $U$, for any closed form
$\alpha$ on $M$, the restriction of $\alpha$ to $U$ is cohomologous to
$\pi^* i^* \alpha$. Hence, if $\alpha$ is a form on $M$ of degree $d-k=\dim N$,
$$
\int_M \alpha \wedge \tau
= \int_U \alpha \wedge \tau
= \int_U \pi^* i^* \alpha \wedge \tau
= \int_N i^* \alpha \wedge \pi_* \tau
= (2\pi)^{k/2} \int_N \alpha.
$$

\separate

\subsection{Equivariant Thom Forms}

Let $\pi: V \to N$ be a $G$-equivariant vector bundle.
We assume that $N$ is compact, both $V$ and $N$ are oriented,
$\dim V = m$ and $\dim N = n =m-k$.
If $k=0$ we also assume that $\pi$ is orientation-preserving.
Since $\pi_*$ is a cochain map for the twisted de Rham complex, it
induces a map on equivariant cohomology
\begin{equation}  \label{10.13}
\pi_* : (H^l_c)_G(V) \to H^{l-k}_G(N).
\end{equation}
Just as in the non-equivariant case, we have:

\begin{thm} [Mathai-Quillen]
The map (\ref{10.13}) is an isomorphism for all $l$.
\end{thm}

Taking the inverse image of $(2\pi)^{k/2} \in H^0_G(N)$ we obtain:

\begin{thm} [Mathai-Quillen, \cite{MQ}]
There exists a form $\tau_{eq} \in \Omega^k_G(V)$ with compact support
(i.e. $\tau_{eq} \in (\BB C [\mathfrak{g}] \otimes \Omega^*_c(V))^G$)
which is equivariantly closed ($d_{eq} \tau_{eq} =0$) and
$$
\pi_* \tau_{eq} = (2\pi)^{k/2}.
$$
Moreover, such a form $\tau_{eq}$ can be chosen in a semi-canonical way.
\end{thm}

This equivariant form $\tau_{eq}$ is called an equivariant Thom form
associated with the vector bundle $\pi: V \to N$. Note that the top component
of $\tau_{eq}$, $(\tau_{eq})_{[k]}$, is a non-equivariant Thom form.
We have another result that mimics the non-equivariant case.

\begin{thm} [Mathai-Quillen, \cite{MQ}]
Let $\tau_{eq}(V) \in (H^k_c)_G(V)$ and $\chi_{eq}(V) \in H^k_G(N)$ be the
equivariant Thom and equivariant Euler classes respectively associated
to the vector bundle $\pi: V \to N$.
Denote by $i_0: N \hookrightarrow V$ the embedding of $N$ into $V$
as the zero section. Then
$$
i_0^* \tau_{eq}(V) = \chi_{eq}(V).
$$
\end{thm}

Now let $M$ be an $d$-dimensional oriented manifold (not necessarily compact)
on which $G$ acts, and let $N$ be a compact oriented $G$-invariant submanifold
of dimension $n=d-k$. Let $i: N \hookrightarrow M$ denote the inclusion map.
Denote by ${\cal N}$ the normal bundle of $N$ in $M$. By the equivariant
tubular neighborhood theorem here exists a $G$-invariant tubular neighborhood
$U$ of $N$ in $M$, and a $G$-invariant diffeomorphism
$$
\gamma: {\cal N} \to U
\qquad \text{such that} \qquad
\gamma \circ i_0 = i,
$$
where $i_0: N \hookrightarrow {\cal N}$ is the embedding of $N$ into ${\cal N}$
as the zero section.
Let $\tau_{eq} \in (\Omega^k_c)_G(U)$ be the equivariant Thom form of
associated to ${\cal N}$. As before, we extend $\tau_{eq}$ by zero to an
equivariant form on $M$, then $\tau_{eq} \in (\Omega^k_c)_G(M)$ is an
equivariant Thom form associated with $M$ in the sense that,
for all $\alpha \in \Omega^{d-k}_G(M)$ with $d_{eq}\alpha=0$,
$$
\int_M \alpha \wedge \tau_{eq} =
\int_M \alpha_{[d-k]} \wedge (\tau_{eq})_{[k]} =
(2\pi)^{k/2} \int_N \alpha_{[d-k]} = 
(2\pi)^{k/2} \int_N \alpha.
$$
The main point here is that the non-equivariant Thom form $\tau$ can be chosen
to be $G$-invariant and ``completed'' to an equivariantly form $\tau_{eq}$
such that $d_{eq} \tau_{eq}=0$.
Moreover, such a form $\tau_{eq}$ can be chosen in a semi-canonical way.

\separate

\begin{section}
{Localization Theorems}
\end{section}

\separate

\subsection{Berline-Vergne Integral Localization Formula}

Since the group $G$ acting on a manifold $M$ is compact, there is a Riemannian
metric $\langle \cdot, \cdot \rangle$ on $M$ which is $G$-invariant.
Such metric may be constructed by picking one particular metric and averaging
it with respect to the Haar measure of $G$.
Most proofs of localization theorems rely on existence of such $G$-invariant
metrics and for this reason fail for non-compact groups.

\begin{prop}  \label{exact}
Let $G$ be a compact Lie group and $\alpha: \mathfrak{g} \to \Omega^*(M)$
be a ${\cal C}^{\infty}$ map such that $d_{eq} \alpha=0$.
(Note that we do not require $\alpha$ to be an equivariant form.)
Pick an element $X \in \mathfrak{g}$ and let $M_0(X)$ be the set of
zeroes of the vector field $L_X$.
Then the differential form $\alpha(X)_{[\dim M]}$ is exact on
$M \setminus M_0(X)$.
\end{prop}

\pf
Fix an $X \in \mathfrak{g}$.
Let $\theta_X \in \Omega^1(M)$ be a $1$-form on $M$ such that
$L_X \theta_X =0$ and $\iota(X) \theta_X \ne 0$ on $M \setminus M_0(X)$.
Then it follows that
$$
(d -\iota_X)^2 \theta_X=0.
$$
Such a differential form can be constructed by taking a $G$-invariant 
Riemannian metric $\langle \cdot, \cdot \rangle$ on $M$ and setting
$$
\theta_X (v) = \langle L_X, v \rangle,
\qquad v \in TM.
$$
Then $\theta_X \in \Omega^1(M)$ and it is easy to check that $L_X \theta_X =0$
and $\iota(X) \theta_X = \|L_X \|^2 \ne 0$ on $M \setminus M_0(X)$.

In general, if $A$ is an commutative algebra and $x$ is a nilpotent element,
then, $(1-x) \in A$ is also invertible and
$$
(1-x)^{-1} = 1+x+x^2+x^3+\dots
$$
(the sum is finite since $x$ is nilpotent).
Slightly more generally, for any invertible element $a \in A$,
$(-a+x) \in A$ is also invertible and
$$
(-a+x)^{-1} = -a^{-1} \cdot (1 - x/a)^{-1} =
-a^{-1} - a^{-2}x - a^{-3}x^2 - a^{-4}x^3 - \dots.
$$

These considerations apply to the element
$(d-\iota_X)\theta_X =  -\iota_X \theta_X + d\theta_X$
which is invertible on the set $M \setminus M_0(X)$:
$$
\frac 1{(d-\iota_X)\theta_X} =
-(\iota_X \theta_X)^{-1} - (\iota_X \theta_X)^{-2} d\theta_X
- (\iota_X \theta_X)^{-3} (d\theta_X)^2- (\iota_X \theta_X)^{-4} (d\theta_X)^3
- \dots
$$
(note that the sum is finite) so that
$$
\frac 1{(d-\iota_X)\theta_X} \wedge ((d-\iota_X)\theta_X)
= ((d-\iota_X)\theta_X) \wedge \frac 1{(d-\iota_X)\theta_X} = 1.
$$
Applying $(d-\iota_X)$ to
$((d-\iota_X)\theta_X) \wedge \frac 1{(d-\iota_X)\theta_X}$ we get
$$
((d-\iota_X)\theta_X) \wedge
(d-\iota_X) \Bigl( \frac 1{(d-\iota_X)\theta_X} \Bigr) =0.
$$
Since $(d-\iota_X)\theta_X$ is invertible, this implies
$$
(d-\iota_X) \Bigl( \frac 1{(d-\iota_X)\theta_X} \Bigr) =0.
$$
Therefore, on $M \setminus M_0(X)$ we have
$$
\alpha(X) = (d -\iota_X)
\Bigr( \frac{\theta_X \wedge \alpha(X)}{(d-\iota_X)\theta_X} \Bigr),
$$
and so
$$
\alpha(X)_{[\dim M]} = d
\Bigr( \frac{\theta_X \wedge \alpha(X)}{(d-\iota_X)\theta_X}
\Bigr)_{[\dim M-1]}.
$$
\qed

\begin{ex}
This example shows that the assumption that $G$ is compact is essential.
Consider $M = S^1 \times S^1$ -- the two-dimensional torus --
with coordinates $x,y \in \BB R / 2\pi \BB Z$.
Let $F$ be the nowhere-vanishing vector field
$(2 + \sin x) \frac{\partial}{\partial y}$, thus we obtain an action
of $G= \BB R$ on $M$. Let
$$
\alpha = -(7 \cos x + \sin 2x) + (1-4 \sin x) dx \wedge dy.
$$
It is easy to verify that
$$
\iota(F) \alpha_{[2]} =
(4\sin^2 x + 7 \sin x -2)dx = d \alpha_{[0]},
$$
and so $(d - \iota(F)) \alpha=0$.
However, $\int_M \alpha = (2\pi)^2$, so that $\alpha_{[2]}$ is not exact.
(By the way, this form $\alpha$ is not $G$-equivariant.)
\end{ex}

Proposition \ref{exact} strongly suggests that the integral of an
equivariantly closed form $\alpha \in \Omega^{\infty,*}_G(M)$ depends only
on the restriction of $\alpha(X)$ to $M_0(X)$. Such an integral localization
formula is due to N.~Berline and M.~Vergne and will be proved in this
subsection.

\begin{rem}
Suppose that the group $G$ in Proposition \ref{exact} is abelian
(hence a torus). Then the assignment $X \mapsto \theta_X$ defines an element
of $\Omega^3_G(M)$. The proof implies that $\alpha$ is equivariantly exact
on $M \setminus M_0(X)$.
Moreover, if $\alpha \in \Omega^*_G(M)$ is polynomial, then $\alpha$
is equivariantly exact on $M \setminus M_0(X)$ as an element of
$\Omega^*_G(M)$.
This observation strongly suggests that the equivariant cohomologies
$H^{\infty,*}_G(M)$ and $H^*_G(M)$ are expressible through the local geometry
at the fixed point set
$M^G = \{x \in M;\: g\cdot x = x \: \forall g \in G \}$.
The precise result is due to M.~Atiyah and R.~Bott and we will discuss it
in Subsection \ref{AB-loc}.
\end{rem}

We are now ready to prove the Berline-Vergne integral localization formula
in the important special case when the zeroes of the vector field $L_X$
-- denoted by $M_0(X)$ -- are isolated (hence a finite set).
Here, at each point $p \in M_0(X)$, the Lie action $L_X$ on $\Gamma(M,TM)$,
$F \mapsto [L_X,F]$, gives rise to a transformation $L(X,p)$ of $T_pM$.
Indeed, since $L_X \bigr|_p =0$, $[L_X,F]=0$ if and only if $F \bigr|_p =0$,
thus the transformation $L(X,p)$ is well-defined.

\begin{lem}
The transformation $L(X,p)$ on $T_pM$ is invertible.
\end{lem}

\pf
Pick a $G$-invariant Riemannian metric
and suppose that $v \in T_pM \setminus \{0\}$ is annihilated by $L(X,p)$.
Taking the exponential map with respect to the metric, we see that all
the points on the geodesic $\exp_p(tv)$, $t \in \BB R$, would be fixed by
$\exp (sX) \in G$, $s \in \BB R$,  contradicting to $p \in M_0(X)$
being isolated.
\qed

Select an abelian Lie subgroup (i.e. a torus) $T \subset G$ such that
$X \in Lie(T)$. Then $T$ fixes $p \in M_0(X)$, and thus we obtain
a representation of $T$ on $T_pM$. Since $T$ is compact, we can pick
a Euclidean metric on $T_pM$ which is preserved by $T$. Then
the transformation $L(X,p)$ belongs to $\mathfrak{so}(T_pM)$.
Since $L(X,p)$ is invertible, it follows that the manifold $M$ is
even-dimensional; let $n = \dim M$.
Recall that the Pfaffian of $L(X,p)$ can be computed as follows:
pick a positively oriented basis of $T_pM$ such that $L(X,p)$
has a block-diagonal form, and each block is a matrix of the type
$\begin{pmatrix} 0 & -c_j  \\ c_j & 0 \end{pmatrix}$ for some $c_j \in \BB R$,
then ${\det}(L(X,p)) = c_1^2 \dots c_{n/2}^2$ and
$$
{\det}^{1/2}(L(X,p)) = \operatorname{Pf}(L(X,p)) = c_1 \dots c_{n/2}.
$$

\begin{thm} [Berline-Vergne, 1982]  \label{BV-loc}
Let $G$ be a compact group with Lie algebra $\mathfrak{g}$ acting on a
compact oriented manifold $M$, and let
$\alpha: \mathfrak{g} \to \Omega^*(M)$
be a ${\cal C}^{\infty}$ map such that $d_{eq} \alpha=0$.
(We do not require $\alpha$ to be an equivariant form.)
Let $X \in \mathfrak{g}$ be such that the vector field $L_X$ has only
isolated zeroes. Then
$$
\int_M \alpha(X) = (-2\pi)^{\dim M/2} \sum_{p \in M_0(X)}
\frac {\alpha(X)_{[0]}(p)}{{\det}^{1/2}(L(X,p))},
$$
where by $\alpha(X)_{[0]}(p)$ we mean the value of the function
$\alpha(X)_{[0]}$ at the point $p \in M$.
\end{thm}

\begin{rem}
Note that the integration map
$\int_M : \alpha \mapsto \int_M \alpha(X)_{[\dim M]}$
and the evaluation map $\alpha \mapsto \alpha(X)_{[0]}(p)$ for $p \in M_0(X)$
both descend to equivariant cohomologies $H^*_G(M)$ and $H^{\infty,*}_G(M)$.
Thus the integral localization formula descends to the cohomological level
as well.
\end{rem}

\begin{ex}
Let us compute both sides of the integral localization formula in the setting
of Example \ref{Bott2} for $\alpha = \tilde \omega = \mu(t) + \omega$.
The orientation of a symplectic manifold is always determined by the highest
power of its symplectic form. Thus $\int_{S^2} (\mu(t) + \omega) = 4\pi$
-- the surface area of a unit sphere. On the other hand, $\mu(t)=tz$,
$M_0(X) = \{(0,0,\pm 1)\}$ (unless $X=0$) and ${\det}^{1/2}(L(X,p)) =\pm t$,
so that
$$
\sum_{p \in M_0(X)} \frac {\tilde \omega(t)_{[0]}(p)}{{\det}^{1/2}(L(X,p))}
= \sum_{p=(0,0,\pm 1)} \frac {\mu(t)(p)}{{\det}^{1/2}(L(X,p))}
= \frac t{-t} + \frac {-t}{t} = -2
$$
in total agreement with the localization formula.
\end{ex}

Now we can give a proof of the Berline-Vergne integral localization formula.
(It is copied from \cite{BGV}.)

\pf
By replacing $G$ with a torus $T$ such that $X \in Lie(T)$ we can assume that
$G$ is abelian. Fix a $G$-invariant Riemannian metric on $M$.
Let $p \in M_0(X)$, then $p$ is fixed by $G$, $G$ acts linearly on $T_pM$,
and there is a positively oriented orthonormal basis of $T_pM$ such that
$L(X,p)$ has a block-diagonal form, and each block is a matrix of the type
$\begin{pmatrix} 0 & -c_j  \\ c_j & 0 \end{pmatrix}$ for some $c_j \in \BB R$,
then ${\det}^{1/2}(L(X,p)) = c_1 \dots c_{n/2}$.
Using the exponential map $\exp_p: T_pM \to M$ we can construct local
coordinates $x_1, \dots,x_n$ around $p$ such that the vector field $L_X$
is expressed by
$$
L_X = c_1 \Bigl( x_2 \frac{\partial}{\partial x_1} -
x_1 \frac{\partial}{\partial x_2} \Bigr) + \dots +
c_{n/2} \Bigl( x_n \frac{\partial}{\partial x_{n-1}} -
x_{n-1} \frac{\partial}{\partial x_n} \Bigr).
$$

Let $\theta_X^p$ be the $1$-form in a neighborhood $U_p$ of $p$ given by
$$
\theta_X^p = c_1^{-1}(x_2 dx_1 - x_1 dx_2) + \dots +
c_{n/2}^{-1}(x_n dx_{n-1} - x_{n-1} dx_n).
$$
Then $\theta_X^p$ has the following properties: $L_X \theta_X^p =0$ and
$\iota_X \theta_X^p = \theta_X^p(L_X) = x_1^2 + \dots + x_n^2 = \|x\|^2$.
Using a $G$-invariant partition of unity subordinate to the covering of $M$
by the open sets $U_p$, $p \in M_0(X)$, and $M \setminus M_0(X)$
(which may be constructed by averaging any partition of unity with respect
to the action of $G$), we can construct a $1$-form $\theta_X$ such that
$L_X \theta_X=0$, $\iota_X \theta_X \ne 0$ outside $M_0(X)$, and such that
$\theta_X$ coincides with $\theta_X^p$ in a neighborhood of each
$p \in M_0(X)$.

Consider the neighborhood $B_{\epsilon}^p \subset M$ of $p$ given by
$B_{\epsilon}^p = \{(x_1,\dots,x_n) ;\: x_1^2 + \dots + x_n^2 \le \epsilon \}$,
and set $S_{\epsilon}^p = \partial B_{\epsilon}^p =
\{(x_1,\dots,x_n) ;\: x_1^2 + \dots + x_n^2 = \epsilon \}$.
We have:
\begin{multline*}
\int_M \alpha(X) = \lim_{\epsilon \to 0+}
\int_{M \setminus \cup_p B_{\epsilon}^p} \alpha(X)
= \lim_{\epsilon \to 0+} \int_{M \setminus \cup_p B_{\epsilon}^p}
d \Bigl( \frac{\theta_X \wedge \alpha(X)}{(d-\iota_X)\theta_X} \Bigr)  \\
= - \sum_{p \in M_0(X)} \lim_{\epsilon \to 0+} \int_{S_{\epsilon}^p}
\frac{\theta_X \wedge \alpha(X)}{(d-\iota_X)\theta_X}
= - \sum_{p \in M_0(X)} \lim_{\epsilon \to 0+} \int_{S_{\epsilon}^p}
\frac{\theta_X^p \wedge \alpha(X)}{(d-\iota_X)\theta_X^p}.
\end{multline*}
By construction, near each $p \in M_0(X)$, $\theta_X=\theta_X^p$, and
the negative sign appears because we regard each $S_{\epsilon}^p$
oriented as the boundary of $B_{\epsilon}^p$, which form the exterior of
$M \setminus \cup_p B_{\epsilon}^p$.
Recall that
$$
\frac 1{(d-\iota_X)\theta_X^p} =
-(\iota_X \theta_X^p)^{-1} - (\iota_X \theta_X^p)^{-2} d\theta_X^p
- (\iota_X \theta_X^p)^{-3} (d\theta_X^p)^2 - \dots
- (\iota_X \theta_X^p)^{-n/2-1} (d\theta_X^p)^{n/2},
$$
and so
$$
\frac 1{(d-\iota_X)\theta_X^p} \biggr|_{S_{\epsilon}^p}
= -\epsilon^{-1} - \epsilon^{-2} d\theta_X^p
- \epsilon^{-3} (d\theta_X^p)^2 - \dots
- \epsilon^{-n/2} (d\theta_X^p)^{n/2-1}.
$$
Therefore,
$$
\frac{\theta_X^p \wedge \alpha(X)}{(d-\iota_X)\theta_X^p}
\biggr|_{S_{\epsilon}^p} =
- \epsilon^{-n/2} \alpha(X)_{[0]} \cdot \theta_X^p \wedge (d\theta_X^p)^{n/2-1}
- \sum_{j=1}^{n/2} \epsilon^{j-n/2}
\theta_X^p \wedge \alpha(X)_{[2j]} \wedge (d\theta_X^p)^{n/2-1-j}.
$$
When we integrate over $S_{\epsilon}^p$ and let $\epsilon \to 0+$
only the first term survives:
\begin{multline*}
- \lim_{\epsilon \to 0+} \int_{S_{\epsilon}^p}
\frac{\theta_X^p \wedge \alpha(X)}{(d-\iota_X)\theta_X^p}
= \lim_{\epsilon \to 0+} \int_{S_{\epsilon}^p} 
\epsilon^{-n/2} \alpha(X)_{[0]} \cdot \theta_X^p \wedge (d\theta_X^p)^{n/2-1}\\
= \alpha(X)_{[0]}(p) \cdot \lim_{\epsilon \to 0+} \int_{S_{\epsilon}^p} 
\epsilon^{-n/2} \cdot \theta_X^p \wedge (d\theta_X^p)^{n/2-1}  \\
= \alpha(X)_{[0]}(p) \cdot \lim_{\epsilon \to 0+} \int_{B_{\epsilon}^p} 
\epsilon^{-n/2} \cdot (d\theta_X^p)^{n/2}.
\end{multline*}
But
$$
d\theta_X^p = -2 c_1^{-1} \cdot dx_1 \wedge dx_2 - \dots -
2 c_{n/2}^{-1} \cdot dx_{n-1} \wedge dx_n
$$
and
$$
(d\theta_X^p)^{n/2} = (-2)^{n/2} (n/2)! (c_1 \dots c_{n/2})^{-1} \cdot
dx_1 \wedge \dots \wedge dx_n.
$$
Since the volume of the $n$-dimensional unit ball (for $n$-even) equals
$\pi^{n/2}/(n/2)!$, we obtain the theorem.
\qed

\separate

\subsection{An Alternative Proof of the Integral Localization Formula}

In this subsection we give an alternative proof of the integral
localization formula (theorem \ref{BV-loc}).

Let $\pi: T^*M \twoheadrightarrow M$ denote the projection map.
We regard $M$ as a submanifold of $T^*M$ via the zero section inclusion.
Recall the canonical equivariantly closed form $\mu + \sigma$ on $T^*M$
constructed in Subsection \ref{sigma}. We can consider the form
\begin{equation}  \label{talpha}
\tilde\alpha(X) = e^{\mu(X) + \sigma} \wedge \pi^* \bigl( \alpha(X) \bigr),
\qquad X \in \g g.
\end{equation}
It is a form on $T^*M$ satisfying $(d-\iota_X) \tilde\alpha(X)=0$
for the reason that it is ``assembled'' from forms satisfying this equation.
Moreover, its restriction to $M$ is just $\alpha(X)$.
We will see later that, in a way, $\tilde\alpha$ is the most natural
equivariant extension of $\alpha$ to $T^*M$.

Fix a $G$-invariant Riemannian metric on $T^*M$ and let $\|.\|$ be the norm.
It induces a vector bundle diffeomorphism $TM \to T^*M$ and a norm
on $TM$ which we also denote by $\|.\|$.
Let $s: M \to T^*M$ denote the image of the vector field
$L_X$ under this diffeomorphism. Observe that
$$
\mu(X)(s(x)) = - \|L_X\|^2 \le 0, \qquad \text{for all $x \in M$},
$$
and that the set
$$
\{ x \in M;\, \mu(X)(s(x))=0 \} = \{ x \in M;\, L_X \bigr|_x =0 \} = M_0(X)
$$
is precisely the set of zeroes of $L_X$.

Pick an $R>0$ and let a parameter $t$ increase from 0 to $R$.
Then the family of sections $\{ts\}_{t\in[0,R]}$
provides a deformation of the initial cycle $M$ into $Rs$.
In other words, there exists a chain $C_R$ in $T^*M$
of dimension $(\dim M +1)$ such that
\begin{equation}  \label{deformation}
\partial C_R = M - Rs.
\end{equation}
Note that the support of $C_R$, $|C_R|$, lies inside the set
\begin{equation}  \label{<0}
\{ \xi \in T^*M;\, \mu(X)(\xi) = -\langle \xi, X_M \rangle \le 0 \}.
\end{equation}
This will ensure good behavior of the integrand when we let $R$
tend to infinity.

Observe that if we let $R \to \infty$, then we obtain a conic Borel-Moore
chain $C_{\infty}$ such that
$$
\partial C_{\infty} = M - \sum_{p \in M_0(X)} T^*_pM,
$$
where each cotangent space $T^*_pM$ is provided with
an appropriate orientation.
But we will not use this observation in the proof because one cannot
interchange the order of integration and taking limit as $R \to \infty$
(indeed, $\tilde\alpha \bigr|_p =0$).

\begin{lem}
$\bigl( d \tilde\alpha(X)_{[\dim M]} \bigr) \Bigr|_{|C_R|} =0$.
\end{lem}

\begin{rem}
The form $\tilde\alpha(X)_{[\dim M]}$ itself need not be closed; this lemma
only says that $\tilde\alpha(X)_{[\dim M]}$ becomes closed when restricted to
$|C_R|$.

This lemma is the key part of the proof and it is the only place
where we use that the form $\alpha$ is equivariantly closed and that
our Riemannian metric is $G$-invariant.
\end{rem}

\pf
Since the form $\tilde\alpha$ is equivariantly closed,
$(d-\iota_X) \tilde\alpha(X)=0$.
By abuse of notation, let $L_X$ denote the vector field generated by $X$
on $T^*M$. Since we chose our Riemannian metric to be $G$-invariant,
the vector field $L_X$ is tangent to $|C_R|$.
Therefore,
$$
(d-\iota_X) \bigl( \tilde\alpha(X) \bigr|_{|C_R|} \bigr) =
\bigl( (d-\iota_X) \tilde\alpha(X) \bigr) \bigr|_{|C_R|} =0.
$$
But the highest differential form degree component of 
$\tilde\alpha(X) \bigr|_{|C_R|}$ is $\tilde\alpha(X)_{[\dim M]}$,
hence $d \bigl( \tilde\alpha(X) \bigr|_{|C_R|} \bigr) =0$.
\qed

Therefore,
\begin{equation} \label{R}
\int_M \alpha(X) = \int_M \tilde\alpha(X) = \int_{Rs} \tilde\alpha(X),
\qquad \forall R \in \BB R.
\end{equation}
For each $p \in M_0(X)$ we choose a small neighborhood $U_p$ of $p$.
Since $M_0(X)$ is discrete, we can make these neighborhoods
$U_p$'s small enough so that they do not overlap.
In (\ref{R}) we let $R$ tend to infinity.
The $e^{\mu(X)+\sigma}$ component in $\tilde\alpha$
(defined by (\ref{talpha})) provides
exponential decay of $\tilde\alpha(X)$ on the set (\ref{<0}), and we obtain:
\begin{multline*}
\int_M \alpha(X) = \lim_{R \to \infty} \int_{Rs} \tilde\alpha(X)
= \lim_{R \to \infty}
\int_{Rs \cap \bigcup_{p \in M_0(X)} \pi^{-1} U_p} \tilde\alpha(X)  \\
= \sum_{p \in M_0(X)} \lim_{R \to \infty}
\int_{Rs \cap \pi^{-1} U_p} \tilde\alpha(X).
\end{multline*}
That is, when we integrate over $Rs$, only the covectors whose basepoint
lies in some $U_p$ (i.e. ``near'' $M_0(X)$) count.
Thus we see that the integral is localized around the
set of zeroes $M_0(X)$.
It remains to calculate the contribution of each point $p$:
$$
\lim_{R \to \infty} \int_{Rs \cap \pi^{-1} U_p} \tilde\alpha(X).
$$

As in the original proof, we can assume that $G$ is abelian.
Let $p \in M_0(X)$, then $p$ is fixed by $G$, $G$ acts linearly on $T_pM$,
and there is a positively oriented orthonormal basis $\{e_1,\dots,e_n\}$
of $T_pM$ such that $L(X,p)$ has a block-diagonal form,
and each block is a matrix of the type
$\begin{pmatrix} 0 & -c_j  \\ c_j & 0 \end{pmatrix}$ for some $c_j \in \BB R$,
then ${\det}^{1/2}(L(X,p)) = c_1 \dots c_{n/2}$.
Making $U_p$ smaller and using the exponential map $\exp_p: T_pM \to M$
we can construct local coordinates $x_1, \dots,x_n$ on $U_p$
centered at $p$ such that
$$
\frac{\partial}{\partial x_1}(p) = e_1,\dots,
\frac{\partial}{\partial x_n}(p) = e_n
$$
and the vector field $L_X$ is expressed by
$$
L_X = c_1 \Bigl( x_2 \frac{\partial}{\partial x_1} -
x_1 \frac{\partial}{\partial x_2} \Bigr) + \dots +
c_{n/2} \Bigl( x_n \frac{\partial}{\partial x_{n-1}} -
x_{n-1} \frac{\partial}{\partial x_n} \Bigr).
$$
Making $U_p$ smaller once more, we can assume that
the coordinate functions $(x_1,\dots, x_n)$ send $U_p$ diffeomorphically onto
a cube
$$
\{ (t_1,\dots,t_n);\, |t_1| \le \epsilon, \dots, |t_n| \le \epsilon \}
$$
for some $\epsilon>0$.

Recall that
$$
\Bigl\{ \frac {\partial}{\partial x_1} = e_1, \dots,
\frac {\partial}{\partial x_n} =e_n \Bigr\}
$$
is an orthonormal basis of $T_pM$.
Hence the Riemannian metric around each $p$ can be
written in terms of the frame on $U_p$
$$
\Bigl( \frac {\partial}{\partial x_1}, \dots,
\frac {\partial}{\partial x_n} \Bigr)
$$
as $I + E_p(x)$, where $I$ is the $n \times n$ identity matrix
and $E_p(x)$ is a symmetric $n \times n$ matrix which depends on $x \in M$
and vanishes at the point $p$.

Let $(x_1,\dots, x_n;\xi_1,\dots,\xi_n)$ be the system of coordinates on
$T^*U_p$ associated to the coordinate system $(x_1,\dots, x_n)$
so that every element of $T^*U_p$ is uniquely expressed as
$(x_1,\dots, x_n; \xi_1 dx_{p,1}+\dots+\xi_n dx_n)$.
In these coordinates,
$$
s = \bigl( x_1,\dots, x_n, c_1(x_2 dx_1 - x_1 dx_2) + \dots
+ c_{n/2} (x_n dx_{n-1} - x_{n-1} dx_n) + O(x^2) \bigr),
$$
where $O(x^2)$ denotes terms of quadratic order in $x_1,\dots, x_n$.
We also have:
$$
\sigma = d\xi_1 \wedge dx_1 + \dots + d\xi_n \wedge dx_n,
$$
and
$$
\mu(X) : \: (x_1,\dots, x_n;\xi_1,\dots,\xi_n)  \\
\mapsto
- c_1(x_2 \xi_1 - x_1 \xi_2) - \dots - c_{n/2}( x_n \xi_{n-1} - x_{n-1} \xi_n).
$$

Next we introduce a coordinate system for $Rs \cap \pi^{-1} U_p$ by taking
$$
(y_1=x_1 \circ \pi, \dots, y_n = x_n \circ \pi);
$$
it is positively oriented. In these coordinates
$$
dx_1=dy_1, \dots, dx_n=dy_n,
$$
\begin{multline*}
d\xi_1 = Rc_1 dy_2 + R \cdot O(y),  \qquad
d\xi_2 = -Rc_1 dy_1 + R \cdot O(y), \dots,  \\
d\xi_{n-1} = Rc_{n/2} dy_n + R \cdot O(y), \quad
d\xi_n = -Rc_{n/2} dy_{n-1} + R \cdot O(y),
\end{multline*}
$$
\mu(X) \bigr|_{Rs} = -R \cdot \bigl( c_1^2 (y_1^2 + y_2^2) + \dots
+ c_{n/2}^2 (y_{n-1}^2 + y_n^2) \bigr) + R \cdot O(y^3).
$$
Thus $\sigma \bigr|_{Rs}$ in these coordinates becomes
$$
\sigma \bigr|_{Rs} = -2Rc_1 dy_1 \wedge dy_2 - \dots
-2R c_{n/2} dy_{n-1} \wedge dy_n + R \cdot O(y).
$$
Hence
$$
(e^{\sigma})_{[n]} \bigr|_{Rs} =
(-2R)^{n/2} c_1 \dots c_{n/2} dy_1 \wedge \dots \wedge dy_n
+ R^{n/2} \cdot O(y).
$$
On the other hand, the expression of $\pi^* \alpha(X) \bigr|_{Rs}$ in these
coordinates does not depend on $R$. We have:
$$
\tilde\alpha(X)_{[n]}
= \bigl( e^{\mu(X)+\sigma} \bigr)_{[n]} \alpha(X)_{[0]} +
\bigl( e^{\mu(X)+\sigma} \bigr)_{[n-2]} \alpha(X)_{[2]}+ \dots +
\bigl( e^{\mu(X)+\sigma} \bigr)_{[0]} \alpha(X)_{[n]}.
$$
Hence
\begin{multline*}
\tilde\alpha(X)_{[n]} \bigr|_{Rs}
= \bigl( e^{\mu(X)+\sigma} \bigr)_{[n]} \alpha(X)_{[0]} \bigr|_{Rs}  \\
+ e^{\mu(X)}
\bigl( \text{terms depending on $R$ through powers of $R$
strictly less than $n/2$} \bigr)  \\
= e^{\mu(X)} (e^{\sigma})_{[n]}
\alpha(X)_{[0]} \bigr|_{Rs} + e^{\mu(X)}o(R^{n/2}) = \\
e^{- R \bigl( c_1^2 (y_1^2 + y_2^2) + \dots
+ c_{n/2}^2 (y_{n-1}^2 + y_n^2) + O(y^3) \bigr)}
\cdot \bigl( (e^{\sigma})_{[n]} \alpha(X)_{[0]} \bigr|_{Rs}
+ o(R^{n/2}) \bigr).
\end{multline*}

Recall that the integration takes place over $U_p$,
hence the range of the variables is
$|y_1| \le \epsilon, \dots, |y_n| \le \epsilon$.
We perform a change of coordinates
$$
z_1 = \sqrt{R} y_1, \dots, z_n = \sqrt{R} y_n, \qquad
|z_1| \le \epsilon \sqrt{R}, \dots, |z_n| \le \epsilon\sqrt{R}.
$$
In these new coordinates,
$$
\tilde\alpha(X)_{[n]} \bigr|_{Rs} =  
e^{- \bigl( c_1^2 (z_1^2 + z_2^2) + \dots
+ c_{n/2}^2 (z_{n-1}^2 + z_n^2) + O(z^3)/O(\sqrt{R}) \bigr)}
\cdot \bigl( (e^{\sigma})_{[n]} \alpha(X)_{[0]} \bigr|_{Rs}
+ 1/O(\sqrt{R}) \bigr),
$$
and
$$
(e^{\sigma})_{[n]} \bigr|_{Rs} =
(-2)^{n/2} c_1 \dots c_{n/2} dz_1 \wedge \dots \wedge dz_n
+ O(z)/O(\sqrt{R}).
$$
As $R$ tends to infinity, the integrand tends pointwise to
$$
(-2)^{n/2} c_1 \dots c_{n/2}
e^{- c_1^2 (z_1^2 + z_2^2) - \dots - c_{n/2}^2 (z_{n-1}^2 + z_n^2)}
\cdot \alpha(X)_{[0]}(p) dz_1 \wedge \dots \wedge dz_n.
$$
Hence by the Lebesgue dominated convergence theorem,
\begin{multline*}
\lim_{R \to \infty}
\int_{Rs \cap \pi^{-1} U_p} \tilde\alpha(X)  \\
= (-2)^{n/2} c_1 \dots c_{n/2} \alpha(X)_{[0]}(p) \cdot
\int_{\BB R^n}  e^{- c_1^2 (z_1^2 + z_2^2) - \dots
- c_{n/2}^2 (z_{n-1}^2 + z_n^2)}
dz_1 \wedge \dots \wedge dz_n  \\
= (-2\pi)^{n/2} \frac {\alpha(X)_{[0]}(p)}{c_1\dots c_{n/2}}
= (-2\pi)^{n/2} \frac {\alpha(X)_{[0]}(p)}{{\det}^{1/2}(L(X,p))}.
\end{multline*}
This finishes our proof of Theorem \ref{BV-loc}.

\separate

\subsection{The General Integral Localization Formula}

In this subsection we give a statement of the general localization
formula.

\begin{prop}
Fix $X \in \mathfrak{g}$, and let $M_0(X)$ be the set of zeroes of the
vector field $L_X$ in $M$. Then $M_0(X)$ is a smooth submanifold of $M$
which may have several components of different dimension.
The normal bundle ${\cal N}$ of $M_0(X)$ in $M$ is an orientable vector
bundle with even-dimensional fibers.
\end{prop}

\pf
Fix a $G$-invariant Riemannian structure on $M$. Let $p \in M$ be a zero of
$L_X$, then the one parameter group $\exp (tX)$, $t \in \BB R$, acts on $T_pM$.
Note that a tangent $v \in T_pM$ is fixed by $\exp (tX)$ if and only if
the geodesic tangent $v$ is contained in $M_0(X)$.
Since the exponential map $\exp_p: T_pM \to M$ is a diffeomorphism near $p$,
this proves that $M_0(X)$ is a smooth submanifold of $M$.

The normal bundle ${\cal N} \to M_0(X)$ is defined fiberwise by setting
${\cal N}_p$ to be the orthogonal complement of $T_pM_0(X)$ inside $T_pM$,
for all $p \in M_0(X)$.
The linear transformation $L(X,p)$ preserves the decomposition
$T_pM = {\cal N}_p \oplus T_pM_0(X)$, it is identically zero on $T_pM_0(X)$
and invertible anti-symmetric on ${\cal N}_p$. This proves that the fibers
of ${\cal N}$ are even-dimensional.
The positive orientation of fibers ${\cal N}_p$ can be determined by requiring
$$
{\det}^{1/2} (L(X,p) \bigr|_{{\cal N}_p}) >0.
$$
\qed

Let $\operatorname{rk}({\cal N}): M_0(X) \to \BB Z$ be the locally constant
function on $M_0(X)$ which gives the codimension of each component.
Since $T_pM = {\cal N}_p \oplus T_pM_0(X)$ and $M$, ${\cal N}$ are oriented,
we get a canonical orientation on $M_0(X)$.
Thus a basis of $T_pM_0(X)$ is positive if an only if, when completed by
a positive basis of ${\cal N}_p$, it becomes a positive basis of $T_pM$.
Let ${\mathfrak g}_0(X) = \{ Y \in \mathfrak{g} ;\: [Y,X]=0 \}$
be the centralizer of $X \in \mathfrak{g}$, and let $G_0(X)$ be the connected
Lie subgroup of $G$ whose Lie algebra is ${\mathfrak g}_0(X)$.
Then $G_0(X)$ preserves the manifold $M_0(X)$ and acts on the normal bundle
${\cal N}$. We choose a $G_0(X)$-invariant Euclidean structure and a
$G_0(X)$-invariant metric connection $\nabla^{\cal N}$ on ${\cal N}$.
From this we get an equivariant curvature form
$$
\Theta^{\cal N}_{eq}(Y) = \Theta^{\cal N} + \mu^{\cal N}(Y),
\qquad \mu^{\cal N}(Y) = L^{\cal N}_Y - \nabla^{\cal N}_Y,
\qquad Y \in {\mathfrak g}_0(X).
$$
In particular, since the vector field $L_X$ vanishes on $M_0(X)$, we have
$$
\mu^{\cal N}(X) = L^{\cal N}_X
\qquad \text{and} \qquad
\mu^{\cal N}(X) \bigr|_p = L(X,p) \bigr|_{{\cal N}_p}
$$
which is an invertible transformation on ${\cal N}_p$.
We have already selected an orientation on ${\cal N}$, hence we can consider
the $G_0(X)$-equivariant Euler class associated to ${\cal N}$:
$$
\chi_{eq}({\cal N})(Y) = \operatorname{Pf}(-\Theta^{\cal N}_{eq}(Y))
= {\det}^{1/2}(-\Theta^{\cal N}_{eq}(Y)).
$$
Since
$$
\chi_{eq}({\cal N})(Y)_{[0]} = {\det}^{1/2}(-\mu^{\cal N}(Y)) \ne 0
$$
for $Y$ sufficiently close to $X$, the equivariant differential form
$\chi_{eq}({\cal N})(Y)$ is invertible for those $Y$.
We can now state the general integral localization formula.

\begin{thm} [Berline-Vergne]
Let $G$ be a compact group with Lie algebra $\mathfrak{g}$ acting on a
compact oriented manifold $M$, and let
$\alpha: \mathfrak{g} \to \Omega^*(M)$
be a ${\cal C}^{\infty}$ map such that $d_{eq} \alpha=0$.
(We do not require $\alpha$ to be an equivariant form.)
Pick an $X \in \mathfrak{g}$, let $M_0(X)$ be the set of zeroes of
$L_X$ in $M$, and let ${\cal N}$ be the normal bundle of $M_0(X)$ in $M$.
Then for $Y$ in the centralizer of $X$ in $\mathfrak{g}$ and sufficiently
close to $X$ (so that $\chi_{eq}({\cal N})(Y)$ is invertible), we have
$$
\int_M \alpha(Y) = (2\pi)^{\operatorname{rk}({\cal N})/2} \int_{M_0(X)}
\frac {\alpha(Y)}{\chi_{eq}({\cal N})(Y)},
$$
where $\chi_{eq}({\cal N})$ is the equivariant Euler form of the normal bundle.
In particular,
$$
\int_M \alpha(X) = (-2\pi)^{\operatorname{rk}({\cal N})/2} \int_{M_0(X)}
\frac {\alpha(X)}{{\det}^{1/2}(L_X^{{\cal N}}+ \Theta^{\cal N})}.
$$
\end{thm}

The proof of this localization formula is similar to the proof of the
isolated zeroes case, but more involved.

\separate

\subsection{Duistermaat-Heckman's Exact Stationary Phase Approximation}

One important application of the integral localization formula is the
``exact stationary phase approximation'' of Duistermaat-Heckman \cite{DH}.

\begin{thm} [Duistermaat-Heckman, 1982]  \label{DH-thm}
Let $(G \lefttorightarrow M, \omega, \mu)$ be a Hamiltonian system
with $M$ compact.
If $X \in \mathfrak{g}$ is such that $M_0(X)$ is finite, then
$$
\int_M e^{i\mu(X)} \, \frac{\omega^{\dim M/2}}{(\dim M/2)!}
= (2\pi i)^{\dim M/2} \sum_{p \in M_0(X)}
\frac {e^{i\mu(X)(p)}}{{\det}^{1/2}(L(X,p))}.
$$
(The square root of ${\det}(L(X,p))$ is computed with respect to the
canonical orientation of $T_pM$ given by $\omega^{\dim M/2}$.)
\end{thm}

\pf
Recall that the equivariant symplectic form $\tilde\omega = \mu + \omega$
is an equivariantly closed form on $M$, and so is $e^{i \tilde\omega(X)}$.
Since
$$
(e^{i \tilde\omega(X)})_{[\dim M]} = (e^{i\mu(X)} e^{i\omega})_{[\dim M]}
= e^{i\mu(X)} \, \frac{\omega^{\dim M/2}}{(\dim M/2)!},
$$
we have
$$
\int_M e^{i\mu(X)} \, \frac{\omega^{\dim M/2}}{(\dim M/2)!} =
\int_M e^{i \tilde\omega(X)},
$$
and the result follows from the integral localization formula.
\qed

This result can be interpreted as follows.
The Liouville form $\omega^{\dim M/2}/(\dim M/2)!$ determines a measure
on $M$. Then the moment map $\mu : M \to \mathfrak{g}^*$ determines the
push-forward of this measure $d\nu$ called the Duistermaat-Heckman measure:
$$
\int_{\mathfrak{g}^*} f \,d\nu =_{def}
\int_M (f \circ \mu) \cdot \frac{\omega^{\dim M/2}}{(\dim M/2)!},
\qquad f \in {\cal C}^{\infty}_0(\mathfrak{g}^*).
$$
Clearly, the support of this measure lies in $\mu(M) \subset \mathfrak{g}^*$.
In this context, the expression
$$
\int_M e^{i\mu(X)} \, \frac{\omega^{\dim M/2}}{(\dim M/2)!}
= \int_{\mathfrak{g}^*} e^{i \langle \,\cdot\, , X \rangle} \,d\nu
$$
is nothing but the Fourier transform of the Duistermaat-Heckman measure.
Since the Fourier transform can be inverted, the measure $d\nu$ is uniquely
determined by its Fourier transform.

Suppose now that the group $G$ is a torus, and introduce the set of
fixed points $M^G = \{p \in M;\: G \cdot p =p \}$.
This set need not be finite, but it has finitely many connected components,
and the moment map $\mu$ is constant on each connected component.
Hence the set $\mu(M^G)$ is finite.

Then we have the following convexity theorem
(for reference see \cite{A}, \cite{GS2} and a very nice exposition \cite{GGK}):

\begin{thm} [Atiyah, Guillemin and Sternberg, 1982]
If $M$ is connected, the image $\mu(M)$ is a convex polytope in
the vector space $\mathfrak{g}^*$.
Moreover, $\mu(M)$ is the convex hull of $\mu(M^G)$.
\end{thm}
(Note that it is possible for some $\mu(p)$, $p \in M^G$, to lie in the
interior of $\mu(M)$.)

Let $\mu(M)^{reg}$ denote the set of regular values of $\mu$ in $\mu(M)$.
The connected components of $\mu(M)^{reg}$ are known to be open convex
polytopes.
Inverting the Fourier transform one can prove the following result.

\begin{thm} [Duistermaat-Heckman]
Let $(G \lefttorightarrow M, \omega, \mu)$ be a Hamiltonian system where
the group $G$ is a torus acting on a compact manifold $M$.
Then the Duistermaat-Heckman measure $d\nu$ on $\mathfrak{g}^*$
is locally polynomial relatively to the Lebesgue measure
on $\mathfrak{g}^*$. That is, on each connected component of 
$\mu(M)^{reg}$ the measure $d\nu$ is polynomial.
\end{thm}

\separate

\subsection{Kirillov's Character Formula}

In this subsection we give a proof due to N.~Berline and M.~Vergne
of Kirillov's character formula using the integral localization formula.

Let $G$ be a connected compact Lie groups with maximal torus $T$.
The exponential map $\exp : \mathfrak{t} \to T$ is not only
a surjective homomorphism, but also locally bijective,
hence a covering homomorphism,
\begin{equation}
\exp : \mathfrak{t} /\Lambda' \ \tilde{\longrightarrow} \ T,
\qquad \Lambda' = \{ X \in \mathfrak{t} ;\: \exp X = e \}.
\end{equation}
That makes $\Lambda' \subset \mathfrak{t}$ a discrete, cocompact subgroup.
In other words, $L'$ is a lattice, the so-called {\em unit lattice}.
Let $\widehat T$ denote the group of characters, i.e. the group of
homomorphisms from $T$ to the unit circle $S^1 = \{ z \in \BB C ;\: |z|=1 \}$.
Then
\begin{equation*}
\Lambda =_{def} \{ \lambda \in i\mathfrak{t}^* ;\:
\langle \lambda, \Lambda' \rangle \subset 2\pi i \BB Z \}
\ \ \overset{\sim}{\longrightarrow} \ \ \widehat T,
\qquad \Lambda \ni \lambda \mapsto e^{\lambda} \in \widehat T,
\end{equation*}
with $e^{\lambda} : T \to \{ z \in \BB C ;\: |z|=1 \}$
defined by $e^{\lambda}(\exp X) = e^{\langle \lambda, X \rangle}$,
for any $X \in \mathfrak{t}$.
One calls $\Lambda \subset i\mathfrak{t}^*$ the {\em weight lattice};
except for the factor $2\pi i$ in its definition, it is the lattice dual
of the unit lattice $\Lambda' \subset \mathfrak{t}$.

We have a notion of roots $\Phi \subset \Lambda$; these are the
weights of the adjoint representation of $G$.
The Weyl group $W$ acts on $\mathfrak{t}$, $i\mathfrak{t}^*$,
$\Lambda$ and $\Phi$. For an element $w \in W$, we denote by
$\epsilon(w) = \pm 1$ the determinant of the transformation on $\mathfrak{t}$.
Select a positive root system $\Phi^+ \subset \Phi$.
Fixing a $G$-invariant positive definite metric $( \,\cdot\,,\,\cdot\,)$ on
$i\mathfrak{g}^*$, we get dominant weights:
$$
\Lambda^+ = \{ \lambda \in \Lambda ;\:
(\lambda, \alpha) \ge 0 \: \forall \alpha \in \Phi^+ \}.
$$
The highest weight theorem asserts that there is a bijection between
the dominant elements $\lambda$ of $\Lambda$ and irreducible
representations $\pi_{\lambda}$ of $G$ so that $\lambda$ occurs as the
highest weight of $\pi_{\lambda}$.

Define an element
$$
\rho = \frac 12 \sum_{\alpha \in \Phi^+} \alpha \quad \in i\mathfrak{t}^*
$$
(evidently $2\rho$ is a weight, but $\rho$ itself may or may not be a weight).
We can state the Weyl character formula:

\begin{thm} [Weyl]
Let $\lambda \in \Lambda^+$ be a dominant weight, and let
$(\pi_{\lambda}, V_{\lambda})$ be the irreducible representation of $G$
with highest weight $\lambda$. Then, for all $X \in \mathfrak{t}$,
$$
\tr \pi_{\lambda}(\exp X) = \frac
{\sum_{w \in W} \epsilon(w) e^{\langle w (\lambda+\rho), X \rangle}}
{\prod_{\alpha \in \Phi^+} (e^{\langle \alpha, X \rangle/2} -
e^{-\langle \alpha, X \rangle/2})}.
$$
\end{thm}

The decomposition
$$
\mathfrak{g} = \mathfrak{t} \oplus (\text{root spaces} \cap \mathfrak{g})
= \mathfrak{t} \oplus [\mathfrak{t}, \mathfrak{g}]
$$
provides a splitting
$$
i \mathfrak{g}^* = i \mathfrak{t}^* \oplus i [\mathfrak{t}, \mathfrak{g}]^*
$$
which allows us to think of $\lambda$ and $\rho$ as elements of
$i \mathfrak{g}^*$.
Let ${\cal O}_{-i(\lambda+\rho)} \subset \mathfrak{g}^*$ denote the coadjoint
$G$-orbit of $-i(\lambda+\rho)$.
Recall that ${\cal O}_{-i(\lambda+\rho)}$ has a canonical structure of a
Hamiltonian system
$(G \lefttorightarrow {\cal O}_{-i(\lambda+\rho)}, \sigma,\mu)$
given by Proposition \ref{symplectic-form},
with the moment map $\mu$ being the inclusion
${\cal O}_{-i(\lambda+\rho)} \hookrightarrow \mathfrak{g}^*$.

Kirillov's character formula involves the square root of a function 
$j_{\mathfrak{g}}$ on $\mathfrak{g}$ defined by
$$
j_{\mathfrak{g}}(X) = \det \exp_*(X) =
\det \biggl( \frac{\sinh (ad \, X/2)}{ad \, X/2} \biggr),
\qquad X \in \mathfrak{g}.
$$
If $X \in \mathfrak{t}$ we get
$$
j_{\mathfrak{g}}(X) = \prod_{\alpha \in \Phi}
\frac {e^{\langle \alpha, X \rangle/2} - e^{-\langle \alpha, X \rangle/2}}
{\langle \alpha, X \rangle}
= \biggl( \prod_{\alpha \in \Phi^+}
\frac {e^{\langle \alpha, X \rangle/2} - e^{-\langle \alpha, X \rangle/2}}
{\langle \alpha, X \rangle} \biggr)^2.
$$
Chevalley's Theorem asserts that the restriction map on functions induced
by the inclusion $\mathfrak{t} \hookrightarrow \mathfrak{g}$ provides
isomorphisms of the following spaces of functions:
$$
{\cal C}^{\infty}(\mathfrak{g})^G \simeq {\cal C}^{\infty}(\mathfrak{t})^W,
\qquad
{\cal C}^{\omega}(\mathfrak{g})^G \simeq {\cal C}^{\omega}(\mathfrak{t})^W,
\qquad
\BB C [\mathfrak{g}]^G \simeq \BB C [\mathfrak{t}]^W.
$$
The analytic function $\prod_{\alpha \in \Phi^+}
\frac {e^{\langle \alpha, X \rangle/2} - e^{-\langle \alpha, X \rangle/2}}
{\langle \alpha, X \rangle}$ on $\mathfrak{t}$ is easily seen to be
$W$-invariant and hence by Chevalley's Theorem defines a $G$-invariant
function on all of $\mathfrak{g}$ which we denote by $j_{\mathfrak{g}}^{1/2}$.

\begin{thm} [Kirillov]
Let $\lambda \in \Lambda^+$ be a dominant weight, and let
$(\pi_{\lambda}, V_{\lambda})$ be the irreducible representation of $G$
with highest weight $\lambda$. Then, for all $X \in \mathfrak{g}$,
$$
j_{\mathfrak{g}}^{1/2}(X) \cdot \tr \pi_{\lambda}(\exp X) =
(2\pi)^{-n/2} \cdot
\int_{{\cal O}_{-i(\lambda+\rho)}} e^{i\mu(X)} \,dm,
$$
where $dm = \frac{\sigma^{n/2}}{(n/2)!}$ is the Liouville
measure on ${\cal O}_{-i(\lambda+\rho)}$ and
$n= \dim {\cal O}_{-i(\lambda+\rho)}$.
\end{thm}

The expression on the left hand side is a function on the Lie algebra
$\mathfrak{g}$ and is often called the character of
$(\pi_{\lambda}, V_{\lambda})$ on the Lie algebra.
The expression on the right hand side is the Fourier transform
of coadjoint orbit and denoted by $\widehat {\cal O}_{-i(\lambda+\rho)}$.

The proof uses the integral localization formula and essentially matches
the contribution from each zero with the corresponding term in the Weyl
character formula.
Both sides are $G$-invariant functions on $\mathfrak{g}$.
Thus it is sufficient to check the equality when $X \in \mathfrak{t}$.
Moreover, by continuity we can assume that $X$ is regular, i.e. that
$\langle \alpha, X \rangle \ne 0$ for every $\alpha \in \Phi$.
We expand the left hand side to get:
\begin{multline*}
j_{\mathfrak{g}}^{1/2}(X) \cdot \tr \pi_{\lambda}(\exp X)  \\
= \prod_{\alpha \in \Phi^+}
\frac {e^{\langle \alpha, X \rangle/2} - e^{-\langle \alpha, X \rangle/2}}
{\langle \alpha, X \rangle} \cdot 
\frac {\sum_{w \in W} \epsilon(w) e^{\langle w (\lambda+\rho), X \rangle}}
{\prod_{\alpha \in \Phi^+} (e^{\langle \alpha, X \rangle/2} -
e^{-\langle \alpha, X \rangle/2})}  \\
= \frac {\sum_{w \in W} \epsilon(w) e^{\langle w (\lambda+\rho), X \rangle}}
{\prod_{\alpha \in \Phi^+} \langle \alpha, X \rangle}
\end{multline*}
thus obtaining the Weyl formula for the character on the Lie algebra.
Kirillov's character formula will then follow from Corollary
\ref{Harish-Chandra} and an observation that the subset
$\Phi^+_{{-i(\lambda+\rho)}} \subset \Phi$ which will be defined later
coincides with $\Phi^+$.

Let $\mathfrak{g}_{\BB C} = \BB C \otimes \mathfrak{g}$ and
$\mathfrak{t}_{\BB C} = \BB C \otimes \mathfrak{t}$ be the complexifications
of $\mathfrak{g}$ and $\mathfrak{t}$ respectively. We have
$$
\mathfrak{g}_{\BB C} = \mathfrak{t}_{\BB C} \oplus
\bigoplus_{\alpha \in \Phi} (\mathfrak{g}_{\BB C})_{\alpha},
$$
where $(\mathfrak{g}_{\BB C})_{\alpha}$ is the root space corresponding to
$\alpha$:
$$
(\mathfrak{g}_{\BB C})_{\alpha} = \{ X \in \mathfrak{g}_{\BB C} ;\:
[H,X] = \langle \alpha, X \rangle \cdot X, \:
\forall H \in \mathfrak{t}_{\BB C} \}.
$$
Then, for each $\alpha \in \Phi$, $-\alpha$ is also in $\Phi$ and
$$
\dim (\mathfrak{g}_{\BB C})_{\alpha} = \dim (\mathfrak{g}_{\BB C})_{-\alpha}
= \dim [(\mathfrak{g}_{\BB C})_{\alpha}, (\mathfrak{g}_{\BB C})_{-\alpha}] =1,
\qquad 
[(\mathfrak{g}_{\BB C})_{\alpha}, (\mathfrak{g}_{\BB C})_{-\alpha}]
\subset \mathfrak{t}_{\BB C}.
$$
Furthermore, there exists a unique element
$H_{\alpha} \in
[(\mathfrak{g}_{\BB C})_{\alpha}, (\mathfrak{g}_{\BB C})_{-\alpha}]$
such that $\langle \alpha, H_{\alpha} \rangle =2$.
The Killing form $B(X,Y)= \tr (ad X \cdot ad Y)$ is positive definite
on the real span
$\oplus_{\alpha \in \Phi} \BB H_{\alpha} \subset i\mathfrak{t}$.
We introduce vectors
$$
\tilde H_{\alpha} = \frac{2 H_{\alpha}}{B(H_{\alpha}, H_{\alpha})}
\in i\mathfrak{t}, \qquad \alpha \in \Phi.
$$
These vectors are the unique vectors in $\mathfrak{t}_{\BB C}$
such that
$$
B(\tilde H_{\alpha}, H) = \langle \alpha, H \rangle,
\qquad \forall H \in \mathfrak{t}_{\BB C}.
$$

Every coadjoint orbit of a compact group $G$ has the form
${\cal O}_{\nu} = G \cdot \nu$, with $\nu \in \mathfrak{t}^*$.
In other words, every coadjoint orbit in $\mathfrak{g}^*$ meets
$\mathfrak{t}^*$.
Note that we are considering a slightly general situation than in the
setting of Kirillov's character formula and, for instance, do not require
$\nu$ to be a weight nor regular.
The orbit ${\cal O}_{\nu}$ can be identified with the homogeneous space
$G/G(\nu)$, where $G(\nu)$ is the stabilizer of $\nu$.
For $\nu \in \mathfrak{t}^*$, let $\Phi^+_{\nu}$ be the set of roots
$$
\Phi^+_{\nu} = \{ \alpha \in \Phi ;\: \langle \nu, i H_{\alpha} \rangle > 0 \},
$$
and set
$$
\tau_{\nu} =
\mathfrak{g} \cap \Bigl( \bigoplus_{\alpha \in \Phi^+_{\nu}}
(\mathfrak{g}_{\BB C})_{\alpha} \oplus  (\mathfrak{g}_{\BB C})_{-\alpha}
\Bigr).
$$
The space $\tau_{\nu}$ can be identified with the tangent space
$T_e (G/G(\nu))$ at the base point $e= G(\nu)$.

In particular, if $\langle \nu, i H_{\alpha} \rangle \ne 0$ for every
$\alpha \in \Phi$, then $G(\nu)=T$ and
$\tau_{\nu} = [\mathfrak{t}, \mathfrak{g}]$.
Such a point $\nu$ is called regular and its orbit has maximal dimension.

The restriction of minus the Killing form to $\tau_{\nu}$ is positive
definite and determines a $G(\nu)$-invariant inner product on
$T_e (G/G(\nu))$. Thus the homogeneous space $G/G(\nu)$ is a
Riemannian manifold, with a $G$-invariant metric which coincides with
$-B \bigr|_{\tau_{\nu}}$ at the point $e$.
Denote by $d\bar g$ the corresponding $G$-invariant measure on $G/G(\nu)$.
The Liouville measure $dm_{\nu}$ of ${\cal O}_{\nu}$ is another
$G$-invariant measure on $G/G(\nu)$, so the two measures $dm_{\nu}$ and
$d\bar g$ must be proportional. The following lemma provides the
proportionality coefficient.

\begin{lem}
Let $n = \dim {\cal O}_{\nu}$.
The Liouville measure $dm_{\nu} = \frac{\sigma_{\nu}^{n/2}}{(n/2)!}$
of ${\cal O}_{\nu}$ is given by the formula
$$
\int_{{\cal O}_{\nu}} \phi \,dm_{\nu} = \prod_{\alpha \in \Phi^+_{\nu}}
\langle \nu, i \tilde H_{\alpha} \rangle
\int_{G/G(\nu)} \phi(g\nu) \,d\bar g.
$$
\end{lem}

\pf
On $\mathfrak{g}_{\BB C} = \BB C \otimes \mathfrak{g}$ we have a
complex conjugation relative to $\mathfrak{g}$.
For each $\alpha \in \Phi^+_{\nu}$ we pick a element
$X_{\alpha} \in (\mathfrak{g}_{\BB C})_{\alpha}$, then it is easy to see that
$$
X_{-\alpha}=_{def} - \overline{X_{\alpha}} \in (\mathfrak{g}_{\BB C})_{-\alpha}
\qquad \text{and} \qquad
[X_{\alpha}, X_{-\alpha}] \in i \mathfrak{t}.
$$
Hence $[X_{\alpha}, X_{-\alpha}]$ is a real multiple of $H_{\alpha}$ and
the real span of $\{ X_{\alpha}, X_{-\alpha}, H_{\alpha} \}$
is a Lie algebra isomorphic to $\mathfrak{sl}(2,\BB R)$.
By assumption, $\langle \alpha, H_{\alpha} \rangle =2 >0$, hence
$[X_{\alpha}, X_{-\alpha}]$ is a positive real multiple of $H_{\alpha}$.
Rescaling $X_{\pm\alpha}$ we can assume that 
$$
[X_{\alpha}, X_{-\alpha}] = H_{\alpha}.
$$
Then it follows that $2B(X_{\alpha},X_{-\alpha}) = B(H_{\alpha}, H_{\alpha})$.

Set
$$
(\tau_{\nu})_{\alpha} =
\mathfrak{g} \cap \bigl( (\mathfrak{g}_{\BB C})_{\alpha}
\oplus  (\mathfrak{g}_{\BB C})_{-\alpha} \bigr),
$$
then the elements $\{ e_{\alpha}, f_{\alpha} \}$ defined by
$$
e_{\alpha} = \frac {X_{\alpha}- X_{-\alpha}}{B(H_{\alpha},H_{\alpha})^{1/2}}
\qquad \text{and} \qquad
f_{\alpha} = \frac {iX_{\alpha}+ iX_{-\alpha}}{B(H_{\alpha},H_{\alpha})^{1/2}}
$$
lie in $(\tau_{\nu})_{\alpha}$ and form an orthonormal basis of
$(\tau_{\nu})_{\alpha}$ relative to $-B$.
Moreover,
$$
[e_{\alpha}, f_{\alpha}]
= \frac {2i [X_{\alpha}, X_{-\alpha}]}{B(H_{\alpha},H_{\alpha})}
= i\tilde H_{\alpha}.
$$
Therefore,
$$
\sigma_{\nu}(L_{f_{\alpha}}, L_{e_{\alpha}}) = \nu([e_{\alpha}, f_{\alpha}])
= \langle \nu, i \tilde H_{\alpha} \rangle,
$$
and the result follows.
\qed

\begin{cor}  \label{f-e-basis}
The pairs of vectors $\{ f_{\alpha}, e_{\alpha} \}_{\alpha \in \Phi^+_{\nu}}$
form a positively oriented orthonormal basis of
$\tau_{\nu} \simeq T_e (G/G(\nu))$ relative to $-B$.
Moreover, for each $X \in \mathfrak t$, the linear transformation
$L(X,\nu)$ on $T_e (G/G(\nu))$ is given by
$$
L(X,\nu) f_{\alpha} = i\alpha(X) e_{\alpha}, \qquad
L(X,\nu) e_{\alpha} = -i\alpha(X) f_{\alpha},
$$
i.e. $L(X,\nu)$ is block diagonal with each block equal
$$
\begin{pmatrix} 0 & -i\alpha(X) \\ i\alpha(X) & 0 \end{pmatrix}.
$$
\end{cor}

The Fourier transform of the measure
$dm_{\nu} = \frac{\sigma_{\nu}^n}{(n/2)!}$ of ${\cal O}_{\nu}$ is a
$G$-invariant analytic function $\widehat {\cal O}_{\nu}(X)$ on $\mathfrak{g}$
defined by
$$
\widehat {\cal O}_{\nu}(X) = \int_{{\cal O}_{\nu}}
e^{i \langle l, X \rangle} \,dm_{\nu}(l),
\qquad X \in \mathfrak{g}.
$$
Clearly, it is completely determined by its restriction to $\mathfrak{t}$.

\begin{thm} [Harish-Chandra]  \label{H-Ch-thm}
Given $\nu \in \mathfrak{t}^*$, let $W_{\nu} = \{ w \in W ;\: w\nu = \nu \}$
be the stabilizer of $\nu$ in the Weyl group $W$.
For $X \in \mathfrak{t}$, $X$ regular, the Fourier transform
$\widehat {\cal O}_{\nu}$ is given by the formula
$$
\widehat {\cal O}_{\nu}(X) = (2\pi)^{n/2} \cdot
\sum_{w \in W/W_{\nu}} \frac {e^{i \langle w\nu, X \rangle}}
{\prod_{\alpha \in \Phi^+_{\nu}} \langle w\alpha, X \rangle}.
$$
\end{thm}

\pf
Since $X \in \mathfrak{t}$ is regular, the zero set of the vector field
generated by the action of $\exp tX$ on $\mathfrak{g}^*$ is the subspace
$\mathfrak{t}^*$ of $\mathfrak{g}^*$ which is fixed by the whole of $T$.
Thus the zero set of the vector field $L_X$ on ${\cal O}_{\nu}$ is
the finite set
$$
{\cal O}_{\nu} \cap \mathfrak{t}^* = \{ w \nu ;\: w \in W/W_{\nu} \}.
$$
By the integral localization formula or the
Duistermaat-Heckman Theorem \ref{DH-thm} it suffices
to compute ${\det}^{1/2} (L(X,w\nu))$.
From Corollary \ref{f-e-basis} we see that
${\det}^{1/2} (L(X,\nu)) = \prod_{\alpha \in \Phi^+_{\nu}} (i\alpha(X))$
and it follows that
${\det}^{1/2} (L(X,w\nu)) = \prod_{\alpha \in \Phi^+_{\nu}} (iw\alpha(X))$,
for all $w \in W$.
\qed

\begin{cor} [Harish-Chandra]  \label{Harish-Chandra}
If ${\cal O}_{\nu}$ is a regular orbit of the coadjoint representation,
then for $X \in \mathfrak{t}$, $X$ regular,
$$
\widehat {\cal O}_{\nu}(X) = (2\pi)^{n/2} \cdot
\frac {\sum_{w \in W} \epsilon(w) e^{i \langle w\nu, X \rangle}}
{\prod_{\alpha \in \Phi^+_{\nu}} \langle \alpha, X \rangle}.
$$
\end{cor}

\pf
If $\nu$ is regular, $\Phi^+_{\nu}$ is a positive root system in $\Phi$,
and up to a sign, the denominator
$$
{\det}^{1/2} (L(X,w\nu)) = \prod_{\alpha \in \Phi^+_{\nu}}
\langle w\alpha, iX \rangle
$$
is the same at each fixed point, and the sign is exactly $\epsilon(w)$.
\qed

\separate

\subsection{The Fourier Transform of Coadjoint Orbits of Reductive Groups}

If the group $G$ is not compact, the coadjoint orbits of $G$ in
$\mathfrak{g}^*$ are Hamiltonian systems and their Fourier transforms
may be defined as distributions (generalized functions) on $\mathfrak{g}$.
In this subsection we discuss an extension of Theorem \ref{H-Ch-thm} to the
case of real semisimple Lie groups.
Everything stated here is valid for the group $G$ equal the special linear
group $SL(2,\BB R)$ or $SU(p,q)$.

In this subsection only we assume that $G$ is a connected real semisimple
Lie group and possibly not compact, $\mathfrak{g}=Lie(G)$.
The coadjoint orbit ${\cal O}$ is a closed submanifold of $\mathfrak{g}$
if and only if it is semisimple, i.e. it is of the form $G \cdot \nu$ with
$\nu \in \mathfrak{g}^*$ corresponding to a semisimple element in
$\mathfrak{g}$ via a $G$-equivariant isomorphism
$\mathfrak{g} \simeq \mathfrak{g}^*$.
As before, $\sigma$ is the canonical symplectic form on ${\cal O}$,
$n = \dim {\cal O}$, $dm = \frac{\sigma^{n/2}}{(n/2)!}$ is the Liouville
measure, and the moment map $\mu: {\cal O} \hookrightarrow \mathfrak{g}^*$
is the inclusion. For an $X \in \mathfrak{g}$ the integral
$\int_{{\cal O}_{\nu}} e^{i\mu(X)} \,dm$ will not converge at all.
However, the Fourier transform $\widehat {\cal O}$ can be defined as
a distribution on $\mathfrak{g}$.
For every test function $\phi(X) \in {\cal C}^{\infty}_0(\mathfrak{g})$,
its Fourier transform
$$
\hat \phi(\xi) = \int_{\mathfrak{g}} e^{i\langle \xi, X \rangle}
\cdot \phi(X) \,dX,
\qquad \xi \in \mathfrak{g}^*,
$$
where $dX$ is a fixed Lebesgue measure on $\mathfrak{g}$, decays rapidly
as $\xi \to \infty$. The decay at infinity is so fast that the integral
$\int_{{\cal O}_{\nu}} \hat \phi(\xi) \,dm(\xi)$ converges absolutely.
By definition, the Fourier transform of the coadjoint orbit ${\cal O}$
is the distribution on $\mathfrak{g}$ such that
$$
\widehat {\cal O} : \qquad {\cal C}^{\infty}_0(\mathfrak{g}) \ni \phi
\quad \longmapsto \quad 
\int_{{\cal O}_{\nu}} \hat \phi(\xi) \,dm(\xi) =
\int_{{\cal O}_{\nu}} \biggl(
\int_{\mathfrak{g}} e^{i\langle \xi, X \rangle} \cdot \phi(X)
\,dX \biggr) \,dm(\xi).
$$

Let $K \subset G$ be a maximal compact subgroup and $T \subset K$
be a maximal torus in $K$; denote their Lie algebras respectively by
$\mathfrak{k}$ and $\mathfrak{k}$.
In general, $\mathfrak{t}_{\BB C}$ is a Cartan subalgebra in
$\mathfrak{k}_{\BB C}$, but may fail to be a Cartan subalgebra of
$\mathfrak{g}_{\BB C}$ being too small. In this subsection we assume that
$G$ and $K$ have {\em equal rank}, i.e.
$\mathfrak{t}_{\BB C} \subset \mathfrak{g}_{\BB C}$ {\em is} a Cartan
subalgebra. Examples of such groups $G$ are $SU(p,q)$ including an
important special case $SU(1,1) \simeq SL(2,\BB R)$.
We use the decomposition
$\mathfrak{g} = \mathfrak{t} \oplus [\mathfrak{t}, \mathfrak{g}]$
to identify $\mathfrak{t}^*$ as a linear subspace of $\mathfrak{g}^*$.

We will need the notions of Cartan involution on $\mathfrak{g}$
and Cartan decomposition of $\mathfrak{g}$.
The Killing form $B(X,Y)= \tr (ad X \cdot ad Y)$ is negative definite on
the Lie algebras of compact semisimple Lie groups, but never definite on
the Lie algebras of non-compact semisimple groups. A Lie algebra involution
$\theta: \mathfrak{g} \to \mathfrak{g}$ is called a {\em Cartan involution}
if the form $B_{\theta}(X,Y) =_{def} -B(X,\theta Y)$ is positive definite.
Then $\mathfrak{g}$ decomposes into the direct sum of $\pm 1$ eigenspaces
of $\theta$.
A general fact from the structure theory of semisimple Lie algebras asserts
that there always exists a Cartan involution $\theta$ such that the
$+1$ eigenspace is $\mathfrak{k}$. Let $\mathfrak{p}$ be the $-1$ eigenspace
of $\theta$. Thus we get the Cartan decomposition of $\mathfrak{g}$:
$\mathfrak{g} = \mathfrak{t} \oplus \mathfrak{p}$.

Let $\Phi = \Phi(\mathfrak{g}_{\BB C}, \mathfrak{t}_{\BB C})$ be the
root system of $K$, and set
$$
\Phi_{\mathfrak{p}} = \{ \alpha \in \Phi ;\:
(\mathfrak{g}_{\BB C})_{\alpha} \subset \mathfrak{p}_{\BB C} \};
$$
a root belonging to $\Phi_{\mathfrak{p}}$ is called non-compact.
Fix a $\nu \in \mathfrak{t}^*$, let
$\Phi^+_{\nu} = \{ \alpha \in \Phi ;\:
\langle \nu, i H_{\alpha} \rangle > 0 \}$, as before, and let $n_{\nu}$
be the number of non-compact roots contained in $\Phi^+_{\nu}$.
Define the set of regular elements in $\mathfrak{t}$ by
$$
\mathfrak{t}^{reg} =
\{ H \in \mathfrak{t} ;\: \langle \alpha, H \rangle \ne 0 \:
\forall \alpha \in \Phi \}.
$$
Then $G \cdot \mathfrak{t}^{reg}$ is an open set of $\mathfrak{g}$, but,
unless the group $G$ is compact, $G \cdot \mathfrak{t}^{reg}$ is never dense
in $\mathfrak{g}$.
The Fourier transform $\widehat {\cal O}$ is a distribution on $\mathfrak{g}$,
its restriction to $G \cdot \mathfrak{t}^{reg}$ is given by
integration against an analytic function and this function in turn is
completely determined by its restriction to $\mathfrak{t}^{reg}$.
The following result due to W.~Rossmann \cite{Ro} describes the restriction of
$\widehat {\cal O}$ to $G \cdot \mathfrak{t}^{reg}$.

\begin{thm} [Rossmann, 1978]
Let ${\cal O}$ is a closed orbit of the coadjoint representation of a connected
real semisimple Lie group $G$ such that $G$ and $K$ have equal rank.
Let $W = W(\mathfrak{k}_{\BB C}, \mathfrak{t}_{\BB C})$ be the compact
Weyl group. Then for $X \in \mathfrak{t}^{reg}$, we have the following results:
\begin{enumerate}
\item
If ${\cal O} \cap \mathfrak{t}^* = \varnothing$,
then $\widehat {\cal O} (X) =0$.

\item
If ${\cal O} = G \cdot \nu$ with $\nu \in \mathfrak{t}^*$, letting
$W_{\nu}$ be the stabilizer of $\nu$ in $W$,
$$
\widehat {\cal O} (X) = (-1)^{n_{\nu}}(2\pi)^{\dim {\cal O}/2} 
\sum_{w \in W/W_{\nu}} \frac {e^{i \langle w\nu, X \rangle}}
{\prod_{\alpha \in \Phi^+_{\nu}} \langle w\alpha, X \rangle}.
$$
\end{enumerate}
\end{thm}

In the non-compact setting we need to deal with two significant obstacles.
First of all, the group $G$ is not compact and the integral localization
formula is no longer applicable. This is dealt with by restricting to
the split group case and focusing all attention to the set
$G \cdot \mathfrak{t}^{reg}$. Roughly speaking, this is equivalent to
replacing $G$ with a compact subgroup $K$, so this result should not be
confused with a localization with respect to a non-compact group action.
Secondly, the manifold ${\cal O}$ is non-compact and extra care must be taken
to ensure convergence of the integrals and making sure there is no
contribution from infinity. This part is done by rewriting
everything as distributions which is a very general and
powerful method for dealing with non-converging integrals.

\separate

\subsection{Atiyah-Bott Abstract Localization Theorem}  \label{AB-loc}

In this subsection we follow the original M.~Atiyah and R.~Bott's paper
\cite{AB} and the book \cite{GS}.
Let $G$ be a compact group acting on a manifold $M$, and let $T \subset G$ be
a maximal torus.
Note that the equation (\ref{restriction_map}) applied to the inclusion
$T \hookrightarrow G$ gives a canonical map $H^*_G(M) \to H^*_T(M)$.
We will show later that this map is injective and its image is
$\bigl( H^*_T(M) \bigr)^W$.
Thus $H^*_G(M)$ can be recovered from $H^*_T(M)$, and
in this subsection we only work with a compact abelian group $T$.
Let $F$ be the set of points in $M$ fixed by $T$:
$$
F = \{ p \in M ;\: T \cdot p = p \}.
$$
The inclusion $i: F \hookrightarrow M$ induces a map on cohomology
$i^*: H^*_T(M) \to H^*_T(F)$. While this map need not be an isomorphism,
the Atiyah-Bott localization theorem asserts that $i^*$ is an
``isomorphism modulo torsion.''

Recall that $H^*_T(pt) = \BB C[\mathfrak{t}] = \BB C[X_1,\dots,X_l]$
is the polynomial algebra with each $X_j$ having degree $2$.
The map $M \to pt$ induces the cohomology map $H^*_T(pt) \to H^*_T(M)$ which
makes $ H^*_T(M)$ a $\BB C[\mathfrak{t}]$-module.
Let $l = \dim T$, then $\BB C[\mathfrak{t}] \simeq \BB C[X_1,\dots,X_l]$.

If $H$ is a $\BB C[X_1,\dots,X_l]$-module, we can make it a module over
the field of fractions $\BB C(X_1,\dots,X_l)$ by ``localizing'' it:
$$
H^{loc} =_{def} H \otimes_{\BB C[X_1,\dots,X_l]} \BB C(X_1,\dots,X_l).
$$
The passage from $H$ to $H^{loc}$ kills the torsion submodule $H^{tor}$:
$$
H^{tor} =_{def} \{ h \in H ;\: \text{$fh=0$ for some
$f \in \BB C[X_1,\dots,X_l]$, $f \ne 0$} \}.
$$
On the other hand, $H^{loc}$ keeps track of the ``free'' part of $H$.
We define the {\em rank} of $H$ to be the dimension of the vector space
$H^{loc}$ over $\BB C(X_1,\dots,X_l)$.

In general, given a $\BB C[X_1,\dots,X_l]$-module $H$, we can associate to it
a sheaf ${\cal H}$ of modules over $\BB C^l$. Here $\BB C^l$ is equipped with
Zariski topology so that the closed sets are the hypersurfaces
$V_f = \{ x \in \BB C^l ;\: f(x) = 0 \}$, $f \in \BB C[X_1,\dots,X_l]$,
and their finite intersections, and the open sets are finite unions
of $U_f = \BB C^l \setminus V_f = \{ x \in \BB C^l ;\: f(x) \ne 0 \}$.
Given an $f \in \BB C[X_1,\dots,X_l]$, we can form a new ring
$\BB C[X_1,\dots,X_l]_f$ by inverting $f$, i.e. $\BB C[X_1,\dots,X_l]_f$
consists of all rational fractions whose denominators are powers of $f$.
We also consider the corresponding module
$$
H_f = H \otimes_{\BB C[X_1,\dots,X_l]} \BB C[X_1,\dots,X_l]_f.
$$
By definition,
$$
{\cal H}(U_f) = H_f, \qquad
{\cal H}(U_{f_1} \cup \dots \cup U_{f_k}) = H \otimes_{\BB C[X_1,\dots,X_l]}
\BB C[X_1,\dots,X_l]_{f_1,\dots,f_k},
$$
where $\BB C[X_1,\dots,X_l]_{f_1,\dots,f_k}$ denotes the ring of fractions
whose denominators are powers of $f_1,\dots,f_k$.
Recall that in general the support of a sheaf ${\cal F}$ is the complement
of the largest open set $U$ such that ${\cal F} \bigr|_U =0$.
We define the {\em support} of a $\BB C[X_1,\dots,X_l]$-module $H$ to be the
support of the sheaf ${\cal H}$. Thus
$$
\supp H = \supp {\cal H} = \bigcap_f V_f \text{ over all
$f \in \BB C[X_1,\dots,X_l]$ with $f^j H=0$ for some $j \in \BB N$}.
$$
In other words,
$$
\supp H \subset V_f \quad \Longleftrightarrow \quad H_f =0
\quad \Longleftrightarrow \quad \text{$f^j H=0$ for some $j \in \BB N$.}
$$
Thus a free module has the whole space $\BB C^l$ as its support, while
the support of a torsion module is a proper subset of $\BB C^l$.

Here is another way of expressing $\supp H$. First we define
the annihilating ideal
$$
Ann(H) = \{ f \in \BB C[X_1,\dots,X_l];\: fH=0 \},
$$
then
$$
\supp H = Var(Ann(H)) =
\{ x \in \BB C^l ;\: f(x)=0 \: \forall f \in Ann(H) \}.
$$

\begin{lem}  \label{support-lemma}
If $H' \to H \to H''$ is an exact sequence of $\BB C[X_1,\dots,X_l]$-modules,
then $H'_f \to H_f \to H''_f$ is an exact sequence of
$\BB C[X_1,\dots,X_l]_f$-modules, for any $f \in \BB C[X_1,\dots,X_l]$.
Hence we get an exact sequence of the corresponding sheaves of modules
$$
{\cal H}' \to {\cal H} \to {\cal H}''
\qquad \text{and} \qquad
\supp H \subset \supp H' \cup \supp H''.
$$
\end{lem}

For a {\em $\BB Z$-graded} $\BB C[X_1,\dots,X_l]$-module $H$ where
$\deg X_j =2$, the ideal $Ann(H)$ is graded, i.e. it is generated by
homogeneous polynomials. Hence the variety $\supp H \subset \BB C^l$ is
{\em conic}, i.e. invariant under the scaling action of $\BB C^{\times}$:
$$
x \in \supp H \quad \Longleftrightarrow \quad \lambda x \in \supp H,
$$
for all $\lambda \in \BB C^{\times}$.
Note that for $l=1$ this implies that
$\supp H^{tor} \subset \{0\} \subset \BB C$.
Since the variables $X_j$ have degree $2$, the graded module $H$
is the direct sum of two submodules:
$$
H = H^{even} \oplus H^{odd}.
$$
These can be localized separately and, in particular,
we can define their ranks.
Note that the integer grading essentially gets lost on localization.

For the situation at hand, $H^*_T(M)$ is a $\BB Z$-graded 
$\BB C[\mathfrak{t}]$-module, and we can talk about
$\supp H^*_T(M)$ which is a subset (variety) in
$\mathfrak{t}_{\BB C} \simeq \BB C^l$.

\begin{lem}
Suppose there is a closed subgroup $K \subset T$ and a
$T$-equivariant map $M \to T/K$, then
$$
\supp H^*_T(M) \subset \mathfrak{k}_{\BB C}.
$$
In particular, if $K \ne T$ then
$H^*_T(M)$ is a torsion module over $H^*_T(pt)$.
\end{lem}

\pf
The $T$-equivariant maps $M \to T/K \to pt$ give ring homomorphisms
$H^*_T(M) \leftarrow H^*_T(T/K) \leftarrow H^*_T(pt)$. But
$$
H^*_T(T/K) \simeq H^* ( (T/K \times E)/T ) \simeq H^* (E/K) \simeq H^*_K(pt)
\simeq H^*_{K^0}(pt),
$$
where $K^0$ denotes the connected component of identity in $K$.
Thus $H^*_T(M)$ is effectively a module over $H^*_{K^0}(pt)$ which becomes
a module over $H^*_T(pt)$ by restriction from $T$ to the sub-torus $K^0$.
It is obvious now that, for any $f \in \BB C[\mathfrak{t}]$
with $f \bigr|_{\mathfrak{k}_{\BB C}} =0$, $f \cdot H^*_T(M) =0$ and
$\supp H^*_T(M) \subset V_f$. It follows that
$\supp H^*_T(M) \subset \mathfrak{t}_{\BB C}$.
\qed

One situation in which the lemma applies is the following.
Let $p$ be a point in $M$ with isotropy (stabilizer) group $K$.
Let $N$ be the orbit of $T$ through $p$ and ${\cal N}$ be the conormal
(or normal) bundle of $N$ in $M$.
The equivariant tubular neighborhood theorem asserts that there exists
a $T$-invariant tubular neighborhood $U \subset M$ of $N$ and a
$T$-equivariant diffeomorphism $\gamma: {\cal N} \to U$ such that
$$
\gamma \circ j_{N \hookrightarrow {\cal N}} = j_{N \hookrightarrow M},
$$
where $j_{N \hookrightarrow {\cal N}}$ and $j_{N \hookrightarrow M}$
are the inclusions of $N$ respectively into ${\cal N}$ as the zero section
and into $M$.
Being a tubular neighborhood, $U$ has a $T$-equivariant projection onto $N$.
Therefore, identifying $N$ with $T/K$ we get a $T$-equivariant map
$\pi: U \to T/K$. Moreover, for any $T$-invariant open subset $U' \subset U$,
we can restrict $\pi$ to $U'$ and regard it as a $T$-invariant map of
$U'$ onto $T/K$. Thus we have proved

\begin{lem}
There exists a $T$-invariant neighborhood $U$ of $p$ with the property that
for every $T$-invariant neighborhood $W$ of $p$,
$$
\supp H^*_T (U \cap W) \subset \mathfrak{k}_{\BB C}.
$$
\end{lem}

Note that in this situation the isotropy group of each point in $U$
is contained in $K$. Moreover, since $K$ acts linearly on the fiber of
${\cal N}$ over $p$, there are only a finite number of distinct
closed subgroups of $T$ that can occur as isotropy groups of points in $U$.
Since the manifold $M$ is compact and can be covered by finitely many
such tubular neighborhoods, it proves

\begin{lem}
The set of all subgroups of $T$ which occur as isotropy groups of points of
$M$ is finite.
\end{lem}

Just as we defined the equivariant cohomology $H^*_G(M) = H^*((M \times E)/G)$,
given any $G$-invariant subset $N \subset M$, we can define the relative
equivariant cohomology
$$
H^*_G(M,N) = H^*((M \times E)/G, (N \times E)/G).
$$

\begin{prop}  \label{AB-prop}
Let $T$ act smoothly on the compact manifold $M$ and
$N$ a closed $T$-invariant submanifold. Then the supports of the modules
$H^*_T(M \setminus N)$ and $H^*_T(M, N)$ are contained in the set
\begin{equation}  \label{set_k}
\bigcup_{K \subset T} \mathfrak{k}_{\BB C}
\end{equation}
where $K$ runs over the finite set of all subgroups of $T$ which occur
as isotropy groups of points of $M \setminus N$.
\end{prop}

\pf
Let $U$ be a $T$-invariant tubular neighborhood of $N$. Then
$$
H^*_T(M \setminus N) \simeq H^*_T(M \setminus \overline{U})
\quad \text{and} \quad
H^*_T(M, N) \simeq H^*_T(M, \overline{U}),
$$
as can be seen by an equivariant deformation argument.
Thus it is enough to prove the assertion for
$H^*_T(M \setminus \overline{U})$ and $H^*_T(M, \overline{U})$.
Since $M \setminus U$ is compact, one can find $T$-invariant open
sets $U_i$, $i=1,\dots,k$, covering $M \setminus \overline{U}$
and equivariant maps $\pi_i : U_i \to G/K_i$, each $K_i$
being an isotropy group of a point in $M \setminus N$.
Let $V_r = U_1 \cup \dots \cup U_{r-1}$. We will show by induction that
the support of $H^*_T(V_r)$ is contained in the set (\ref{set_k}).
We have the equivariant Mayer-Vietoris sequence of
$\BB C [\mathfrak{t}]$-modules:
$$
H^*_T(U_r \cap V_r) \to H^{*+1}_T(V_{r+1}) \to
H^{*+1}_T(U_r) \oplus H^{*+1}_T(V_r).
$$
In this sequence the end terms are supported in (\ref{set_k}) by induction.
Hence by Lemma \ref{support-lemma} the middle terms are as well.
For $r=k+1$ we get
$$
\supp H^*_T(M \setminus \overline{U}) \subset
\bigcup_{i=1}^k (\mathfrak{k}_i)_{\BB C}
$$
which proves the first part of the proposition.

Excising $N$ we get an isomorphism
$$
H^*_T (M, \overline{U}) \simeq H^*_T (M \setminus N, \overline{U} \setminus N).
$$
Now the second part of the proposition follows from the exact sequence
$$
H^*_T (\overline{U} \setminus N) \to
H^{*+1}_T (M \setminus N, \overline{U} \setminus N)
\to H^{*+1}_T (M \setminus N),
$$
the first part of the proposition and Lemma \ref{support-lemma}.
\qed

\begin{cor}  \label{AB-cor}
In the setting as above let $j_{N \hookrightarrow M} : N \hookrightarrow M$
be the inclusion map. Then the kernel and cokernel of the map on equivariant
cohomology
$$
(j_{N \hookrightarrow M})^* : H^*_T(M) \to H^*_T(N)
$$
are supported in the set (\ref{set_k}).
\end{cor}

\pf
From the exact sequence
$$
H^*_T(M, N) \to H^*_T(M) \to H^*_T(N) \to H^{*+1}_T(M,N)
$$
we observe that the kernel and cokernel of $(j_{N \hookrightarrow M})^*$
are a quotient module and a submodule of $H^*_T(M, N)$ respectively.
Thus the supports of both of these modules are contained in (\ref{set_k}).
\qed

Letting $N$ be the set of fixed points $F$ we obtain the Atiyah-Bott
localization theorem.

\begin{thm} [Atiyah-Bott]  \label{AB-loc-thm}
The kernel and cokernel of
$$
i^*: H^*_T(M) \to H^*_T(F)
$$
have supports in
$$
\bigcup_{K \subsetneq T} \mathfrak{k}_{\BB C},
$$
where $K$ runs over the finite set of all isotropy subgroups $\ne T$.
In particular both modules have the same rank.
\end{thm}

Since $H^*_T(F) \simeq H^*(F) \otimes H^*_T(pt)$ is a free module,
the last part of the theorem implies

\begin{cor}
The rank of $H^*_T(M)$ equals $\dim H^*(F)$ and the same holds
for the individual odd and even ranks.
\end{cor}

The theorem also implies that $H^*_T(M)$ becomes a free module once we
localize to an open set $U_f \subset \mathfrak{t}_{\BB C}$ where $f$ is any
polynomial which vanishes on the set (\ref{set_k}).

The result takes an especially simple form when $T=S^1$ is a circle.
Then (unless $M=F$) $\mathfrak{k}_{\BB C} = 0$, so that
the kernel and cokernel of $i^*$ have supports at
$0 \in \BB C \simeq \mathfrak{t}_{\BB C}$ and $H^*_{S^1}(M)$ becomes free
on $\BB C \setminus \{0\}$. Equivalently the torsion subgroup is annihilated
by a power of $X_1$.
As we saw earlier, this follows directly from the fact that
$H^*_T(M)$ is a graded module. In the general torus case however the grading
simply tells us that that the support of the torsion submodule of
$H^*_T(M)$ is a proper cone, while Theorem \ref{AB-loc-thm} is more precise.

If $f: N \to M$ is a map of compact oriented manifolds, associated to this map
we have a push-forward map on ordinary cohomology
$$
f_* : H^*(N) \to H^{*+q}(M),
\qquad q = \dim M - \dim N,
$$
obtained as a composition of Poincare duality maps on $M$ and $N$ and the
push-forward map on homology:
$$
H^j(N) \simeq H_{\dim N-j}(N) \to H_{\dim N-j}(M) \simeq
H^{\dim M - \dim N+j}(M).
$$
Recall that this cohomology push-forward map has the following properties:
\begin{itemize}
\item
It is functorial: $(f \circ g)_* = f_* \circ g_*$;

\item
It is a homomorphism of $H^*(M)$-modules:
$$
f_* ( v f^*u) = (f_*v)u,
\qquad \forall u \in H^*(M),\: v \in H^*(N);
$$

\item
If $f$ is a fibering,
\begin{equation}  \label{2.14}
\text{$f_*$ corresponds to integration over the fiber;}
\end{equation}

\item
When $f: N \hookrightarrow M$ is the inclusion of a submanifold,
$f_*$ factors through the Thom isomorphism: that is, in the diagram
$$
\begin{matrix}
H^{*-1}(M \setminus N) & \xrightarrow{\delta} & H^*(M, M \setminus N) &
\xrightarrow{j^*} & H^*(M)  \\
\\
& & \qquad \bigr \uparrow \wr \: Thom_N & & \\
\\
& & H^{*-q}(N) & &
\end{matrix}
$$
we have
\begin{equation}  \label{2.17}
f_* = j^* \circ Thom_N,
\end{equation}
with $Thom_N$ the Thom isomorphism $H^*(N) \simeq H^{*+q}(M, M \setminus N)$.
\end{itemize}

Here of course $H^*(M, M \setminus N)$ is, by excision, a purely $N$-local
object, so that in the differentiable category this group can be identified
with $H^*({\cal N}, {\cal N} \setminus N)$, where ${\cal N}$ is
the conormal (or normal) bundle to $N$ in $M$, and in turn with the compactly
supported cohomology of ${\cal N}$:
$$
H^*(M, M \setminus N) \simeq H^*({\cal N}, {\cal N} \setminus N)
\simeq H^*_c({\cal N}).
$$
The image $Th_{\cal N} = Thom_N (1) \in H^*_c({\cal N})$ of $1 \in H^*(N)$
is called the {\em Thom class} of the normal bundle and that its restriction
to $N$ is precisely the {\em Euler class} $\chi({\cal N}) \in H^*(N)$:
$$
Th_{\cal N} \bigr |_N = \chi({\cal N}).
$$
In other words, for the inclusion
$j_{N \hookrightarrow M} : N \hookrightarrow M$ we have
$$
(j_{N \hookrightarrow M})^* \circ (j_{N \hookrightarrow M})_* 1 =
\chi({\cal N}).
$$
Let $\pi : {\cal N} \to N$ be the projection map.

On the level of differential forms the map
$Thom_N : H^*(N) \to H^{*+\operatorname{rk}({\cal N})}_c({\cal N})$ is given by
$$
H^*(N) \ni \quad \alpha \mapsto Th_{\cal N} \wedge \pi^* \alpha \quad
\in H^{*+\operatorname{rk}({\cal N})}_c({\cal N}).
$$
Its inverse map $H^*_c({\cal N}) \to H^{*-\operatorname{rk}({\cal N})}(N)$
is given by
$$
H^*_c({\cal N}) \ni \quad \beta \mapsto
(2\pi)^{-\operatorname{rk}({\cal N})/2} \cdot \pi_* \beta
\quad \in H^{*-\operatorname{rk}({\cal N})}(N),
$$
where $\pi_*$ is the integration over the fiber map.
Moreover, for every closed differential form $\beta \in \Omega^*({\cal N})$,
\begin{equation}  \label{thom-int}
\int_{\cal N} Th_{\cal N} \wedge \beta =
(2\pi)^{\operatorname{rk}({\cal N})/2} \cdot \int_N \beta \bigr|_N.
\end{equation}

This push-forward is now seen to extend word-for-word to the equivariant
situation. To verify this recall that by the graph construction every map
can be factored into an inclusion followed by a fibering (and in fact product)
projection. Hence it is sufficient to check the properties
(\ref{2.14}) and (\ref{2.17}) in the equivariant setting.
Now when $f$ is a fibering, so is the induced mapping
$f^G : (N \times E)/G \to (M \times E)/G$,
and integration over the fiber is well-defined in any fibering with an
oriented compact manifold as fiber.
Similarly, the usual Thom isomorphism, but now applied to bundles over
$(M \times E)/G$, extends the classical one to the equivariant theory.

It should be clear from this discussion that this equivariant push-forward
-- denoted by $f^G_*$ -- preserves the $H^*_G$-module structure, and finally
that the push-forward $(\pi^G_{M \twoheadrightarrow pt})_*$ of the map
$M \xrightarrow{\pi} pt$ corresponds to integration over fiber in the fibering
$(M \times E)/G \to E/G$.

We have the following exact sequence of $\BB C[\mathfrak{t}]$-modules:
$$
\begin{matrix}
H_T^{*-1}(M \setminus F) & \to & H_T^*(M, M \setminus F) &
\to & H_T^*(M) & \to & H_T^*(M \setminus F) \\
\\
& & \qquad \bigr \uparrow \wr \: Thom_F & \nearrow_{i_*} & & & \\
\\
& & H^{*-q}(F) & & & &
\end{matrix}
$$
Thus from Proposition \ref{AB-prop} and Lemma \ref{support-lemma} we conclude:

\begin{thm} [Atiyah-Bott]
The kernel and cokernel of the push-forward map
$$
i_*: H^*_T(F) \to H^{*+\dim M - \dim F}_T(M)
$$
have supports in
$$
\bigcup_{K \subsetneq T} \mathfrak{k}_{\BB C},
$$
where $K$ runs over the finite set of all isotropy subgroups $\ne T$.
\end{thm}

This can be seen more directly and explicitly as follows.
Let ${\cal N}$ be the conormal (or normal) bundle of $F$ in $M$ and
$\chi_T({\cal N})$ its equivariant Euler class, then the composition
$i^* \circ i_* : H^*_T(F) \to H^*_T(F)$ is multiplication by the
equivariant Euler class:
$$
i^* \circ i_* u = \chi_T({\cal N}) \cdot u,
\qquad \forall u \in H^*_T(F).
$$
It follows that $\chi_T({\cal N})$ must be invertible in the localized module
$$
H^*_T(F)_f = H^*(F) \otimes \BB C[\mathfrak{t}]_f,
$$
where $f \in \BB C[\mathfrak{t}]$ is a suitable polynomial.
Let $F = \coprod F_j$ be the decomposition of $F$ into its connected
components, and we may consider separately each $\chi_T({\cal N}_j)$.
Recall that the component of $\chi_T({\cal N}_j)$ in
$H^0(F_j) \otimes \BB C[\mathfrak{t}]$ is
$$
\chi_T({\cal N}_j)(Y)_{[0]} = {\det}^{1/2} (-L^{{\cal N}_j}_Y)
$$
which is (up to a sign) the product of weights of the linear action
of $T$ on the fiber of ${\cal N}_j$. This action has no fixed non-zero
vectors (because the only fixed directions are tangental to $F_j$).
Thus $\chi_T({\cal N}_j)_{[0]}$ is a homogeneous polynomial
on $\mathfrak{t}_{\BB C}$ of degree $\operatorname{rk}({\cal N}_j)/2$.
If we denote this polynomial by $f_j$, it follows that $i^* \circ i_*$
becomes invertible after localizing with respect to
$$
f_F = \prod_j f_j.
$$

Working over the localized ring $\BB C[\mathfrak{t}]_{f_F}$ (or over the
full field of rational functions) we see that the operator
$$
Q = \sum_j \frac {(j_{F_j \hookrightarrow M})^*}{\chi_T({\cal N}_j)}
: \quad H^*_T(M) \to H^*_T(F)
$$
is inverse to $i_*: H^*_T(F) \to H^*_T(M)$.
Thus for any $\alpha \in H^*_T(M)$ we have (after localizing)
$$
\alpha = i_* \circ Q \alpha = \sum_j \frac
{(j_{F_j \hookrightarrow M})_* \circ (j_{F_j \hookrightarrow M})^* \alpha}
{\chi_T({\cal N}_j)}
= \sum_j \frac
{(j_{F_j \hookrightarrow M})_* \bigl(\alpha \bigr|_{F_j} \bigr)}
{\chi_T({\cal N}_j)}.
$$
Applying the push-forward to a point
$$
(\pi^T_{M \twoheadrightarrow pt})_*: H^*_T(M) \to H^*_T(pt)
= \BB C[\mathfrak{t}]
$$
to both sides, using the functoriality of push-forwards and (\ref{thom-int})
we deduce the Berline-Vergne integral localization formula:
$$
(\pi^T_{M \twoheadrightarrow pt})_* \alpha =
\sum_j (2\pi )^{\operatorname{rk}({\cal N}_j)/2} \cdot
(\pi^T_{F_j \twoheadrightarrow pt})_* \biggl( \frac
{\alpha \bigr|_{F_j}} {\chi_T({\cal N}_j)} \biggr)
$$
or
$$
\int_M \alpha
= \sum_j (2\pi )^{\operatorname{rk}({\cal N}_j)/2} \cdot
\int_{F_j} \frac {\alpha \bigr|_{F_j}} {\chi_T({\cal N}_j)}.
$$
Note that the left hand side is a polynomial on $\mathfrak{t}_{\BB C}$
while the right hand side appears to be a rational function. Part of the
statement is that the denominators cancel when the sum is simplified.

\separate

\section{Equivariant Cohomology via Spectral Sequences}

The goal of this section is to prove the equivariant cohomology isomorphism
$$
H^*_G(M) \simeq \bigl( H^*_T(M) \bigr)^W,
$$
where $M$ is a compact manifold, $G$ is a connected compact Lie group,
$T \subset G$ is a maximal torus and $W$ is its Weyl group.
We will do this by realizing the equivariant cohomology as the cohomology
of a double complex and applying some basic spectral sequence techniques.
Along the way we derive some other useful properties of $H^*_G(M)$.
We closely follow Chapter 6 of \cite{GS}.

\separate

\subsection{The Cartan Model as a Double Complex}

A {\em double complex} is a bigraded vector space
$$
C = \bigoplus_{p,q \in \BB Z} C^{p,q}
$$
with coboundary operators
$$
d: C^{p,q} \to C^{p,q+1}, \qquad \delta: C^{p,q} \to C^{p+1,q}
$$
satisfying
$$
d^2=0, \qquad d\delta + \delta d =0, \qquad \delta^2=0.
$$
The {\em associated total complex} is defined by
$$
C^n =_{def} \bigoplus_{p+q=n} C^{p,q}
$$
with coboundary $d+\delta: C^n \to C^{n+1}$.
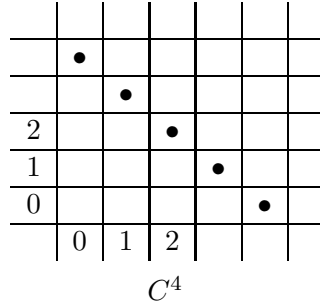
\begin{center}
\begin{tabular} {r|c|c|c|c|c|c}
 & & & & & & \\
\hline
 & $\bullet$ & & & & & \\
\hline
 & & $\bullet$ & & & & \\
\hline
2 & & & $\bullet$ & & & \\
\hline
1 & & & & $\bullet$ & & \\
\hline
0 & & & & & $\bullet$ & \\
\hline
 & 0 & 1 & 2 & & & \\
\end{tabular}

\separate

$C^4$
\end{center}
Note that $(d+\delta)^2=0$, and we can define
\begin{align*}
Z^n &= \ker(d+\delta: C^n \to C^{n+1}) = \{ z \in C^n ;\: (d+\delta)z=0 \},
\qquad n \in \BB Z,  \\
B^n &= \operatorname{im} (d+\delta: C^{n-1} \to C^n) = (d+\delta)C^{n-1},  \\
H^n(C^*,d+\delta) &= Z^n/B^n.
\end{align*}

The twisted de Rham complex
$$
(\Omega^*_G(M), d_{eq}) =
\bigl( ( \BB C [\mathfrak{g}] \otimes \Omega^*(M) )^G, d-\iota \bigr)
$$
can be thought of as a double complex with bigrading
\begin{equation}  \label{bigrading}
C^{p,q} = ( \BB C [\mathfrak{g}]^p \otimes \Omega^{q-p}(M) )^G,
\end{equation}
where $\BB C [\mathfrak{g}]^p$ denotes the space of homogeneous polynomials
on $\mathfrak{g}$ of degree $p$,
and with vertical and horizontal operators given by
$$
d = 1 \otimes d
\quad \text{(the de Rham operator on $\Omega^*(M)$)}
$$
and
$$
\delta = -\iota
\quad \text{(contraction by the vector field $-L_X$).}
$$
Notice that in the bigrading (\ref{bigrading}) the subspace
$( \BB C [\mathfrak{g}]^p \otimes \Omega^m(M) )^G$ has bidegree $(p,m+p)$
and hence total degree $2p+m$ which is the grading we have been using on
$\Omega^*_G(M)$.
\begin{center}
\begin{tabular} {r|c|c|c|c|c|c}
 & & & & & & \\
\hline
 & $\Omega^4(M)^G$ & $( \BB C [\mathfrak{g}]^1 \otimes \Omega^3(M) )^G$ &
$( \BB C [\mathfrak{g}]^2 \otimes \Omega^2(M) )^G$ &
$( \BB C [\mathfrak{g}]^3 \otimes \Omega^1(M) )^G$ &
$( \BB C [\mathfrak{g}]^4 \otimes \Omega^0(M) )^G$ & \\
\hline
 & $\Omega^3(M)^G$ & $( \BB C [\mathfrak{g}]^1 \otimes \Omega^2(M) )^G$ &
$( \BB C [\mathfrak{g}]^2 \otimes \Omega^1(M) )^G$ &
$( \BB C [\mathfrak{g}]^3 \otimes \Omega^0(M) )^G$ & 0 & \\
\hline
2 & $\Omega^2(M)^G$ & $( \BB C [\mathfrak{g}]^1 \otimes \Omega^1(M) )^G$ &
$( \BB C [\mathfrak{g}]^2 \otimes \Omega^0(M) )^G$ & 0 & 0 &  \\
\hline
1 & $\Omega^1(M)^G$ & $( \BB C [\mathfrak{g}]^1 \otimes \Omega^0(M) )^G$
& 0  & 0 & 0 & \\
\hline
0 & $\Omega^0(M)^G$ & 0 & 0 & 0 & 0 & \\
\hline
 & 0 & 1 & 2 & & & \\
\end{tabular}

\separate

$\Omega^*_G(M)$
\end{center}
This diagram clearly shows that $C^{i,j} =0$ if either $i<0$ or $j<0$.

\separate

\subsection{Spectral Sequences of Double Complexes}

Let $C = \bigoplus_{p,q \in \BB Z} C^{p,q}$ be a double complex.
In this subsection we discuss the properties of the spectral sequence
associated to this double complex. Of course, we are primarily interested
in the twisted de Rham complex realized as a double complex.

Let
$$
C^n_k =_{def} \bigoplus_{p+q=n, \: p \ge k} C^{p,q},
$$
so $C^n_k$ consists of all elements of $C^n$ which live to the right of a
vertical line.
\begin{center}
\begin{tabular} {r|c|c|c|c|c|c}
 & & & & & & \\
\hline
 & & & & & & \\
\hline
 & & & & & & \\
\hline
2 & & & $\bullet$ & & & \\
\hline
1 & & & & $\bullet$ & & \\
\hline
0 & & & & & $\bullet$ & \\
\hline
 & 0 & 1 & 2 & & & \\
\end{tabular}

\separate

$C^4_2$
\end{center}

Let
$$
Z^n_k = Z^n \cap C^n_k, \qquad B^n_k = B^n \cap C^n_k \subset Z^n_k.
$$
The map $Z^n_k \to Z^n_k / B^n_k = H^n_k$ gives a decreasing filtration
$$
\dots \supset H^n_{k-1} \supset H^n_k \supset H^n_{k+1} \supset \dots
$$
of the cohomology group $H^n(C^*,d+\delta)$.
We denote the successive quotients by $H^{k, n-k}$ and let
\begin{equation}  \label{H^{p,q}}
\operatorname{gr} H^n =_{def} \bigoplus_k H^{k,n-k},
\qquad H^{k,n-k} =_{def} H^n_k / H^n_{k+1}.
\end{equation}
The spectral sequence of a double complex is a scheme for computing these
quotients $H^{k, n-k}$ starting from the cohomology groups of the
``vertical complexes'' $(C^{r,*},d)$.
In general, a spectral sequence is a sequence of complexes $(E^*_r,\delta_r)$
such that each $E^*_{r+1}$ is the cohomology of the preceding one,
$E^*_{r+1} = H^*(E^*_r,\delta_r)$, and (under suitable hypotheses)
the ``limit'' of these complexes are the quotients $H^{p,q}$.

Any element $\tilde a \in C^n$ has all its components on the (anti-)diagonal
line $\bigoplus_{i+j=n} C^{i,j}$ and will have a ``leading term'' at position
$(p,q)$ where $p$ denotes the largest $i$ such that $\tilde a \in C^n_i$
(or, equivalently, $p$ is the smallest $i$ such that the
$(i,n-i)$ component of $\tilde a$ does not vanish).
Let $Z^{p,q}$ denote the set of such leading terms of cocyles at
position $(p,q)$. In other words, $Z^{p,q}$ denotes the set of all
$a \in C^{p,q}$ with the property that the system of equations
\begin{equation}  \label{zigzag}
\begin{matrix}
da = 0  \\
\delta a = -d a_1  \\
\delta a_1 = -da_2  \\
\delta a_2 = -da_3  \\
\vdots
\end{matrix}
\end{equation}
admits a solution $(a_1, a_2, a_3, \dots)$ where $a_i \in C^{p+i, q-i}$.
In other words, $a \in Z^{p,q} \subset C^{p,q}$ can be ratcheted by a
sequence of zigzags to any position below $a$ on the (anti-)diagonal line
$l$ through $a$, where $l = \{(i,j) ;\: i+j=p+q \}$:
$$
\begin{matrix}
0  \\
\uparrow  \\
a & \to & 0  \\
& & \uparrow  \\
& & a_1 & \to & 0   \\
& & & & \uparrow  \\
& & & & a_2 & \to & 0  \\
& & & & & &\uparrow  \\
& & & & & & a_3 & \to
\end{matrix}
\begin{matrix}
\text{the vertical arrows are $d$ and}  \\
\text{the horizontal arrows are $\delta$}
\end{matrix}
$$
Note that $Z^{p,q}$ can be canonically identified with
$Z^{p+q}_p/Z^{p+q}_{p+1}$.

In the case of twisted de Rham complex (\ref{bigrading}) we have $C^{i,j} =0$
if either $i<0$ or $i-j<0$. In particular, there is an $m_l \in \BB N$
(which we may take equal $p+q$) such that
\begin{equation}  \label{boundedness}
C^{i,j}=0 \quad \text{for $i+j=p+q$, $|i-j|>m_l$.}
\end{equation}
Hence the system (\ref{zigzag}) will be solvable for all $i$ provided that
it is solvable for $i \le m_l+1$.

Let $B^{p,q} \subset C^{p,q}$ consist of the $(p,q)$ components of
$B^{p+q}_p = B^{p+q} \cap C^{p+q}_p$. In other words, $B^{p,q}$ denotes
the set of all $b \in C^{p,q}$
with the property that the system of equations
\begin{equation}  \label{zigzag2}
\begin{matrix}
dc_0 + \delta c_{-1} = b  \\
dc_{-1} + \delta c_{-2} = 0  \\
dc_{-2} + \delta c_{-3} = 0  \\
dc_{-3} + \delta c_{-4} = 0  \\
\vdots
\end{matrix}
\end{equation}
admits a solution $(c_0,c_{-1},c_{-2}, \dots)$ with $c_{-i} \in C^{p-i,q+i-1}$:
$$
\begin{matrix}
\uparrow  \\
c_{-3} & \to & 0  \\
& & \uparrow  \\
& & c_{-2} & \to & 0  \\
& & & & \uparrow  \\
& & & & c_{-1} & \to & b \\
& & & & & &\uparrow  \\
& & & & & & c_0 & \to & ?
\end{matrix}
\begin{matrix}
\text{the vertical arrows are $d$ and}  \\
\text{the horizontal arrows are $\delta$}
\end{matrix}
$$
Once again, because of the boundedness condition (\ref{boundedness})
it suffices to solve these equations for $|i| \le m_l+1$.

Recall that the quotients $H^{p,q}$ defined in (\ref{H^{p,q}}) are
$$
H^{p,q} = H^{p+q}_p / H^{p+q}_{p+1}, \quad \text{where} \quad
H^{p+q}_p = Z^{p+q}_p / B^{p+q}_p, \qquad
H^{p+q}_{p+1} = Z^{p+q}_{p+1} / B^{p+q}_{p+1}.
$$
Hence
\begin{equation}  \label{6.5}
H^{p,q} = Z^{p+q}_p / (Z^{p+q}_{p+1} + B^{p+q}_p) = Z^{p,q} / B^{p,q}.
\end{equation}
We try to compute these quotients by solving the system (\ref{zigzag})
inductively.
Let $Z^{p,q}_r \subset C^{p,q}$ consist of those $a \in C^{p,q}$ for which
the first $r-1$ of the equations (\ref{zigzag}) can be solved.
In other words, $a \in Z^{p,q}_r$ if and only if it can be joined by
a sequence of zigzags to an element $a_{r-1} \in C^{p+r-1,q-r+1}$.
When can such an $a$ be joined by a sequence of $r$ zigzags to an element of
$C^{p+r,q-r}$? In order to do so, we may have to retrace our steps and replace
the partial solution $(a_1,\dots,a_{r-1})$ by a different partial solution
$(a_1',\dots,a_{r-1}')$ so that the differences $a_i'' = a_i' - a_i$ satisfy
$$
\begin{matrix}
da_1'' = 0  \\
\delta a_1'' + d a_2'' = 0  \\
\vdots  \\
\delta a_{r-2}'' + d a_{r-1}'' = 0
\end{matrix}
$$
and to zigzag one step further down we need an $a_r'' \in C^{p+r,q-r}$
such that
$$
\delta a_{r-1}' = \delta a_{r-1} + \delta a_{r-1}'' = - d a_r''.
$$
Let us set
$$
b = \delta a_{r-1}, \quad c_0 = - a_r'', \quad c_{-i} = - a_{r-i}'',
\qquad i = 1, \dots, r-1,
$$
and $c_{-i}=0$ for $i \ge r$.
So let us define
$$
B^{p,q}_r \subset B^{p,q}
$$
be the set of all $b \in C^{p,q}$ for which there is a solution of
(\ref{zigzag2}) with $c_{-i}=0$ for $i \ge r$. Then we have proved

\begin{thm}  \label{spectral}
Let $a \in Z^{p,q}_r$. Then
$$
a \in Z^{p,q}_{r+1} \qquad \Longleftrightarrow \qquad
\delta a_{r-1} \in B^{p+r,q-r+1}_r
$$
for any solutions $(a_1, \dots, a_{r-1})$ of the first $r-1$ equations of
(\ref{zigzag}).
\end{thm}

Notice that since $\delta a_{r-1}$ satisfies the system of equations
$$
\begin{matrix}
\delta a_{r-1} = b  \\
\delta a_{r-2} + d a_{r-1} = 0  \\
\vdots  \\
\delta a + d a_1 = 0  \\
da =0
\end{matrix}
$$
we see that $b=\delta a_{r-1}$ itself is in $B^{p+r,q-r+1}_{r+1}$.
Indeed, every element of $B^{p+r,q-r+1}_{r+1}$ can be written as a sum
of the form
\begin{equation}  \label{6.9}
\delta a_{r-1} + dc
\end{equation}
with $c \in C^{p+r,q-r}$ and
$(a,a_1,\dots,a_{r-1})$ a solution of the first $r-1$ equations of
(\ref{zigzag}).
From this we can draw a number of conclusions:

\begin{enumerate}
\item
Let
$$
E^{p,q}_r = H^{p,q}_r = Z^{p,q}_r / B^{p,q}_r.
$$
Since
$$
\delta a_{r-1} \in B^{p+r,q-r+1}_{r+1} \subset Z^{p+r, q-r+1} \subset 
Z^{p+r, q-r+1}_r,
$$
we see that $\delta a_{r-1}$ projects onto an element
$\delta_r a \in E^{p+r, q-r+1}_r$.
Theorem \ref{spectral} can be rephrased as saying that an $a \in Z^{p,q}_r$
lies in $Z^{p,q}_{r+1}$ if and only if $\delta_r a =0$.

\item
Observe that if we add to $a$ an element from $B^{p,q}$
(i.e. the $(p,q)$ component of an element in $B^{p+q} \cap C^{p,q}_p$),
then the solution $(a,a_1,\dots,a_{r-1})$ of the first $r-1$ equations of
(\ref{zigzag}) will be changed by an element from $B^{p+q} \cap C^{p,q}_p$,
and $\delta a_{r-1}$ will be changed by an element from $B^{p+r,q-r+1}_r$.
Hence $\delta_r a$ only depends on the class of $a$ modulo $B^{p,q}$.
Since $B^{p,q} \supset B^{p,q}_r$ we can consider $\delta_r$ to be a map
$$
\delta_r : E^{p,q}_r \to E^{p+r, q-r+1}_r.
$$

\item
By (\ref{6.9}), the image of this map is the projection of
$B^{p+r, q-r+1}_{r+1}$ into $E^{p+r, q-r+1}_r$ and, by Theorem \ref{spectral},
the kernel of this map is the projection of $Z^{p,q}_{r+1}$ into $E^{p,q}_r$.
Thus the sequence
$$
\dots \xrightarrow{\delta_r} E^{p,q}_r \xrightarrow{\delta_r} \dots
$$
has the property that $\ker \delta_r \supset \operatorname{im} \delta_r$ and
$$
(\ker \delta_r) / (\operatorname{im} \delta_r) = E^{p,q}_{r+1}
$$
in position $(p,q)$.
\end{enumerate}

In other words, the sequence of complexes $(E^*_r, \delta_r)$, $r=1,2,3,\dots$,
has the property that $H^*(E^*_r, \delta_r)=E^*_{r+1}$.
By construction, these complexes are bigraded and $\delta_r$ is of bidegree
$(r, -r+1)$.
Moreover, if condition (\ref{boundedness}) is satisfied for all diagonals,
the ``spectral sequence'' eventually stabilizes with
$$
E^{p,q}_r = E^{p,q}_{r+1} = E^{p,q}_{r+2} = \dots
$$
for $r$ large enough (depending on $p$ and $q$). Moreover, the ``limiting''
complex, according to (\ref{6.5}), is given by
$$
E^{p,q}_{\infty} = \lim E^{p,q}_r = H^{p,q}.
$$
But keep in mind, however, that the cohomology of the double complex
$H^*(C^*, d+\delta)$ often has a multiplicative or a module structure.
Knowledge of the vector spaces $H^{p,q}$ alone does not tell much about
the multiplicative or the module structure.

The case $r=1$ is easy to describe. By definition,
$$
E^{p,q}_1 = H^q (C^{p,*}, d)
$$
is the vertical cohomology of each column. Moreover, since
$d \delta + \delta d =0$, one gets from $\delta$ the induced map on cohomology
$$
H^q (C^{p,*}, d) \to H^q(C^{p+1,*}, d)
$$
and this is the induced map $\delta_1$. So we have described $(E^*_1,\delta_1)$
and hence $E^*_2$.
The $\delta_r$ for $r \ge 2$ are more complicated, but a lot of useful
information can be obtained just from the knowledge of $(E^*_1,\delta_1)$.

\separate

\subsection{Functorial Behavior}

Let $(C, d, \delta)$ and $(C', d', \delta')$ be double complexes, and
$\rho: C \to C'$ a morphism of double complexes of bidegree $(m.n)$ which
intertwines $d$ with $d'$, and $\delta$ with $\delta'$.
This gives rise to a cochain map
$$
\rho: (C^*, d+\delta) \to (C'^*, d'+\delta')
$$
of degree $m+n$. It induces a map $\rho_{\#}$ on the total cohomology:
$$
\rho_{\#}: H^* (C^*, d+\delta) \to H^* (C'^*, d'+\delta')
$$
of degree $m+n$ and consistent with the filtrations on both sides.
Similarly $\rho$ maps the cochain complex $(C^{p,*},d)$ into
the cochain complex $(C^{p+m,*},d)$ and hence induces a map on cohomology
$$
\rho_1 : E^*_1 \to E'^*_1
$$
of bidegree $(m,n)$ which intertwines $\delta_1$ with $\delta_1'$.
Inductively we get maps
$$
\rho_r: (E^*_r, \delta_r) \to (E'^*_r, \delta_r').
$$
Here $\rho_{r+1}$ is the map on cohomology induced from $\rho_r$, where,
we recall, $E^*_{r+1} = H^* (E^*_r, \delta_r)$. It is also clear that

\begin{thm}
If the two spectral sequences converge, then
$$
\lim \rho_r = \operatorname{gr} \rho_{\#}.
$$
\end{thm}

In particular,

\begin{thm}  \label{spectral-iso}
If $\rho_r$ is an isomorphism for some $r=r_0$, then it is an isomorphism
for all $r>r_0$ and so, if both spectral sequences converge, then $\rho_{\#}$
is an isomorphism.
\end{thm}

\separate

\subsection{Gaps}

Sometimes a pattern of zeroes among the $E^{p,q}_r$ allows for easy
conclusions. Here is a typical example:

\begin{thm}
Suppose that $E^{p,q}_r =0$ when $p+q$ is odd. Then the spectral sequence
``collapses at the $E^*_r$ stage'',
i.e. $E^*_r = E^*_{r+1} = E^*_{r+2} = \dots $.
\end{thm}

\pf
The map $\delta_r : E^{p,q}_r \to E^{p+r,q-r+1}_r$ changes the parity of $p+q$.
Thus either its domain or its range is $0$. So $\delta_r=0$.
\qed

\separate

\subsection{Switching Rows and Columns}

The point of this technique is that switching $p$ and $q$,
and hence $d$ and $\delta$ does not change the total complex,
but the spectral sequence of the switched double complex can be
quite different from that of the original. We will use this technique
to study the spectral sequence that computes equivariant cohomology.
Another illustration is a simple proof that the de Rham cohomology
equals the {\v C}ech cohomology.

\begin{thm}
Let $(C, d, \delta)$ be a double complex all of whose columns are exact except
the $p=0$ column, and all of whose rows are exact except the $q=0$ row. Then
$$
H^*(C^{0,*}, d) = H^*(C^{*,0}, \delta).
$$
\end{thm}

\pf
The $E_1$ term of the spectral sequence associated with $(C, d, \delta)$ has
only one non-zero column, the column $p=0$, and in that column the entries
are the cohomology groups of $(C^{0,*},d)$. Hence $\delta_r=0$ for
$r \ge 2$ and
$$
H^* (C^*, d+\delta) = H^*(C^{0,*}, d).
$$
Switching rows and columns yields
$$
H^* (C^*, d+\delta) = H^*(C^{*,0}, \delta).
$$
Putting these two facts together produces the isomorphism stated in the
theorem.
\qed

\separate

\subsection{The Cartan Model as a Double Complex Revisited}

Recall that we realized the twisted de Rham complex as a double complex
with $C^{p,q} = ( \BB C [\mathfrak{g}]^p \otimes \Omega^{q-p}(M) )^G$,
$d = 1 \otimes d$ (the de Rham operator on $\Omega^*(M)$) and
$\delta = -\iota$ (contraction by the vector field $-L_X$).
In this subsection we apply the general results about spectral sequences
of double complexes to the twisted de Rham complex.

We begin by evaluating the $E_1$ term:

\begin{thm}
The $E_1$ term in the spectral sequence of (\ref{bigrading}) is
$$
( \BB C [\mathfrak{g}] \otimes H^*(M) )^G.
$$
More explicitly
$$
E^{p,q}_1 = ( \BB C [\mathfrak{g}]^p \otimes H^{q-p}(M) )^G.
$$
\end{thm}

\pf
The complex 
$C = \bigoplus_{p,q} ( \BB C [\mathfrak{g}]^p \otimes \Omega^{q-p}(M) )^G$
with boundary operator $d = 1 \otimes d$ sits inside the complex
$$
(C', 1 \otimes d) =
\Bigl ( \bigoplus_{p,q} \BB C [\mathfrak{g}]^p \otimes \Omega^{q-p}(M),
1 \otimes d \Bigr)
$$
and (by averaging over the group) the cohomology groups of $(C, 1 \otimes d)$
are just the $G$-invariant components of the cohomology of $(C', 1 \otimes d)$
which are the appropriately graded components of
$\BB C [\mathfrak{g}] \otimes H^*(M)$.
\qed

\begin{lem}
The connected component of the identity element in $G$ acts trivially
on $H^*(M)$.
\end{lem}

So we get

\begin{thm}  \label{E_1-thm}
If $G$ is connected then
$$
E^{p,q}_1 = (\BB C [\mathfrak{g}]^p)^G \otimes H^{q-p}(M).
$$
\end{thm}

We now may apply the gap method to conclude

\begin{thm}  \label{odd-formal}
If $G$ is connected and $H^p(M)=0$ for $p$ odd, then the spectral sequence of
the twisted de Rham complex collapses at the $E^*_1$ stage.
\end{thm}

\pf
By Theorem \ref{E_1-thm} $E^{p,q}_1=0$ when $p+q$ is odd.
\qed

\separate

\subsection{$H^*_G(M)$ as a $\BB C [\mathfrak{g}]^G$-Module}

We assume from now on that $G$ is compact and connected so that
Theorem \ref{E_1-thm} holds.
If $f \in (\BB C [\mathfrak{g}]^m)^G$, the multiplicative operator
$$
u \otimes a \mapsto fu \otimes a
$$
is a morphism of the double complex $(C, d, \delta)$ given by (\ref{bigrading})
of bidegree $(m,m)$ and so it induces a map of $H^*_G(M)$ into itself.
This map on $H^*_G(M)$ is exactly the map $\alpha \mapsto f\alpha$ given by
the structure of $H^*_G(M)$ as a module over
$\BB C [\mathfrak{g}]^G = H^*_G(pt)$.
Also, all the $E^*_r$'s in the spectral sequence also become 
$\BB C [\mathfrak{g}]^G$-modules.
From Theorem \ref{E_1-thm} we see that this module structure of $E^*_1$
is just multiplication on the left factor of 
$E^{p,q}_1 = (\BB C [\mathfrak{g}]^p)^G \otimes H^{q-p}(M)$, which shows that
$E^*_1$ is a free $\BB C [\mathfrak{g}]^G$-module.
The ring $\BB C [\mathfrak{g}]^G$ is known to be Noetherian, so if $H^*(M)$
is finite dimensional, all of its subquotients, in particular all the
$E^*_r$'s  are finitely generated as $\BB C [\mathfrak{g}]^G$-modules.
Since the spectral sequence converges to a graded version of $H^*_G(M)$,
we conclude

\begin{thm}
If $\dim H^*(M)$ is finite, then $H^*_G(M)$ is finitely generated as a
$\BB C [\mathfrak{g}]^G$-module.
\end{thm}

Another useful fact that we can extract from this argument is:

\begin{thm}  \label{free-mod}
If the spectral sequence of the Cartan double complex collapses at the
$E^*_1$ stage, then $H^*_G(M)$ is a free $\BB C [\mathfrak{g}]^G$-module.
\end{thm}

\pf
Theorem \ref{E_1-thm} shows that $E_1$ is free as a
$\BB C [\mathfrak{g}]^G$-module, and if the spectral sequence collapses at the
$E^*_1$ stage, then $E^*_1 \simeq \operatorname{gr} H^*_G(M)$ and the
isomorphism is an isomorphism of $\BB C [\mathfrak{g}]^G$-modules by
Theorem \ref{spectral-iso}. So $\operatorname{gr} H^*_G(M)$ is a free
$\BB C [\mathfrak{g}]^G$-module.
Then $H^*_G(M)$ is freely generated as a $\BB C [\mathfrak{g}]^G$-module
by the preimages of the free generators in $\operatorname{gr} H^*_G(M)$.
\qed

\separate

\subsection{Morphisms on Equivariant Cohomology}

Let $N \to M$ be a smooth map of manifolds intertwining $G$-actions.
Then we get pull back maps on (ordinary) differential forms
$$
\rho: \Omega^*(M) \to \Omega^*(N),
$$
on (de Rham) cohomology
$$
\rho_*: H^*(M) \to H^*(N)
$$
and on equivariant cohomology
$$
\rho_{\#}: H^*_G(M) \to H^*_G(N).
$$
From Theorems \ref{spectral-iso} and \ref{E_1-thm} we conclude:

\begin{thm}
If the induced map $\rho_*$ on ordinary cohomology is bijective, then so is
the induced map $\rho_{\#}$ on equivariant cohomology.
\end{thm}

\separate

\subsection{Restricting the Group}

Suppose that $G$ is a compact connected Lie group and that $K$ is
a closed subgroup of $G$ (not necessarily connected).
We then get an injection of Lie algebras
$\mathfrak{k} \hookrightarrow \mathfrak{g}$ and a projection
$\mathfrak{g}^* \twoheadrightarrow \mathfrak{k}^*$ which extends to a map
$$
\BB C [\mathfrak{g}] \to \BB C [\mathfrak{k}]
$$
and then to a map
$$
( \BB C [\mathfrak{g}] \otimes \Omega^*(M) )^G \to
( \BB C [\mathfrak{k}] \otimes \Omega^*(M) )^K
$$
which is a morphism of complexes, in fact of double complexes
$$
\Omega^*_G(M) \to \Omega^*_K(M).
$$
We thus get a restriction mapping
$$
H^*_G(M) \to H^*_K(M)
$$
and also a restriction morphism at each stage of the corresponding spectral
sequences.
The group $G$ being connected acts trivially on $H^*(M)$,
and $K$ is a subgroup of $G$, and so acts trivially on $H^*(M)$ as well,
even though it need not be connected. In particular the conclusion of
Theorem \ref{E_1-thm} applies to $K$ as well, and hence the restriction
morphism on the $E^*_1$ level is just the restriction applied to the
left hand factors in
$$
\BB C [\mathfrak{g}]^G \otimes H^*(M) \to
\BB C [\mathfrak{g}]^K \otimes H^*(M).
$$
Therefore by Theorem \ref{spectral-iso} we conclude:

\begin{thm}  \label{group-rest}
Suppose that the restriction map 
$\BB C [\mathfrak{g}]^G \to \BB C [\mathfrak{k}]^K$
is bijective. Then the restriction map $H^*_G(M) \to H^*_K(M)$
in equivariant cohomology is bijective.
\end{thm}

Unfortunately, there is only one non-trivial example we know of for which the
hypothesis of the theorem is fulfilled, but this is a very important example.
Let $T$ be a Cartan subgroup (a maximal torus) of $G$ and let
$K = N_G(T)$ be its normalizer. The quotient group $W = K/T$ is the
Weyl group of $G$. It is a finite group, so the Lie algebra of $K$ is the
same as the Lie algebra of $T$. Since $T$ is abelian, its action on
$\mathfrak{t}^*$, and hence on $\BB C [\mathfrak{t}]$, is trivial. So
$$
\BB C [\mathfrak{k}]^K = \BB C [\mathfrak{t}]^K = \BB C [\mathfrak{t}]^W.
$$
Recall that by Chevalley's Theorem the restriction
$$
\BB C [\mathfrak{g}]^G \to \BB C [\mathfrak{t}]^W
$$
is bijective, so Theorem \ref{group-rest} applies:
$$
H^*_G(M) \simeq H^*_K(M).
$$
We can do a bit more: From the inclusion $T \hookrightarrow K$ we get a
morphism of double complexes
$$
\Omega^*_K(M) \to \Omega^*_T(M)^W
$$
which induces a morphism
$$
H^*_K(M) \to H^*_T(M)^W
$$
and also a morphism at each stage of the spectral sequences.
At the $E^*_1$ level this is just the identity morphism
$$
\BB C [\mathfrak{t}]^W \otimes H^*(M) \to 
\BB C [\mathfrak{t}]^W \otimes H^*(M)
$$
and hence another application of Theorem \ref{spectral-iso} yields
$$
H^*_K(M) = H^*_T(M)^W.
$$
Putting this together with the isomorphism coming from
Theorem \ref{group-rest} we obtain the following important result:

\begin{thm}  \label{restriction-thm}
Let $G$ be a connected compact Lie group, $T$ a maximal torus and $W$ its
Weyl group. Then for any smooth compact manifold $M$ with smooth action
$G \lefttorightarrow M$ we have
$$
H^*_G(M) \simeq H^*_T(M)^W.
$$
\end{thm}

This result can actually be strengthened a bit: The tensor product
$$
\Omega^*_K(M) \otimes_{\BB C [\mathfrak{t}]^W} \BB C [\mathfrak{t}]
$$
is also a double complex since the coboundary operators on $\Omega^*_K(M)$
are $\BB C [\mathfrak{t}]^W$-module morphisms.
Moreover there is a canonical morphism
\begin{equation}  \label{6.27}
\Omega^*_K(M) \otimes_{\BB C [\mathfrak{t}]^W} \BB C [\mathfrak{t}]
\to \Omega^*_T(M),
\qquad a \otimes f \mapsto fa.
\end{equation}
The spectral sequence associated to the double complex
$\Omega^*_K(M) \otimes_{\BB C [\mathfrak{t}]^W} \BB C [\mathfrak{t}]$
converges to
$$
H^*_K(M) \otimes_{\BB C [\mathfrak{t}]^W} \BB C [\mathfrak{t}].
$$
From (\ref{6.27}) one gets a morphism of spectral sequences which is an
isomorphism at the $E^*_1$ level. Hence, at the $E^*_{\infty}$ level we have
$$
H^*_K(M) \otimes_{\BB C [\mathfrak{t}]^W} \BB C [\mathfrak{t}]
\simeq H^*_T(M).
$$

As a corollary of Theorem \ref{restriction-thm} we get the so-called
{\em splitting principle} in topology.
Suppose that the group $G$ acts on the manifold $M$ freely.
Let $\pi: M/T \to M/G$ be the projection map.
This is a smooth fibration with fiber $G/T$. One gets from it a map
$$
\pi^*: H^*(M/G) \to H^*(M/T).
$$
Moreover there is a natural action of the Weyl group $W$ on $M/T$ which
leaves fixed the fibers of $\pi$.
Hence $\pi^*$ maps $H^*(M/G)$ into $H^*(M/T)^W$.

\begin{thm}
The map
$$
\pi^*: H^*(M/G) \to H^*(M/T)^W
$$
is a bijection.
\end{thm}

\pf
This follows from Theorem \ref{restriction-thm} and the identification
$$
H^*_G(M) \simeq H^*(M/G), \qquad H^*_T(M) \simeq H^*(M/T).
$$
\qed

\separate

\begin{section}
{Equivariant Formality}
\end{section}

\separate

\subsection{Basic Properties}

Suppose that the group $G$ is connected.
By Theorem \ref{E_1-thm} the $E^*_1$ term of the spectral sequence computing
the $G$-equivariant cohomology of $M$ is
$$
E^{p,q}_1 = (\BB C [\mathfrak{g}]^p)^G \otimes H^{q-p}(M).
$$
We say that $M$ is {\em equivariantly formal} if the spectral sequence
collapses at the $E^*_1$ stage.
Then, as a $\BB C [\mathfrak{g}]^G$-module,
$$
H^*_G(M) \simeq H^*(M) \otimes \BB C [\mathfrak{g}]^G
$$
and, by Theorem \ref{free-mod}, is a free $\BB C [\mathfrak{g}]^G$-module.
Tensoring this identity with the trivial $\BB C [\mathfrak{g}]^G$-module
$\BB C$ gives
$$
H^*(M) = H^*_G(M) \otimes_{\BB C [\mathfrak{g}]^G} \BB C
$$
expressing the ordinary de Rham cohomology of $M$ in terms of its equivariant
cohomology.

For the rest of this section we assume that the group acting on the manifold
$M$ is a compact connected abelian group -- i.e. a torus -- 
which we denote by $T$.
As a corollary of Atiyah-Bott localization theorem (Theorem \ref{AB-loc-thm})
we get:

\begin{thm}
Let $F \subset M$ denote the set of points in $M$ fixed by $T$ and
$i : F \hookrightarrow M$ the inclusion map.
If $M$ is equivariantly formal, the restriction map
$i^*: H^*_T(M) \to H^*_T(F)$ is injective.
\end{thm}

The assumption ``$H^*_T(M)$ is a free $\BB C [\mathfrak{t}]$-module''
has the following important consequence.

\begin{thm}  \label{fixed-point}
For every closed subgroup $H$ of $T$ each of the connected components of
$$
M^H = \{ p \in M ;\: H \cdot p = p \}
$$
(the set of points fixed by $H$) contains a $T$-fixed point.
\end{thm}

\pf
Let $N$ be a connected component of $M^H$. We will assume that
$N \cap F = \varnothing$ and derive a contradiction.
Observe that $T$ preserves $N$ and hence acts on the normal bundle
${\cal N}$. Let $Th_T({\cal N})$ denote the $T$-equivariant
Thom class of ${\cal N}$. Just like the non-equivariant Thom class,
$Th_T({\cal N})$ can be realized by an
equivariant form supported in an arbitrary small tubular neighborhood of $N$.
The fact that $N \cap F = \varnothing$ implies that
the restriction map $H^*_T(M) \to H^*_T(F)$ sends $Th_T({\cal N})$ to zero.
But we have just seen that the restriction map is injective.
Hence $Th_T({\cal N})=0$.
Restricting the group we get a map $H^*_T(M) \to H^*_H(M)$ which sends
$Th_T({\cal N})$ into $Th_H({\cal N})$ -- the $H$-equivariant Thom class
of ${\cal N}$. Thus $Th_H({\cal N})$ is also zero.
But then the $H$-equivariant Euler class of ${\cal N}$
is $\chi_H({\cal N}) = j^* Th_H({\cal N}) =0$, where
$j: N \hookrightarrow {\cal N}$ is the inclusion map of $N$ as the zero
section. Since the $H$-equivariant Euler class of a connected component of
$M^H$ is never zero, we get a contradiction.
\qed

\separate

\subsection{The Chang-Skjelbred Theorem}

We continue our assumptions that the acting group $T$ is compact connected and
abelian and that the manifold $M$ is equivariantly formal.
Recall that $F$ denotes the set of points in $M$ fixed by $T$.
We have just seen that the map
$$
i^* : H^*_T(M) \to H^*_T(F)
$$
embeds $H^*_T(M)$ as a submodule of $H^*_T(F)$. Our goal in this subsection
is to prove the following theorem describing the image of this map.

\begin{thm} [Chang-Skjelbred]  \label{Chang-Skjelbred}
The image of $i^*$ is the set
\begin{equation}  \label{H-intersection}
\bigcap_H (j_{F \hookrightarrow M^H})^* H^*_T(M^H)
\end{equation}
the intersection being taken over all codimension-one subtori $H$ of $T$,
and $j_{F \hookrightarrow M^H}$ being the inclusion map of $F$ into $M^H$.
Moreover, it is sufficient to take the intersection only over the finite set
of those codimension-one subtori $H$ which occur as the identity components
of isotropy subgroups.
\end{thm}

First we state the following consequence of Corollary \ref{AB-cor}:

\begin{prop}
Let $H$ be a subtorus of $T$, and let $M^H$ be the set of
points in $M$ fixed by $H$. Then the restriction map
$$
(j_{M^H \hookrightarrow M})^* : H^*_T(M) \to H^*_T(M^H)
$$
is injective and its cokernel has support in the set
$$
\bigcup Lie(K) \otimes \BB C,
$$
where $K$ runs over the finite set of all isotropy subgroups
not containing $H$.
\end{prop}

\begin{prop}  \label{CS-prop}
Let $H$ be a subtorus of $T$ which occurs as the identity component
of an isotropy group. Then there exists a polynomial
$p=\alpha_1 \cdots \alpha_n$, each $\alpha_i$
being a weight of $T$ which does not vanish on $\mathfrak{h}$,
such that $p$ is in the annihilator of
$\operatorname{coker} (j_{M^H \hookrightarrow M})^*$.
In particular, for every $e \in H^*_T(M^H)$, $pe$ is in the image of
$(j_{M^H \hookrightarrow M})^* : H^*_T(M) \to H^*_T(M^H)$.
\end{prop}

\pf
Since $H$ is connected, the above proposition says that the cokernel of the
restriction map
$$
(j_{M^H \hookrightarrow M})^* : H^*_T(M) \to H^*_T(M^H)
$$
is supported on the set
$$
\bigcup \mathfrak{k}_i \otimes \BB C,
\qquad \mathfrak{k}_i \nsupseteq \mathfrak{h},
$$
where $\mathfrak{k}_i$ runs over the finite set of Lie algebras
of all isotropy subgroups.
So, for every polynomial which vanishes on this set, some power of it
annihilates $\operatorname{coker} (j_{M^H \hookrightarrow M})^*$.
In particular, since the complexified Lie algebras
$\mathfrak{k}_i \otimes \BB C$ do not contain $\mathfrak{h}$,
one can find a weight, $\alpha_i \in \mathfrak{t}_{\BB C}^*$, of $T$
which vanishes on $\mathfrak{k}_i \otimes \BB C$, but not on $\mathfrak{h}$.
By taking the product of these weights and then some power of this product
we obtain the result.
\qed

We can now give the proof of the Chang-Skjelbred Theorem.
The image of $i^*$ is obviously contained in (\ref{H-intersection}),
so it is sufficient to prove containment in the other direction, i.e.
to prove that (\ref{H-intersection}) is contained in the image of $i^*$.
Since $i^*$ embeds $H^*_T(M)$ into $H^*_T(F)$, we can regard $H^*_T(M)$
as being a submodule of $H^*_T(F)$.
Let $e_1, \dots, e_k$ be a basis of $H^*_T(M)$ as a free module over
$\BB C [\mathfrak{t}]$.
First we apply Proposition \ref{CS-prop} with $H=T$. Thus there exists a
monomial $p=\alpha_1 \cdots \alpha_n$, the $\alpha_i$'s being non-zero
weights of $T$ such that for every element $e$ in $H^*_T(F)$,
$pe$ is in $H^*_T(M)$. Hence $pe$ can be written uniquely in the form
$$
f_1e_1 + \dots + f_ke_k, \qquad f_i \in \BB C [\mathfrak{t}].
$$
Dividing by $p$
$$
e = \frac{f_1}p e_1 + \dots + \frac{f_k}p e_k.
$$
If $f_i$ and $p$ have a common factor we can eliminate it to write
\begin{equation}  \label{11.13}
e = \frac{g_1}{p_1} e_1 + \dots + \frac{g_k}{p_k} e_k.
\end{equation}
with $g_i$ and $p_i$ relatively prime.
Note that each $p_i$ is a product of a subset of the weights,
$\alpha_1, \dots, \alpha_n$, and this expression is unique.

Let us now suppose that $e$ is in the image of $H^*_T(M^H)$.
Then by Proposition \ref{CS-prop} there exist weights, $\beta_1,\dots,\beta_r$
of $T$ such that $\beta_i |_{\mathfrak{h}}$ is not identically zero and
$\beta_1 \cdots \beta_r e$ lies in $H^*_T(M)$ and so expressible as a linear
combination of the $e_i$'s.
We can assume without loss of generality that the $\beta_i$'s are contained
among $\alpha_i$'s, and hence, in the ``reduced representation'' (\ref{11.13})
of $e$, the $p_i$'s occurring in the denominators are products of $\beta_i$'s.
Thus we have proved

\begin{lem}  \label{11.5.2}
If the element (\ref{11.13}) of $H^*_T(F)$ is contained in the image of
$H^*_T(M^H)$, none of the $p_i$'s vanish on $\mathfrak{h}$.
\end{lem}

Suppose now that this element is contained in the intersection
(\ref{H-intersection}). Let $p_i = \alpha_{i_1} \cdots \alpha_{i_s}$
and assume $s>0$, i.e. $p_i \ne 1$.
Let $\tilde{\mathfrak{h}}$ be the set
$\{X \in \mathfrak{t} ;\: \alpha_{i_1}(X)=0 \}$.
Since $\alpha_{i_1}$ is a weight, $\tilde{\mathfrak{h}}$ is the Lie algebra of
a subtorus, $\tilde H$, of codimension-one, contradicting Lemma \ref{11.5.2}.
Thus, for all $i$, $p_i=1$, and $e= f_1e_1+ \dots + f_ke_k$.
Hence $e$ is in $H^*_T(M)$.

\separate

\subsection{Criteria for Equivariant Formality}

In this subsection we list without proofs some additional criteria for $M$
to be equivariantly formal.

\begin{thm} [Goresky-Kottwitz-MacPherson]
Suppose the ordinary homology of $M$, $H_k(M, \BB R)$, is generated by classes
which are representable by cycles, $\xi$, each of which is invariant under
the action of $G$. Then $M$ is equivariantly formal.
\end{thm}

Here are some other criteria:
\begin{itemize}
\item
For every compact $G$-manifold
$$
\dim \bigoplus H^i(M^G) \le \dim \bigoplus H^i(M),
$$
and $M$ is equivariantly formal if and only if this inequality is
an equality.

\item
$M$ is equivariantly formal if and only if the canonical restriction map
$H^*_G(M) \to H^*(M)$ is onto.

\item
$M$ is equivariantly formal if it possesses a $G$-invariant Bott-Morse
function whose critical set is $M^G$.

\item
In particular, if $M$ is a Hamiltonian system and $G$ is abelian,
it is equivariantly formal as a consequence of the fact that,
for a generic $X \in \mathfrak{g}$, the $X$-component of the moment map,
$\mu(X): M \to \BB R$, is a $G$-invariant Bott-Morse function whose
critical set is $M^G$.

\item
If $M$ is equivariantly formal as a $G$-manifold, and $K$ is a closed
subgroup of $G$, $M$ is equivariantly formal as a $K$-manifold.
(This follows from the fact that the restriction map, $H^*_G(M) \to H^*(M)$
factors through $H^*_K(M)$.)
\end{itemize}

\separate

\begin{section}
{Two Dimensional $G$-Manifolds}
\end{section}

This section is taken verbatim from \cite{GS}.
There is a standard action of $SO(3)$ on $S^2$ and $\BB R P^2$ and an action of
the two-dimensional torus $T^2$ on itself, and it is well known that these
are basically the {\em only} actions of a compact connected Lie group on a
two dimensional manifold. From this one easily deduces the following.

\begin{thm}
Let $T$ be an $n$-dimensional torus acting non-trivially on an oriented,
compact, connected two-manifold $M$.
If $M^T$ is non-empty, this action has the following properties:
\begin{enumerate}
\item
$M^T$ is finite.

\item
The set of elements of $T$ which act trivially on $M$ is a closed
codimension-one subgroup, $H$.

\item
$M$ is diffeomorphic to $S^2$. Moreover, this diffeomorphism conjugates
the action of $T/H$ on $M$ into the standard action of $S^1$ on $S^2$ given by
rotation about the $z$ axis.
\end{enumerate}
\end{thm}

For the sake of completeness we will sketch a proof of this theorem:
Equip $M$ with a $T$-invariant metric, and let $p$ be a fixed point.
The exponential map
\begin{equation}  \label{11.14}
\exp: T_pM \to M
\end{equation}
is surjective and intertwines the linear isotropy action of $T$ on $T_pM$
with the action of $T$ on $M$. Therefore, if an element of $T$ acts trivially
on $T_pM$ it also acts trivially on $M$.
However, there exists a weight, $\alpha$, of $T$ and a linear isomorphism,
$T_pM \sim\to \BB C$, conjugating the linear isotropy representation with
the linear action of $T$ on $\BB C$ given by
\begin{equation}  \label{11.15}
(\exp X) z = e^{i\alpha(X)} z, \qquad X \in \mathfrak{t}.
\end{equation}
Thus the subgroup $H$ of $T$ with the Lie algebra
$$
\{ X \in \mathfrak{t} ;\: \alpha(X)=0 \}
$$
acts trivially on $T_pM$ (and hence acts trivially on $M$).
Moreover, (\ref{11.14}) maps a neighborhood of 0 diffeomorphically onto a
neighborhood of $p$, and by (\ref{11.15}) $(T_pM)^T = \{0\}$.
Hence $p$ is an isolated fixed point.
Thus we have proved the first two of the assertions above.

To prove the third assertion let $v$ be the infinitesimal generator of the
action of the circle group $T/H$ on $M$. By Assertion 2 the zeroes of this
vector field are isolated, and from (\ref{11.15}) it is easy to compute
the index at $v$ at $p$ and show that it is one.
Thus the Euler characteristic of $M$
$$
\sum_{p \in M^T} \operatorname{ind}_p v
$$
is just the cardinality of $M^T$. Hence M has positive Euler characteristic
and has to be diffeomorphic to $S^2$.
(The torus $T^2$ has Euler characteristic zero and the projective space
$\BB R P^2$ is not orientable.)

Coming back to the invariant Riemann metric on M, by the Korn-Lichtenstein
theorem (cf. \cite{C}) this metric is conformally equivalent to the standard
"round" metric on $S^2$ and hence $T/H$ acts on $M$ as a one dimensional
subgroup of the group of conformal transformations $SL(2, \BB C)$ of $S^2$.
However, the connected compact one-dimensional subgroups of $SL(2,\BB C)$
are all conjugate to each other; so, up to conjugacy, the action of $T/H$ on
$M$ is the standard action of $S^1$ on $S^2$ given by rotation about the
$z$-axis.

We will now use the results above to compute the equivariant cohomology ring,
$H^*_T(M)$. Since $H$ is a closed subgroup of $T$ one has an inclusion map,
$\mathfrak{h} \hookrightarrow \mathfrak{t}$ and hence a restriction map
\begin{equation}  \label{11.16}
r: \BB C [\mathfrak{t}] \to \BB C [\mathfrak{h}].
\end{equation}

\begin{thm}
The equivariant cohomology ring $H^*_T(M)$ is the subring of
$\BB C [\mathfrak{t}] \oplus \BB C [\mathfrak{t}]$ consisting of all pairs
\begin{equation}  \label{11.17}
(f_1,f_2) \in \BB C [\mathfrak{t}] \oplus \BB C [\mathfrak{t}]
\end{equation}
satisfying
\begin{equation}  \label{11.18}
rf_1=rf_2.
\end{equation}
\end{thm}

\pf
Since the ordinary cohomology of $M$ is zero in odd dimensions,
$M$ is equivariantly formal by Theorem \ref{odd-formal}; so the map
$H^*_T(M) \to H^*_T(M^T)$ is injective, and embeds $H^*_T(M)$ into $H^*_T(M^T)$
as a subring. However, $M^T$ consists of two points, $p_1$ and $p_2$, so
$H^*_T(M^T) = \BB C [\mathfrak{t}] \oplus \BB C [\mathfrak{t}]$,
with one copy of $\BB C [\mathfrak{t}]$ for each $p_r$.
Moreover the embedding of $H^*_T(M)$ into $H^*_T(M^T)$ is given explicitly by
\begin{equation}  \label{11.19}
H^*_T(M) \ni c \mapsto (i^*_1c, i^*_2c)
\in \BB C [\mathfrak{t}] \oplus \BB C [\mathfrak{t}],
\end{equation}
$i_r$ being the inclusion map, $i_r : \{p_r\} \hookrightarrow  M$.
Now note that since $H$ acts trivially on $M$
$$
H^*_H(M) = H^*(M) \otimes \BB C [\mathfrak{h}].
$$
so the natural mapping of $H^*_T(M)$ into $H^*_H(M)$ maps the element
(\ref{11.19}) into an element of the form
$$
c' = 1 \otimes g_0 + \sum_{i=1}^l a_1 \otimes g_i
$$
with $a_i \in H^i(M)$, $i > 0$, and $g_0,\dots,g_l$ in $\BB C [\mathfrak{h}]$.
Thus
$$
i^*_1 c' = i^*_2 c' = g_0
$$
giving one the compatibility condition (\ref{11.18}).
Thus we have proved that $H^*_T(M)$ is contained in the subring of
$\BB C [\mathfrak{t}] \oplus \BB C [\mathfrak{t}]$ defined by
(\ref {11.17}) -- (\ref{11.18}).
To prove that it is equal to this ring we note that
\begin{align*}
\dim H^{2k}_T(M) &= \dim \bigl( (H^0(M) \otimes \BB C [\mathfrak{t}]^k)
\oplus (H^2(M) \otimes \BB C [\mathfrak{t}]^{k-1}) \bigr)  \\
&= \dim \BB C [\mathfrak{t}]^k + \dim \BB C [\mathfrak{t}]^{k-1}
\end{align*}
by Theorem \ref{E_1-thm}; so it suffices to check that this dimension is
the same as the dimension of the $2k$-th component of the ring
(\ref{11.17}) -- (\ref{11.18}), viz.
$$
2 \dim \BB C [\mathfrak{t}]^k - \dim \BB C [\mathfrak{h}]^k
$$
i.e., to check that
$$
\dim \BB C [\mathfrak{t}]^k = \dim \BB C [\mathfrak{h}]^k +
\dim \BB C [\mathfrak{t}]^{k-1}.
$$
This, however, follows from the fact that the restriction map, (\ref{11.16}),
is onto and that its kernel is $\alpha \cdot \dim \BB C [\mathfrak{t}]^{k-1}$,
$\alpha$ being an element of $\mathfrak{t}^* \setminus \{0\}$ which vanishes
on $\mathfrak{h}$.
\qed

\separate

\begin{section}
{A Theorem of Goresky-Kottwitz-MacPherson}
\end{section}

This section is taken verbatim from \cite{GS}.
We denote by $T$ an $n$-dimensional torus.
Let $M$ be a compact $T$-manifold having the following three properties
\begin{enumerate}
\item[a)]
$H^*_T(M)$ is a free $\BB C [\mathfrak{t}]$-module.

\item[b)]
$M^T$ is finite.

\item[c)]
For every $p \in M$ the weights
\begin{equation}  \label{11.20}
\{ \alpha_{i,p} \}, \qquad i=1,\dots,d,
\end{equation}
of the isotropy representation of $T$ on $T_pM$ are pairwise linearly
independent: i.e., for $i \ne j$, $\alpha_{i,p}$ is not a linear multiple of
$\alpha_{j,p}$.
\end{enumerate}

The role of properties a) and b) is clear. We will clarify the role of
property c) by proving:

\begin{thm}
Given properties a) and b) property c) is equivalent to: For every
codimension one subtorus $H$ of $T$, $\dim M^H \le 2$.
\end{thm}

\pf
Let $N$ be a connected component of $M^H$ of dimension greater than zero.
By Theorem \ref{fixed-point}, $N$ contains a $T$-fixed point, $p$.
Moreover,
$$
T_pN = (T_pM)^H.
$$
Therefore, since $H$ is of codimension one, its Lie algebra is equal to
\begin{equation}  \label{11.21}
\{ X \in \mathfrak{t} ;\: \alpha_{i,p}(X)=0 \}
\end{equation}
$\alpha_{i,p}$ being one of the weights on the list (\ref{11.20}),
and hence $T_pN$ is the one dimensional (complex) subspace of $T_pM$
associated with this weight.
\qed

\begin{rem}
It is clear from this proof that $\dim M = 2$ if and only if the Lie algebra
$\mathfrak{h}$ of $H$ is the algebra (\ref{11.21}) for some $i$ and $p$.
Hence there are only a finite number of subtori
\begin{equation}  \label{11.22}
\{ H_i \}, \qquad i = l,...,n,
\end{equation}
with the property that $\dim M^{H_i} = 2$, and if $H$ is not one of the
groups on the list, $M^H = M^T$.
\end{rem}

Moreover, if $H$ is one of these exceptional subtori, the connected
components $\Sigma_{i,j}$ of $M^{H_i}$ are two-spheres, and each of these
two-spheres intersects $M^T$ in exactly two points (a "north pole" and a
"south pole"). For $i$ fixed, the $\Sigma_{i,j}$'s cannot intersect each other;
however, for different $i$'s, they can intersect at points of $M^T$;
and their intersection properties can be described by an "intersection graph"
$\Gamma$ whose edges are the $\Sigma_{i,j}$'s and whose vertices are
the points of $M^T$. (Two vertices, $p$ and $q$, of $\Gamma$ are joined by an
edge, $\Sigma$, if $\Sigma \cap M^T = \{p,q\}$.)

Moreover, for each $\Sigma$ there is a unique $H_i$ on the list (\ref{11.22})
for which
\begin{equation}  \label{11.23}
\Sigma \subseteq M^{H_i}
\end{equation}
so the edges of $\Gamma$ are labeled (or colored) by the $H_i$ on this list.

Since $M^T$ is finite
$$
H^*_T(M^T) = H^0(M^T) \otimes \BB C [\mathfrak{t}] =
\operatorname{Maps} (M^T, \BB C [\mathfrak{t}])
$$
and hence
\begin{equation}  \label{11.24}
H^*_T(M^T) = \operatorname{Maps} (V_{\Gamma}, \BB C [\mathfrak{t}])
\end{equation}
where $V_{\Gamma}$ is the set of vertices of $\Gamma$.

\begin{thm} [Goresky-Kottwitz-MacPherson, \cite{GKM}]
An element, $p$, of the ring
$$
\operatorname{Maps} (V_{\Gamma}, \BB C [\mathfrak{t}])
$$
is in the image of the embedding
$$
i^*: H^*_T(M) \to H^*_T(M^T)
$$
if and only if for every edge $\Sigma$ of the intersection graph, $\Gamma$,
it satisfies the compatibility condition
\begin{equation}  \label{11.25}
r_h p(v_1) = r_h p(v_22)
\end{equation}
$v_1$ and $v_2$ being the vertices of $\Sigma$, $\mathfrak{h}$ being the Lie
algebra of the group (\ref{11.23}), and
\begin{equation}  \label{11.26}
r_h : \BB C [\mathfrak{t}] \to \BB C [\mathfrak{h}]
\end{equation}
being the restriction map.
\end{thm}

\pf
By the Chang-Skjelbred Theorem (Theorem \ref{Chang-Skjelbred}) the image
of $i^*$ is the intersection:
$$
\bigcap_k (j_{M^T \hookrightarrow M^{H_k}})^* H^*_T(M^{H_k})
$$
and by (\ref{11.18}) the image of $(j_{M^T \hookrightarrow M^{H_k}})^*$
is the set of elements of
$\operatorname{Maps} (V_{\Gamma}, \BB C [\mathfrak{t}])$
satisfying the compatibility condition (\ref{11.25}) at the vertices
of $\Gamma$ labelled by $H_k$.
\qed

\separate

\begin{section}
{Paradan-Witten Localization}
\end{section}

\separate

\subsection{Introduction}

Let $(G \lefttorightarrow M, \omega, \mu)$ be a Hamiltonian system;
we assume that the acting group $G$ is compact, while the manifold $M$
need not be compact.
Let $\tilde \omega(X) = \omega + \mu(X)$, $X \in \mathfrak{g}$,
be the equivariant symplectic form.
In this section we always assume that the moment map
$$
\mu : M \to \mathfrak{g}^* \quad \text{is proper.}
$$
Recall that for a compact manifold $M$ we have an integration map
$$
\int_M : H^*_G(M) \to \BB C[\mathfrak{g}]^G.
$$
P.-E.~Paradan (\cite{Par2}) defines an analogue of this map for non-compact
symplectic manifolds
\begin{equation}  \label{pushforward}
{\cal P}: H^*_G(M) \to {\cal D}^{-\infty}(\mathfrak{g}^*)^G,
\end{equation}
where ${\cal D}^{-\infty}(\mathfrak{g}^*)^G$ is the space of $G$-invariant
distributions on $\mathfrak{g}^*$.

Fix a $G$-invariant Euclidean norm $\|\cdot\|$ on $\mathfrak{g}^*$.
We will see that the integration map ${\cal P}$ localizes at the set of
critical points of the $G$-invariant function
$$
\|\mu\|^2 : M \to \BB R.
$$
In particular, we will see an interesting technique developed by 
P.-E.~Paradan called the ``partition of unity for equivariant cohomology.''

\separate

\subsection{Coadjoint Orbits}

In this subsection we describe a series of important examples of
Hamiltonian systems with proper moment maps
-- coadjoint orbits of real semisimple Lie groups.
This special case is very important and it was thoroughly studied
in \cite{DHV}, \cite{DV} and \cite{Par1}.
Let $G^{ss}$ be a real semisimple Lie group and $G \subset G^{ss}$
its maximal compact subgroup.
We denote by $\mathfrak{g}^{ss}$ and $\mathfrak{g}$
the Lie algebras of $G^{ss}$ and $G$ respectively.
Let $\lambda \in (\mathfrak{g}^{ss})^*$ be a semisimple element,
and consider its coadjoint orbit
$$
{\cal O}_{\lambda} = G^{ss} \cdot \lambda \subset (\mathfrak{g}^{ss})^*.
$$
When $\lambda \in (\mathfrak{g}^{ss})^*$ is semisimple, the orbit
${\cal O}_{\lambda}$ is a closed submanifold of $(\mathfrak{g}^{ss})^*$.

We saw in Proposition \ref{symplectic-form} that ${\cal O}_{\lambda}$
has a canonical structure of a Hamiltonian system
$(G^{ss} \lefttorightarrow {\cal O}_{\lambda}, \sigma_{\lambda}, \tilde\mu)$
with the moment map given by the inclusion map
$\tilde\mu: {\cal O}_{\lambda} \hookrightarrow (\mathfrak{g}^{ss})^*$.
The $G^{ss}$ action on ${\cal O}_{\lambda}$ restricts to a Hamiltonian
action of $G$, with symplectic moment map
$\mu: {\cal O}_{\lambda} \to \mathfrak{g}^*$
given by the inclusion
${\cal O}_{\lambda} \hookrightarrow (\mathfrak{g}^{ss})^*$
followed by the natural projection
$(\mathfrak{g}^{ss})^* \twoheadrightarrow \mathfrak{g}^*$.

\begin{prop}  \label{proper}
Let $G^{ss}$ be a real semisimple Lie group and take $G \subset G^{ss}$
a maximal compact subgroup.
Denote by $\mathfrak{g}^{ss}$ and $\mathfrak{g}$ their respective Lie algebras.
Let ${\cal O}_\lambda \subset (\mathfrak{g}^{ss})^*$ be a semisimple
coadjoint orbit of $G^{ss}$.
Then the restriction of the natural projection map
$(\mathfrak{g}^{ss})^* \twoheadrightarrow \mathfrak{g}^*$ to
${\cal O}_{\lambda}$ is proper.

Let $T \subset G$ be a maximal torus with Lie algebra $\mathfrak{t}$.
If $\mathfrak{t}$ is also a Cartan algebra in $\mathfrak{g}^{ss}$, then
$({\cal O}_{\lambda})^T$ -- the set of points in
${\cal O}_{\lambda}$ fixed by $T$ -- is compact.
\end{prop}

\pf
The proof uses standard facts about the structure of real semisimple 
Lie groups. For details see, for instance, \cite{Kn}.
Let $B(X,Y)$ be the Killing form on $\mathfrak{g}^{ss}$:
$$
B(X,Y) = \tr (ad(X) \circ ad(Y)), \qquad X, Y \in \mathfrak{g}^{ss}.
$$
The Killing form $B$ is a symmetric $Ad(G^{ss})$-invariant bilinear
form on $\mathfrak{g}^{ss}$ and its restriction to $\mathfrak{g}$
is negative definite.

Let $\mathfrak{g}^{ss} = \mathfrak{g} \oplus \mathfrak{p}$ be the Cartan
decomposition of $\mathfrak{g}^{ss}$,
and let $\theta$ be the Cartan involution:
$$
\theta (X) =
\begin{cases}
X, & \text{if $X \in \mathfrak{g}$};  \\
-X, & \text{if $X \in \mathfrak{p}$}.
\end{cases}
$$
The map $\theta$ is a $G$-equivariant Lie algebra automorphism of
$\mathfrak{g}^{ss}$.
Then $B$ is positive definite on $\mathfrak{p}$, the spaces $\mathfrak{g}$
and $\mathfrak{p}$ are orthogonal with respect to $B$, and
$$
B_{\theta}(X,Y) =_{def} B(-\theta(X), Y)
$$
is a symmetric positive-definite inner product on $\mathfrak{g}^{ss}$
which is $G$-invariant, but not $G^{ss}$-invariant.

The $G^{ss}$-invariant bilinear form $B$ allows us to identify
$\mathfrak{g}^{ss}$ and its dual $(\mathfrak{g}^{ss})^*$.
Under this identification
$(\mathfrak{g}^{ss})^* \ni \lambda \longleftrightarrow
\lambda' \in \mathfrak{g}^{ss}$,
$\lambda'$ is semisimple since $\lambda$ is, and the projection
$(\mathfrak{g}^{ss})^* \twoheadrightarrow \mathfrak{g}^*$ becomes
the projection
$\pi: \mathfrak{g}^{ss} = \mathfrak{g} \oplus \mathfrak{p} \twoheadrightarrow
\mathfrak{g}$.
Thus it is enough to show that $\pi$ is proper on the adjoint
orbit ${\cal O}_{\lambda'}$.

Let $K \subset \mathfrak{g}$ be any compact set. We need to show that
$\pi^{-1}(K) \cap {\cal O}_{\lambda'}$ is compact.
Since $\lambda'$ is semisimple, ${\cal O}_{\lambda'}$ is closed, and it
is sufficient to show that $\pi^{-1}(K) \cap {\cal O}_{\lambda'}$ is bounded,
i.e. there is a constant $l$ such that
$$
B_{\theta}(Z,Z) \le l, \qquad
\forall Z \in \pi^{-1}(K) \cap {\cal O}_{\lambda'}.
$$

Because $B$ is $G^{ss}$-invariant, the quadratic form $B(Z,Z)$ is constant
on ${\cal O}_{\lambda'}$; let us call this constant $c$.
The set $K \subset \mathfrak{g}$ is compact, so there is a constant
$b$ such that
$$
B_{\theta}(X,X) \le b, \qquad \forall X \in K.
$$
Now suppose that $Z = X+Y \in \pi^{-1}(K) \cap {\cal O}_{\lambda'}$,
$X \in K$, $Y \in \mathfrak{p}$, then
$$
c = B(Z,Z) = B(X+Y,X+Y) = B_{\theta}(-X+Y,X+Y)
= -B_{\theta}(X,X) + B_{\theta}(Y,Y),
$$
so $B_{\theta}(Y,Y) = B_{\theta}(X,X) +c \le b+c$.
Therefore,
\begin{multline*}
B_{\theta}(Z,Z) = B_{\theta}(X+Y,X+Y) =
B_{\theta}(X,X) + B_{\theta}(Y,Y) \le 2b+c,  \\
\forall Z= X+Y \in \pi^{-1}(K) \cap {\cal O}_{\lambda'}, \:
X \in K, \: Y \in \mathfrak{p}.
\end{multline*}
This proves that $\pi^{-1}(K) \cap {\cal O}_{\lambda'}$ is bounded,
hence compact.

It remains to show that $({\cal O}_{\lambda})^T$ is compact when
$\mathfrak{t}$ is a Cartan algebra in $\mathfrak{g}^{ss}$.
Since the set $({\cal O}_{\lambda})^T$ is closed,
it is enough to show that $({\cal O}_{\lambda})^T$ or, equivalently,
$({\cal O}_{\lambda'})^T$ is bounded.

Let $X \in ({\cal O}_{\lambda'})^T$, then $Ad(T)X =X$ implies
$[\mathfrak{t}, X] = X$. Since $\mathfrak{t}$ is maximal abelian in
$\mathfrak{g}^{ss}$, it follows that $X \in \mathfrak{t}$.
Therefore, $({\cal O}_{\lambda'})^T \subset \mathfrak{t} \subset \mathfrak{g}$.
Since $B(Z,Z)=c$ for all $Z \in {\cal O}_{\lambda'}$,
$$
({\cal O}_{\lambda'})^T \subset \{X \in \mathfrak{t} ;\: B(X,X)=c \}
= \{X \in \mathfrak{t} ;\: B_{\theta}(X,X)=-c \}.
$$
This proves that $({\cal O}_{\lambda'})^T$ is bounded.
\qed

One might ask if the natural projection
$(\mathfrak{g}^{ss})^* \twoheadrightarrow \mathfrak{t}^*$
is proper on ${\cal O}_{\lambda}$.
E.~Prato (Propositions 2.2 and 2.3 in \cite{Pr}) gives a sufficient condition
for this.
However, in general the projection
$(\mathfrak{g}^{ss})^* \twoheadrightarrow \mathfrak{t}^*$
need not be proper on ${\cal O}_{\lambda}$.

\begin{ex}
Consider $G^{ss} = SL(3,\BB R)$, $G=SO(3, \BB R)$,
$$
T= \Biggl\{
\begin{pmatrix} \cos \theta & -\sin \theta & 0 \\
\sin \theta & \cos \theta & 0 \\
0 & 0 & 1 \end{pmatrix} ;\: \theta \in \BB R \Biggr\},
\qquad
\mathfrak{t} =
\Biggl\{
\begin{pmatrix} 0 & -\theta & 0 \\ \theta & 0 & 0 \\
0 & 0 & 0 \end{pmatrix} ;\: \theta \in \BB R \Biggr\},
$$
and take
$$
\lambda' = \begin{pmatrix} 1 & 0 & 0 \\ 0 & 1 & 0 \\
0 & 0 & -2 \end{pmatrix}
\in \mathfrak{sl}(3,\BB R).
$$
Then the adjoint orbit $SL(3,\BB R) \cdot \lambda' = {\cal O}_{\lambda'}$
contains a subset
$$
S = \Biggl\{
Ad \begin{pmatrix} 1 & 0 & s \\ 0 & 1 & 0 \\ 0 & 0 & 1 \end{pmatrix}
\begin{pmatrix} 1 & 0 & 0 \\ 0 & 1 & 0 \\ 0 & 0 & -2 \end{pmatrix}
;\: s \in \BB R \Biggr\}
= \Biggl\{
\begin{pmatrix} 1 & 0 & -3s \\ 0 & 1 & 0 \\ 0 & 0 & -2 \end{pmatrix}
;\: s \in \BB R \Biggr\}.
$$
Because $S \perp_B \mathfrak{t}$, the $B$-orthogonal projection
$\mathfrak{sl}(3,\BB R) \twoheadrightarrow \mathfrak{t}$ sends $S$ into $0$,
hence it is {\em not} proper on ${\cal O}_{\lambda'}$.
By taking $\lambda \in \mathfrak{sl}(3,\BB R)^*$ corresponding to $\lambda'$
we obtain a coadjoint orbit
${\cal O}_{\lambda} \subset \mathfrak{sl}(3,\BB R)^*$
such that the restriction of the projection map
$\mathfrak{sl}(3,\BB R)^* \twoheadrightarrow \mathfrak{t}^*$
to ${\cal O}_{\lambda}$ is not proper.
\end{ex}

\separate

\subsection{Fourier Transform and the Push Forward Map (\ref{pushforward})}

In this section we are primarily in $G$-equivariant forms
and cohomology on $M$ with distributional coefficients.
Let $\Omega_c^{top}(\mathfrak{g})$ be the space of smooth compactly supported
complex-valued differential forms on $\mathfrak{g}$ of top degree;
it will play the role of the space of test functions.
We equip both $\Omega_c^{top}(\mathfrak{g})$ and $\Omega^*(M)$ with
${\cal C}^{\infty}$ topologies.
By an equivariant form with distributional or ${\cal C}^{-\infty}$
coefficients we mean a continuous $\BB C$-linear $G$-equivariant map
$$
\alpha: \Omega_c^{top}(\mathfrak{g}) \ni \phi \mapsto
\langle \alpha, \phi \rangle_{\mathfrak{g}} \in \Omega^*(M).
$$
We denote the space of such equivariant forms those by $\Omega_G^{-\infty}(M)$
and treat elements of $\Omega_G^{-\infty}(M)$ as $\Omega^*(M)$-valued
$G$-equivariant distributions on $\mathfrak{g}$.

To define the twisted de Rham differential on $\Omega_G^{-\infty}(M)$, we
pick a real vector space basis $\{X_1,\dots,X_{\dim \mathfrak{g}}\}$ of
$\mathfrak{g}$ and let $\{X^1,\dots,X^{\dim \mathfrak{g}}\}$ be the
associated dual basis of $\mathfrak{g}^*$.
We can think of each $X^j$, $1 \le j \le \dim \mathfrak{g}$, as a linear
function on $\mathfrak{g}$;
in particular the product $X^j \phi$ makes sense for all
$\phi \in \Omega_c^{top}(\mathfrak{g})$.
For every $\alpha \in \Omega_G^{-\infty}(M)$, we set
\begin{equation}  \label{d_g}
\langle d_{eq} \alpha, \phi \rangle_{\mathfrak{g}} =_{def}
d \langle \alpha, \phi \rangle_{\mathfrak{g}} -
\sum_{j=1}^{\dim \mathfrak{g}}
\iota_{X_j} \langle \alpha, X^j \phi \rangle_{\mathfrak{g}},
\qquad \phi \in \Omega_c^{top}(\mathfrak{g}).
\end{equation}
Note that $\Omega_G^{\infty}(M) \subset \Omega_G^{-\infty}(M)$ and
(\ref{d_g}) agrees with the old definition of the twisted de Rham
differential on $\Omega_G^{\infty}(M)$.
The differential $d_{eq}$ defined by (\ref{d_g}) is independent of a
particular choice of basis $\{X_1,\dots,X_{\dim \mathfrak{g}}\}$ of
$\mathfrak{g}$, and we still have $(d_{eq})^2=0$.
The cohomology of $(\Omega_G^{-\infty}(M), d_{eq})$ is denoted by
$H_G^{-\infty}(M)$, it is $\BB Z_2$-graded and called the
$G$-equivariant cohomology of $M$ with distributional or generalized,
or ${\cal C}^{-\infty}$ coefficients.

Let ${\cal C}_c^{\infty}(\mathfrak{g}^*)$ be the space of test functions
on $\mathfrak{g}^*$, i.e. the space of smooth compactly supported
complex-valued functions on $\mathfrak{g}^*$ endowed with ${\cal C}^{\infty}$
topology. We also consider the space ${\cal C}^{-\infty}(\mathfrak{g}^*, M)$
of continuous $\BB C$-linear maps from 
${\cal C}_c^{\infty}(\mathfrak{g}^*)$ to $\Omega^*(M)$ and
${\cal C}^{-\infty}(\mathfrak{g}^*, M)^G$ -- the subspace of such
$G$-equivariant maps.
Similarly, we treat elements $\beta \in {\cal C}^{-\infty}(\mathfrak{g}^*,M)^G$
as $\Omega^*(M)$-valued $G$-equivariant distributions on $\mathfrak{g}^*$
and denote the value of $\beta$ at
$\psi \in {\cal C}_c^{\infty}(\mathfrak{g}^*)$
by $\langle \beta, \psi \rangle_{\mathfrak{g}^*} \in \Omega^*(M)$.
To define the differential on ${\cal C}^{-\infty}(\mathfrak{g}^*, M)^G$,
we use a vector space basis $\{X_1,\dots,X_{\dim \mathfrak{g}}\}$ of
$\mathfrak{g}$ and its associated dual basis
$\{X^1,\dots,X^{\dim \mathfrak{g}}\}$ of $\mathfrak{g}^*$.
For every $\beta \in {\cal C}^{-\infty}(\mathfrak{g}^*, M)^G$, we set
\begin{equation}  \label{d_g3}
\langle \widehat{d_{eq}} \beta , \psi \rangle_{\mathfrak{g}^*} =
d \langle \beta, \psi \rangle_{\mathfrak{g}^*} +
i \sum_{j=1}^{\dim \mathfrak{g}}
\iota_{X_j} \langle \beta, \partial_{X^j} \psi \rangle_{\mathfrak{g}^*},
\qquad \psi \in {\cal C}_c^{\infty}(\mathfrak{g}^*),
\end{equation}
where $\partial_{X^j} \psi$ denotes the partial derivative of $\psi$
relative to the basis $\{X^1,\dots,X^{\dim \mathfrak{g}}\}$ of
$\mathfrak{g}^*$.
As before, the differential $\widehat{d_{eq}}$ is independent of a
particular choice of basis $\{X_1,\dots,X_{\dim \mathfrak{g}}\}$ of
$\mathfrak{g}$, and $(\widehat{d_{eq}})^2=0$ on
${\cal C}^{-\infty}(\mathfrak{g}^*, M)^G$.

The complexes $(\Omega_G^{-\infty}(M), d_{eq})$ and
$({\cal C}^{-\infty}(\mathfrak{g}^*, M)^G, \widehat{d_{eq}})$
are related by the Fourier transform.
We denote by
$$
\Omega_{temp}^{-\infty}(\mathfrak{g},M)^G \subset \Omega_G^{-\infty}(M)
\qquad \text{and} \qquad
{\cal C}_{temp}^{-\infty}(\mathfrak{g}^*, M)^G \subset
{\cal C}^{-\infty}(\mathfrak{g}^*,M)^G
$$
the subspaces of tempered $G$-equivariant distributions on
$\mathfrak{g}$ and $\mathfrak{g}^*$ respectively with values in $\Omega^*(M)$.
Fix a Lebesgue measure $d\xi$ on $\mathfrak{g}^*$.
The Fourier transform
${\cal F} : \Omega_{temp}^{-\infty}(\mathfrak{g},M)^G \to
{\cal C}_{temp}^{-\infty}(\mathfrak{g}^*,M)^G$ is defined by
$$
\langle {\cal F}(\alpha), \psi \rangle_{\mathfrak{g}^*} =
\langle \alpha, {\cal F}(\psi) \rangle_{\mathfrak{g}}
$$
and it is normalized so that
$$
\Bigl\langle {\cal F}(\alpha),
\int_{\mathfrak{g}} e^{i\langle \xi,X \rangle} \phi(X)
\Bigr\rangle_{\mathfrak{g}^*} = \langle \alpha, \phi \rangle_{\mathfrak{g}},
\qquad \alpha \in \Omega_{temp}^{-\infty}(\mathfrak{g},M)^G ,\:
\phi \in \Omega_c^{top}(\mathfrak{g}),
$$
which can be formally expressed as
$$
\int_{\mathfrak{g}^*} e^{i \langle \xi, X \rangle} {\cal F}(\alpha) (\xi)
= \alpha(X),
\qquad \alpha \in \Omega_{temp}^{-\infty}(\mathfrak{g},M)^G ,\:
X \in \mathfrak{g} ,\: \xi \in \mathfrak{g}^*.
$$
(We do not need $G$-equivariance to define the Fourier transform,
but we will not be interested in non-equivariant distributions anyway.)
Then the following diagram commutes:
$$
\begin{CD}
\Omega_{temp}^{-\infty}(\mathfrak{g},M)^G @>{d_{eq}}>>
\Omega_{temp}^{-\infty}(\mathfrak{g},M)^G \\
@V{\cal F}VV      @VV{\cal F}V \\
{\cal C}_{temp}^{-\infty}(\mathfrak{g}^*,M)^G   @>{\widehat{d_{eq}}}>>
{\cal C}_{temp}^{-\infty}(\mathfrak{g}^*,M)^G
\end{CD}
$$
i.e. the Fourier transform ${\cal F}$ is a chain map.

For $\gamma \in \Omega^*(M)$ and $N \subset M$ an oriented submanifold
we define
$$
\int_N \gamma \quad =_{def} \quad
\int_N \bigl(\gamma_{[\dim N]} \bigr) \bigr|_N.
$$

\begin{df}
We say that a distribution $\beta \in {\cal C}^{-\infty}(\mathfrak{g}^*,M)^G$
has compact support in {\em $\mathfrak{g}^*$-mean} on $M$ if, for every test
function $\psi \in {\cal C}_c^{\infty}(\mathfrak{g}^*)$, the form
$\langle \beta, \psi \rangle_{\mathfrak{g}^*} \in \Omega^*(M)$ has compact
support in $M$.

For every $\beta \in {\cal C}^{-\infty}(\mathfrak{g}^*,M)^G$ with
compact support in $\mathfrak{g}^*$-mean on $M$ we denote by
$\int_M^{distrib} \beta$ the distribution on $\mathfrak{g}^*$ defined by
$$
{\cal C}_c^{\infty}(\mathfrak{g}^*) \ni \psi \mapsto
\Bigl\langle \int_M^{distrib} \beta, \psi \Bigr\rangle_{\mathfrak{g}^*} =_{def}
\int_M \langle \beta, \psi \rangle_{\mathfrak{g}^*} \in \BB C.
$$
\end{df}

\begin{lem}[Lemma 2.11 in \cite{Par2}]  \label{2.11}
If $\beta \in {\cal C}^{-\infty}(\mathfrak{g}^*,M)^G$ has compact support
in $\mathfrak{g}^*$-mean on $M$, then so does $\widehat{d_{eq}} \beta$,
moreover
$$
\int_M^{distrib} \widehat{d_{eq}} \beta =0.
$$
\end{lem}

\pf
Pick any $\psi \in {\cal C}_c^{\infty}(\mathfrak{g}^*)$, then (\ref{d_g3})
implies that
$$
\supp \langle \widehat{d_{eq}}\beta, \psi \rangle_{\mathfrak{g}^*}
\quad \subset \quad
\supp \langle \beta, \psi \rangle_{\mathfrak{g}^*} \cup
\bigcup_{j=1}^{\dim \mathfrak{g}}
\supp \langle \beta, \partial_{X^j} \psi \rangle_{\mathfrak{g}^*}
\quad \subset M,
$$
hence $\widehat{d_{eq}} \beta$ has compact support in
$\mathfrak{g}^*$-mean too.
We have
$$
\Bigl\langle \int_M^{distrib} \widehat{d_{eq}} \beta,
\psi \Bigr\rangle_{\mathfrak{g}^*} =
\int_M \langle \widehat{d_{eq}} \beta, \psi \rangle_{\mathfrak{g}^*} =
\int_M d \langle \beta, \psi \rangle_{\mathfrak{g}^*}
+ i \sum_{j=1}^{\dim \mathfrak{g}}
\int_M \iota_{X_j} \langle \beta, \partial_{X^j} \psi \rangle_{\mathfrak{g}^*}.
$$
The first term of the right hand side is zero because it is an integral of
an exact form, and each integral
$\int_M \iota_{X_j} \langle \beta,
\partial_{X^j} \psi \rangle_{\mathfrak{g}^*}$
is also zero because the differential forms
$\iota_{X_j} \langle \beta, \partial_{X^j} \psi \rangle_{\mathfrak{g}^*}$
have no components of maximal degree.
\qed

To every polynomial $P(X)$ on $\mathfrak{g}$ we associate a
differential operator with constant coefficients $P(-i\partial_{\xi})$
acting on ${\cal C}^{\infty}(\mathfrak{g}^*)$ so that
$$
P(X) \int_{\mathfrak{g}^*} e^{-i \langle \xi,X\rangle} \psi(\xi) \,d\xi =
\int_{\mathfrak{g}^*} e^{-i \langle \xi,X\rangle}
\bigl[P(-i\partial_{\xi}) \psi \bigr](\xi) \,d\xi,
\qquad X\in \mathfrak{g} ,\: \xi \in \mathfrak{g}^*,
$$
for all test functions $\psi \in {\cal C}_c^{\infty}(\mathfrak{g}^*)$.
We can rewrite this equation as
\begin{equation}  \label{diff-oper}
P(X) {\cal F}(\psi) = {\cal F} \bigl[P(-i \partial_{\xi}) \psi \bigr],
\qquad \forall \psi \in {\cal C}_c^{\infty}(\mathfrak{g}^*).
\end{equation}
This association extends naturally to $\Omega^*(M)$-valued polynomials on
$\mathfrak{g}$ and hence to $\Omega_G^*(M)$; for an
$\alpha(X) \in \Omega_G^*(M)$ we denote by $\alpha(-i\partial_{\xi})$ the
corresponding differential operator with values in $\Omega^*(M)$.

\begin{lem} [Lemma 2.12 in \cite{Par2}]  \label{2.12}
For every equivariant form $\alpha(X) \in \Omega_G^*(M)$ with polynomial
dependence on $X \in \mathfrak{g}$, the $\Omega^*(M)$-valued distribution
$$
\alpha \wedge e^{i\tilde\omega}: \quad
\Omega_c^{top}(\mathfrak{g}) \ni \phi \mapsto
\int_{\mathfrak{g}} \alpha \wedge e^{i\tilde\omega} \wedge \phi \in \Omega^*(M)
$$
is $G$-equivariant and tempered, hence belongs to
$\Omega_{temp}^{-\infty}(\mathfrak{g},M)^G$.
Its Fourier transform ${\cal F}(\alpha \wedge e^{i\tilde\omega})$
has compact support in $\mathfrak{g}^*$-mean on $M$.
Thus we get a distribution
$\int_M^{distrib} {\cal F}(\alpha \wedge e^{i\tilde\omega})$
which is $G$-invariant.
If $\alpha$ is equivariantly exact, i.e. $\alpha = d_{eq} \alpha'$
for some $\alpha' \in \Omega_G^*(M)$, then
$\int_M^{distrib} {\cal F}(\alpha \wedge e^{i\tilde\omega})=0$.
\end{lem}

\pf
From equation (\ref{diff-oper}) we see that, for all $m \in M$,
\begin{multline*}
\langle {\cal F}(\alpha \wedge e^{i\tilde\omega}),
\psi \rangle_{\mathfrak{g}^*} \bigr|_m =
\langle \alpha \wedge e^{i\tilde\omega},
{\cal F}(\psi) \rangle_{\mathfrak{g}} \bigr|_m =
\bigl\langle e^{i\tilde\omega},
{\cal F} \bigl[\alpha(-i\partial_{\xi}) \psi \bigr]
\bigr\rangle_{\mathfrak{g}} \Bigr|_m  \\
= e^{i\omega_m} \int_{\mathfrak{g}} e^{i \langle \mu(m), X \rangle}
{\cal F} \bigl[\alpha(-i\partial_{\xi}) \psi \bigr]
= e^{i\omega_m} \bigl[\alpha(-i\partial_{\xi}) \psi \bigr](\mu(m)),
\qquad \forall \psi \in {\cal C}_c^{\infty}(\mathfrak{g}^*).
\end{multline*}
Hence the differential form
$\langle {\cal F}(\alpha \wedge e^{i\tilde\omega}),
\psi \rangle_{\mathfrak{g}^*}$
is supported inside $\mu^{-1}(\supp \psi)$ which is compact because $\mu$
is proper.
This shows that ${\cal F}(\alpha \wedge e^{i\tilde\omega})$
has compact support in $\mathfrak{g}^*$-mean on $M$.
As the form $\alpha \wedge e^{i\tilde\omega}$ is $G$-equivariant, its
Fourier transform ${\cal F}(\alpha \wedge e^{i\tilde\omega})$ is also
$G$-equivariant, and its integral is a $G$-invariant distribution on
$\mathfrak{g}^*$.

Suppose now that $\alpha$ is exact: $\alpha = d_{eq} \alpha'$ with
$\alpha' \in \Omega^*_G(M)$. Then
$$
\int_M^{distrib} {\cal F}(\alpha \wedge e^{i\tilde\omega}) =
\int_M^{distrib} \widehat{d_{eq}} \bigl(
{\cal F}(\alpha' \wedge e^{i\tilde\omega}) \bigr) =0
$$
by Lemma \ref{2.11}.
\qed

\begin{df}
Following the lemma, we define the push forward map (\ref{pushforward}) as
$$
{\cal P}: H^*_G(M) \to {\cal C}^{-\infty} (\mathfrak{g}^*)^G,
\qquad
\alpha \mapsto
\int_M^{distrib} {\cal F}(\alpha \wedge e^{i\tilde\omega}).
$$
\end{df}

\begin{ex}
Let $\alpha=1$, then $\int_M^{distrib} {\cal F}(e^{i\tilde\omega})$
is just the Duistermaat-Heckman measure on $\mathfrak{g}^*$.
\end{ex}

When the integral $\int_M \alpha \wedge e^{i\tilde\omega}$ defines a
tempered distribution on $\mathfrak{g}$, we can write
$$
{\cal P}(\alpha) = \int_M^{distrib} {\cal F}(\alpha \wedge e^{i\tilde\omega})
= {\cal F} \biggl( \int_M \alpha \wedge e^{i\tilde\omega} \biggr).
$$
It is the case, for example when $\alpha$ has a compact support on $M$.
One interesting example is the case where $M$ is a closed coadjoint orbit
of a connected semisimple Lie group $G^{ss}$. Let $G$ be a maximal compact
connected subgroup of $G^{ss}$. The action of $G$ on $M$ is Hamiltonian and
one can show that $\int_M \alpha \wedge e^{i\tilde\omega}$ defines a
tempered distribution on $\mathfrak{g}$ (see \cite{DV, Par1} for more details).
However, in general the distribution ${\cal P}(\alpha)$ itself need not be
tempered.

The distribution
${\cal P}(\alpha) = \int_M^{distrib} {\cal F}(\alpha \wedge e^{i\tilde\omega})$
has a useful property of  {\em locality}.
For every $G$-invariant open subset ${\cal U}$ of $\mathfrak{g}^*$,
we denote by $\Phi \mapsto \Phi |_{\cal U}$ the restriction map
${\cal C}^{-\infty} (\mathfrak{g}^*)^G \to {\cal C}^{-\infty} ({\cal U})^G$.
Here we have the following equality in ${\cal C}^{-\infty} ({\cal U})^G$:
$$
\int_M^{distrib} {\cal F}(\alpha \wedge e^{i\tilde\omega}) \biggr|_{\cal U}
=
\int_{\mu^{-1}({\cal U})}^{distrib} {\cal F}(\alpha \wedge e^{i\tilde\omega}).
$$
One can check, as we did in Lemma \ref{2.12}, that the restriction
$\int_M^{distrib} {\cal F}(\alpha \wedge e^{i\tilde\omega}) \bigr|_{\cal U}$
depends only on the cohomology class of $\alpha$ in $\mu^{-1}({\cal U})$.

In his paper \cite{Par2} P.-E. Paradan studies the distribution
$\int_M^{distrib} {\cal F}(\alpha \wedge e^{i\tilde\omega})$ on
$\mathfrak{g}^*$ and gives a localization formula for it at
(the connected components of) the critical points of $\|\mu\|^2$.

\separate

\subsection{Partition of Unity in Equivariant Cohomology}

Let $\lambda \in \Omega^1(M)$ be a $G$-invariant 1-form on $M$.
We denote by $\Phi_{\lambda}: M \to \mathfrak{g}^*$ the $G$-invariant
map defined by
$$
\langle \Phi_{\lambda}(m), X \rangle =_{def} \lambda_m (L_X \bigr|_m),
\qquad m \in M, \: X \in \mathfrak{g}.
$$
We will see in a moment that the equivariant form
$$
(d_{eq} \lambda) (X) = d\lambda - \langle \Phi_{\lambda}, X \rangle
$$
is invertible outside the set $\{ \Phi_{\lambda} = 0 \}$ in the space of
equivariant forms with distributional coefficients.

\begin{lem}
For each $G$-invariant differential form $\eta_0$ on $M$,
vanishing in a neighborhood of $\{ \Phi_{\lambda} = 0 \}$,
the limit of equivariant forms
$$
\lim_{a \to +\infty} \eta_0 \cdot \int_{t=0}^a ie^{-it d_{eq} \lambda} \,dt
$$
exists, hence we can define an equivariant form
$$
\eta_0 \cdot \int_{t=0}^{\infty} ie^{-it d_{eq} \lambda} \,dt
\quad \in \Omega^{-\infty}_G(M).
$$
This form satisfies
\begin{equation}  \label{inverse}
\eta_0 \biggl( \int_{t=0}^{\infty} ie^{-it d_{eq} \lambda} \,dt \biggr)
\wedge d_{eq} \lambda = \eta_0
\qquad \text{in $\Omega^{-\infty}_G(M)$.}
\end{equation}
In particular, the form $d_{eq} \lambda$ is invertible outside the set
$\{ \Phi_{\lambda} = 0 \}$ in the space of equivariant forms with
distributional coefficients.
\end{lem}

\pf
Consider a test form $\phi \in \Omega_c^{top}(\mathfrak{g})$ and a point
$m \in M$, we have:
\begin{multline*}
\biggl\langle \eta_0 \cdot \int_{t=0}^a ie^{-it d_{eq} \lambda(X)} \,dt,
\phi \biggr\rangle_{\mathfrak{g}} \biggr|_m
= \eta_0 \bigr|_m \cdot \int_{t=0}^a \biggl( \int_{\mathfrak{g}}
ie^{-it d_{eq} \lambda_m(X)} \phi \biggr) \,dt  \\
= \eta_0 \bigr|_m \cdot \int_{t=0}^a ie^{-itd\lambda_m}
\biggl( \int_{\mathfrak{g}}
e^{it \langle \Phi_{\lambda}, X \rangle_m} \phi \biggr) \,dt  \\
= \sum_{2j \le \dim M} (-i)^{j-1} (d \lambda_m)^j \eta_0 \bigr|_m \cdot
\int_{t=0}^a \biggl( t^j \cdot \int_{\mathfrak{g}}
e^{it \langle \Phi_{\lambda}, X \rangle_m} \phi \biggr) \,dt.
\end{multline*}
Since $\eta_0$ vanishes in a neighborhood of $\{ \Phi_{\lambda} = 0 \}$
and since the Fourier transform 
$\int_{\mathfrak{g}} e^{it \langle \Phi_{\lambda}, X \rangle_m} \phi$
decays rapidly whenever $\Phi_{\lambda} \ne 0$, the limit of integrals
$$
\lim_{a \to +\infty} 
\eta_0 \bigr|_m \cdot \int_{t=0}^a \biggl( t^j \cdot \int_{\mathfrak{g}}
e^{it \langle \Phi_{\lambda}, X \rangle_m} \phi \biggr) \,dt
$$
exists for each $m \in M$ and defines a differential form on $M$.

Next we observe that, for each $a>0$,
$$
\eta_0 \biggl( \int_{t=0}^a ie^{-it d_{eq} \lambda} \,dt \biggr)
\wedge d_{eq} \lambda = \eta_0 (1 - e^{-ia d_{eq} \lambda}).
$$
Since
$$
\lim_{a \to + \infty} \langle \eta_0 \cdot e^{-ia d_{eq} \lambda},
\phi \rangle_{\mathfrak{g}} \bigr|_m
= \lim_{a \to + \infty}
\sum_{2j \le \dim M} (-i)^j (d \lambda_m)^j \eta_0 \bigr|_m \cdot
a^j \cdot \int_{\mathfrak{g}}
e^{ia \langle \Phi_{\lambda}, X \rangle_m} \phi =0,
$$
$\lim_{a \to + \infty} \eta_0 \cdot e^{-ia d_{eq} \lambda} =0$
as a distribution on $\mathfrak{g}$, and formula (\ref{inverse}) follows.

Finally, taking $\eta_0$ a $G$-invariant function equal 1 outside of a
small neighborhood of $\{ \Phi_{\lambda} = 0 \}$, we see that
$d_{eq} \lambda$ is invertible outside $\{ \Phi_{\lambda} = 0 \}$ in
$\Omega^{-\infty}_G(M)$.
\qed

Let $\eta \in {\cal C}^{\infty}(M)^G$ equal 1 in a neighborhood of
$\{ \Phi_{\lambda} = 0 \}$. Then the form $d\eta$ vanishes in a neighborhood
of $\{ \Phi_{\lambda} = 0 \}$, and we can define an equivariant form
$$
\operatorname{P}_{\lambda} = \eta +
d\eta \cdot \biggl( \int_{t=0}^{\infty} ie^{-it d_{eq} \lambda} \,dt \biggr)
\wedge \lambda
\quad \in \Omega^{-\infty}_G(M).
$$
From (\ref{inverse}) we get an identity which P.-E.~Paradan calls
``partition of unity.''

\begin{prop} [Proposition 2.1 in \cite{Par2}]
The equivariant form $\operatorname{P}_{\lambda}$ is closed
($d_{eq} \operatorname{P}_{\lambda}=0$), and we have the identity
$$
1_M = \operatorname{P}_{\lambda} + d_{eq} \delta,
$$
where $1_M$ denotes the constant function 1 on $M$ and
$$
\delta =
(1-\eta) \biggl( \int_{t=0}^{\infty} ie^{-it d_{eq} \lambda} \,dt \biggr)
\wedge \lambda
\quad \in \Omega^{-\infty}_G(M).
$$
\end{prop}

We will decompose the equivariant form $\operatorname{P}_{\lambda}$ into a
sum of forms, each summand corresponding to a component of
$\{ \Phi_{\lambda} = 0 \}$.

\begin{df}
A subset $C \subset M$ is called a {\em component} of
$\{ \Phi_{\lambda} = 0 \}$ if there exists a $G$-invariant open neighborhood
$U$ of $C$ in $M$ such that
$$
\overline{U} \cap \{ \Phi_{\lambda} = 0 \} =C.
$$
This condition implies
$$
\partial U \cap \{ \Phi_{\lambda} = 0 \} = \varnothing.
$$

A $G$-invariant open set $U$ of $M$ which satisfies the last condition
will be called {\em good} for the 1-form $\lambda$; it intersects
$\{ \Phi_{\lambda} = 0 \}$ in the interior of $\overline{U}$.

Let $U subset M$ be a good open set for the 1-form $\lambda$.
The intersection $\overline{U} \cap \{ \Phi_{\lambda} = 0 \}$ (possibly empty!)
is called a {\em component} of $\{ \Phi_{\lambda} = 0 \}$.

Let $U$ be a good open set for the 1-form $\lambda$, and
$C = \overline{U} \cap \{ \Phi_{\lambda} = 0 \}$ be the corresponding
component of $\{ \Phi_{\lambda} = 0 \}$.
We denote by $\operatorname{P}_{\lambda}^U$ (or $\operatorname{P}_{\lambda}^C$)
the equivariant form
$$
\operatorname{P}_{\lambda}^U = \eta_U +
d\eta_U \cdot \biggl( \int_{t=0}^{\infty} ie^{-it d_{eq} \lambda} \,dt \biggr)
\wedge \lambda
\quad \in \Omega^{-\infty}_G(M),
$$
where $\eta_U$ is a $G$-invariant real-valued function on $M$ equal 1
in a neighborhood of $C$ and with support in $U$.
\end{df}

Observe that $\operatorname{P}_{\lambda}^U$ belongs to
$\Omega^{-\infty}_G(U)$. Hence if $\overline U$ is compact, then the form
$\operatorname{P}_{\lambda}^U$ has compact support on $M$.
Also the equivariant cohomology class of $\operatorname{P}_{\lambda}^U$ does
not depend on the choice of the function $\eta_U$.

Consider an open covering $\bigcup_j U_j$ of $\{ \Phi_{\lambda} = 0 \}$,
where the open sets $U_j$ are $G$-invariant. Suppose furthermore that
$U_i \cap U_j = \varnothing$ if $i \ne j$. Then the open sets $U_j$ are good
for the 1-form $\lambda$ and
$$
\operatorname{P}_{\lambda} = \sum_j \operatorname{P}_{\lambda}^{U_j}.
$$
This sum is well defined, even if it is infinite, because the equivariant forms
$\operatorname{P}_{\lambda}^{U_j}$ have disjoint supports.
This decomposition enables us to study the equivariant form
$\operatorname{P}_{\lambda}$ in the neighborhood of each component of
$\{ \Phi_{\lambda} = 0 \}$.

Next we consider two $G$-invariant 1-forms $\lambda_0$, $\lambda_1$ on $M$,
and a $G$-invariant open set $U$ of $M$.

\begin{prop} [Proposition 2.6 in \cite{Par2}]
Suppose there exists a smooth map $f: M \to \mathfrak{g}$ and a real
$\rho >0$ such that the functions
$$
l_0 = \langle \Phi_{\lambda_0}, f \rangle
\qquad \text{and} \qquad
l_1 = \langle \Phi_{\lambda_1}, f \rangle
$$
are bounded from below by $\rho$ on $\partial U$. Then the open set $U$
is good for $\lambda_0$ and $\lambda_1$, and
$$
\operatorname{P}_{\lambda_0}^U = \operatorname{P}_{\lambda_1}^U
\qquad \text{in $H^{-\infty}_G(M)$.}
$$
\end{prop}

\separate

\subsection{Localization of the Push Forward ${\cal P}$ at
$\operatorname{Cr} (\|\mu\|^2)$}

Choose a $G$-invariant Euclidean norm $\|\cdot\|$ on $\mathfrak{g}^*$,
and let ${\cal H}$ be the Hamiltonian vector field on $M$ associated to
the function $\frac 12 \|\mu\|^2$. (This means that ${\cal H}$ is the unique
vector field on $M$ satisfying
$\frac 12 d(\|\mu\|^2) = \iota({\cal H}) \omega$.)
Fix a $G$-invariant Riemannian metric $(\cdot,\cdot)_M$ on $M$.
The Witten 1-form is a $G$-invariant form on $M$ defined by
$$
\lambda =_{def} ({\cal H}, \cdot )_M.
$$

\begin{lem}
$\{ \Phi_{\lambda} = 0 \}$ is precisely the set of critical points
of $\|\mu\|^2$:
$$
\{ \Phi_{\lambda} = 0 \} = \operatorname{Cr} (\|\mu\|^2).
$$
\end{lem}

\pf
Let $( \cdot, \cdot)$ be a $G$-invariant Euclidean metric on $\mathfrak{g}^*$
such that $\|\mu\|^2 = (\mu,\mu)$.
Define a map $\mu': M \to \mathfrak{g}$ so that
$$
\langle \xi, \mu'(m) \rangle = (\xi, \mu(m)),
\qquad \forall m \in M, \: \xi \in \mathfrak{g}^*.
$$
Then
$$
\iota({\cal H}) \omega = \frac 12 d(\|\mu\|^2) = (d\mu, \mu)
= \langle d\mu, \mu' \rangle = \iota(L_{\mu'})\omega.
$$
Hence ${\cal H} = L_{\mu'}$ which shows that the vector field ${\cal H}$
is always tangent to the $G$-orbits in $M$.
Therefore,
$$
\{ \Phi_{\lambda} = 0 \} = \{ {\cal H} =0 \} = \operatorname{Cr} (\|\mu\|^2).
$$
\qed

Associated to the Witten 1-form $\lambda$ we have a partition of unity
in equivariant cohomology: $1_M = \operatorname{P}_{\lambda} + d_{eq} \delta$.
For every equivariantly closed form $\alpha \in \Omega^*_G(M)$, the partition
of unity decomposes the equivariant form $\alpha \wedge e^{i\tilde\omega}$ as
$$
\alpha \wedge e^{i\tilde\omega} =
\operatorname{P}_{\lambda} \wedge \alpha \wedge e^{i\tilde\omega}
+ d_{eq} (\delta \wedge \alpha \wedge e^{i\tilde\omega}).
$$

\begin{prop} [Proposition 3.4 in \cite{Par2}]
Let $\alpha(X) \in \Omega^*_G(M)$ be an equivariant form depending
polynomially on $X \in \mathfrak{g}$. The equivariant forms
$\operatorname{P}_{\lambda} \wedge \alpha \wedge e^{i\tilde\omega}$ and
$\delta \wedge \alpha \wedge e^{i\tilde\omega}$ are tempered and their
Fourier transforms
${\cal F} (\operatorname{P}_{\lambda} \wedge \alpha \wedge
e^{i\tilde\omega})$ and
${\cal F}(\delta \wedge \alpha \wedge e^{i\tilde\omega})$ have compact supports
in $\mathfrak{g}^*$-mean on $M$.

Moreover, let $U$ be a good open set for the 1-form $\lambda$.
Then the same result holds for $\operatorname{P}_{\lambda}^U$ in place
of $\operatorname{P}_{\lambda}$.
\end{prop}

Recall that
${\cal F} \circ d_{eq} (\delta \wedge \alpha \wedge e^{i\tilde\omega})
=\widehat{d_{eq}} \circ {\cal F}
(\delta \wedge \alpha \wedge e^{i\tilde\omega})$.
Hence from Lemma \ref{2.11} we get

\begin{cor}
Let $\alpha(X) \in \Omega^*_G(M)$ be an equivariantly closed form depending
polynomially on $X \in \mathfrak{g}$. We have the following equality of
$G$-invariant distributions on $\mathfrak{g}^*$:
$$
{\cal P}(\alpha) = \int_M^{distrib} {\cal F}(\alpha \wedge e^{i\tilde\omega})
= \int_M^{distrib} {\cal F}
(\operatorname{P}_{\lambda} \wedge \alpha \wedge e^{i\tilde\omega}).
$$
\end{cor}

Recall that the equivariant form $\operatorname{P}_{\lambda}$ is supported
in a $G$-invariant neighborhood of
$\{ \Phi_{\lambda} = 0 \} = \operatorname{Cr} (\|\mu\|^2)$.
This neighborhood can be made arbitrary small and $G$-equivariantly identified
with the normal bundle of $\operatorname{Cr} (\|\mu\|^2)$ in $M$.
Thus ${\cal P}(\alpha)$ localizes at the components of
$\operatorname{Cr} (\|\mu\|^2)$.
In his paper \cite{Par2} P.-E.~Paradan computes the contributions of
the components of $\operatorname{Cr} (\|\mu\|^2)$ to the localization formula
for ${\cal P}(\alpha)$ and describes their qualitative properties.

\separate

\separate

\noindent
{\em E-mail address:} {matvei.libine@yale.edu}

\noindent
{\em Department of Mathematics, Yale University,
P.O. Box 208283, New Haven, CT 06520-8283}

\end{document}